\documentclass[10pt]{amsart}

\usepackage{graphicx}
\usepackage{times}
\usepackage[colorlinks=true,linkcolor=blue]{hyperref}%

\usepackage{xcolor}
\usepackage{hyperref}
\hypersetup{
    linktoc=all,
    colorlinks=true
}
\usepackage{stmaryrd}
\usepackage{pdfrender,xcolor}

\newtheorem{theorem}{Theorem}[section]
\newtheorem{lemma}[theorem]{Lemma}
\newtheorem{corollary}[theorem]{Corollary}
\newtheorem{prop}[theorem]{Proposition}

\theoremstyle{definition}
\newtheorem{definition}[theorem]{Definition}

\theoremstyle{remark}
\newtheorem{remark}[theorem]{Remark}
\numberwithin{equation}{section}

\usepackage[top=0.75in, bottom=0.75in, left=0.75in, right=0.75in]{geometry}

\usepackage[scr]{rsfso}
\usepackage[english]{babel}
\usepackage{fancyhdr}
\usepackage{amsmath}
\usepackage{amsthm}
\usepackage{amssymb}
\usepackage{eucal}
\usepackage{comment}
\usepackage{enumitem}
\setlist{leftmargin=*}
\usepackage[integrals]{wasysym}

\newcommand\nc{\newcommand}
\nc{\on}{\operatorname}
\nc{\E}{\mathbf{E}}
\nc{\R}{\mathbb R}
\nc{\C}{\mathbb C}
\nc{\Q}{\mathbb Q}
\nc{\Z}{\mathbb Z}
\nc{\N}{\mathbb N}
\nc{\F}{\mathbb F}
\nc{\wt}{\widetilde}
\nc{\ol}{\overline}
\nc{\short}[3]{0 \longrightarrow #1 \longrightarrow #2 \longrightarrow #3 \longrightarrow 0}
\nc{\pd}[2]{\frac{\partial #1}{\partial #2}}
\nc{\rnc}{\renewcommand}
\nc{\e}{\varepsilon}
\nc{\DMO}{\DeclareMathOperator}
\nc{\grad}{\nabla}
\nc{\Exp}{\mathbf{Exp}}
\nc{\fsp}{\fontdimen2\font=2.17pt}

\rnc{\t}{\mathrm{t}}
\nc{\s}{\mathrm{s}}
\rnc{\r}{\mathrm{r}}
\nc{\x}{\mathrm{x}}
\nc{\y}{\mathrm{y}}
\nc{\z}{\mathrm{z}}
\nc{\w}{\mathrm{w}}
\rnc{\and}{\quad\mathrm{and}\quad}

\rnc{\leq}{\leqslant}
\rnc{\geq}{\geqslant}
\rnc{\d}{\mathrm{d}}
\rnc{\O}{\mathrm{O}}
\rnc{\exp}{\mathbf{Exp}}
\newenvironment{nouppercase}{%
  \renewcommand{\uppercasenonmath}[1]{}}{}
\def\XXint#1#2#3{{\setbox0=\hbox{$#1{#2#3}{\int}$ }
\vcenter{\hbox{$#2#3$ }}\kern-.6\wd0}}

\title{\fontdimen2\font=1.7pt\Large Non-Stationary KPZ equation from ASEP with slow bonds \vspace{-0.3cm}}
\author{ \large Kevin Yang}
\usepackage{setspace}
\begin{document}
%\pdfrender{StrokeColor=gray,TextRenderingMode=2,LineWidth=0.01pt}
\setstretch{0.99}
\fsp
\raggedbottom
\begin{nouppercase}
\maketitle
\end{nouppercase}
\vspace{-15pt}
\begin{center}
\today
\end{center}
\begin{abstract}
\fsp We prove the height functions for a class of non-integrable and non-stationary particle systems converge to the KPZ equation, thereby making progress on the universality of the KPZ equation. The models herein are ASEP \cite{BG} with a mesoscopic family of slow bonds, thus we partially extend \cite{FGS} to non-stationary models and add to the almost empty set of non-integrable, non-stationary interacting particle systems for which universality is established. To do this, we develop further the strategy of \cite{Y,Y20} introduce a method to establish a novel principle that builds upon the classical hydrodynamic limits of \cite{GPV} and that we call \emph{local hydrodynamics}.
\end{abstract}

{\hypersetup{linkcolor=black}
\setcounter{tocdepth}{1}
\tableofcontents}

%%%
\section{Introduction}
%%%
The \emph{Kardar-Parisi-Zhang} (KPZ) equation is a stochastic PDE (SPDE) that was derived in physics literature \cite{KPZ} as a \emph{canonical} and \emph{universal} model for interface fluctuations. Examples of interface fluctuations that are supposed to be described by the KPZ equation include flame propagation, bacterial colony and tumor growth, crack formation, and many others \cite{KPZ}. We now write the KPZ equation formally as the following SPDE that we pose on the $1+1$-dimensional space-time set $\R_{\geq0}\times\R$, in which $\xi$ is the space-time Gaussian white noise with covariance kernel $\E\xi_{\t,\x}\xi_{\s,\y}=\delta_{\t=\s}\delta_{\x=\y}$ to be interpreted in a generalized function sense:
\begin{align}
\partial_{\t}\mathbf{h}^{\infty} \ = \ 2^{-1}\Delta\mathbf{h}^{\infty} - 2^{-1}|\grad\mathbf{h}^{\infty}|^{2} + \xi. \label{eq:kpz}
\end{align}
We briefly note that the superscript $\infty$ is meant to suggest our perspective on the KPZ equation and its solution as a scaling limit of some kind, specifically in the context of interacting particle systems as we soon discuss. Although universality of the KPZ equation for interface fluctuations of many kinds is formally justified in \cite{KPZ} by means of a non-rigorous perturbative renormalization group, this leads to an explanation of universality that produces incorrect numerics \cite{HQ}. This is because the KPZ equation is a \emph{singular} stochastic PDE in that if interpreted in any classical PDE sense, it is ill-posed. This ill-posedness issue was resolved by introducing a physically meaningful intrinsic solution to \eqref{eq:kpz} in work of Hairer \cite{Hai13}, which eventually led to the theory of regularity structures \cite{Hai14} and was used to rigorously prove universality of the KPZ equation for interface models defined by SPDE generalizations of \eqref{eq:kpz}. The main purpose of this article, however, is to prove universality for what is known as a \emph{non-integrable} and \emph{non-stationary} family of interacting particle systems distinguished by \emph{slow bonds}. At a technical level, we develop a novel and general method for \emph{local} hydrodynamics, which are much more robust than the global approach to establishing hydrodynamic limits, which has been intensively employed continually throughout the past several decades, beginning with \cite{GPV}. This method of local hydrodynamics adds onto previous work \cite{Y} of the author, which focuses on local equilibrium methods for KPZ universality, with estimates of hydrodynamic limits of general observables in general particle systems based upon \emph{analytic regularity estimates} for a microscopic version of the height function, namely the height function associated to the interacting particle system itself. In a somewhat more analytic direction, this paper additionally implements the breakthrough Nash-ian approach to heat kernel estimates for analyzing heat equations in the context of interacting particle systems when more classical methods, which are employed in earlier literature \cite{BG,CST,CT,DT,KL}, are inapplicable. This Nash-ian approach, to the author's knowledge, had only been used in this context in \cite{Y20} for studying heat kernels with boundary conditions. Here, we further develop this Nash-ian approach to rough and degenerate heat operators that ultimately come from slow bonds.
%%%
\subsection{Background}
%%%
We first provide a historical review of universality of the KPZ equation for interacting particle system height functions before we present the details of the main result of this paper and the key innovations herein. This is because our statement of the main result will depend on previous work. To start, in the pioneering work of \cite{BG}, the authors introduced a solution to \eqref{eq:kpz}, which requires much less technology than the theory of regularity structures, that is known as the \emph{Cole-Hopf solution}. It is defined via $\mathbf{h}^{\infty}=-\log\mathbf{Z}^{\infty}$, where $\mathbf{Z}^{\infty}$ is the solution to the following \emph{multiplicative noise stochastic heat equation} (SHE); the solutions of the following SHE are easily shown to be continuous in space-time and positive, if given positive initial data, with probability 1 \cite{Mu}, thereby making the Cole-Hopf solution well-defined; we note \eqref{eq:she} is interpreted via the Duhamel mild form:
\begin{align}
\partial_{\t}\mathbf{Z}^{\infty} \ = \ 2^{-1}\Delta\mathbf{Z}^{\infty} - \mathbf{Z}^{\infty}\xi. \label{eq:she}
\end{align}
Although a formal application of the Ito formula with respect to the space-time white noise relates $\mathbf{Z}^{\infty}$ and $\mathbf{h}^{\infty}$ in exactly the way the Cole-Hopf solution suggests, such a calculation is purely formal and work is needed to show the Cole-Hopf solution is actually meaningful; indeed, one cannot simply avoid the singular issues of KPZ with a formal definition without justification. The work \cite{BG}, fortunately, justifies the Cole-Hopf solution as physically relevant by showing the height function associated to an interacting particle system that is known as \emph{asymmetric simple exclusion process} (ASEP) converges to the Cole-Hopf solution of \eqref{eq:kpz} after appropriate rescaling and renormalization. This began an massive effort of finding more interacting particle systems whose height functions similarly converge to the Cole-Hopf solution of KPZ by applying the method of \cite{BG}; see \cite{CST,CT}, for example. However, such results are not \emph{universality} results, as they depend heavily on important algebraic features of the interacting particle systems that are by no means perturbative or general; such features reflect a self-duality for the interacting particle system. To address this, therefore open, universality problem for the KPZ equation, another approach to defining solutions of \eqref{eq:kpz} was introduced in \cite{GJ15}. This is the energy solution theory, and its perspective and technology allowed \cite{GJ15,GJS15}, among many other papers, to confirm the universality of \eqref{eq:kpz} for height functions of many interacting particle systems. However, such an approach requires said interacting particle systems to be globally stationary with invariant measures that are sufficiently explicit. Under this assumption, microscopic correlations in the particle system are basically explicit, which is a key ingredient for energy solution theory. Actually, by using the energy solution theory, the authors of \cite{FGS} were able to study perturbations of the ASEP model from \cite{BG} with a slow bond. One of the main results, if not the primary result, of this paper \cite{FGS} was in establishing KPZ equation scaling limits for the height function, or some basically equivalent version of this statement, with a main interest being that introducing a slow bond to ASEP destroys its integrability. However, we reiterate that \cite{FGS} requires the slow bond model to be globally stationary. In this paper, we will remove completely this global stationary assumption for a family of these slow bond models in \cite{FGS}. We will actually be able to generalize \cite{FGS} not just to non-stationary models but also to other perturbations of ASEP with a mesoscopic number of slow bonds instead of a single slow bond, and we will not require any explicit knowledge of the invariant measure, though for cleaner presentation we will not address this last point in detail but rather discuss it briefly when relevant. This concludes our preliminary introduction and background discussion. When we have introduced the main constructions and results, we will return to an additional discussion on background and history of KPZ universality, non-stationary models, slow bond models, and the techniques of this paper, which we recall amount to employing the local equilibrium method developed in \cite{Y} to establish a \emph{local hydrodynamic} principle. Actually, we remark that slow bond models, in particular, admit an extremely rich and interesting history/literature in ways that are related to KPZ statistics and that have actually generated controversy between physical and mathematical communities/predictions, which have only recently been started to be resolved \cite{BSaS,BSiS}, providing additional motivation for this paper. We explore this more shortly. We reemphasize that the slow bond models considered here are not only non-stationary but also not integrable. Therefore, the work of this paper adds to the incredibly modest set of such interacting particle systems \cite{DT} for which universality of the KPZ equation is proven and adds to the ongoing development of a general method \cite{Y,Y20} to prove universality of the KPZ equation in general.
%%%
\subsection{The Model}
%%%
Throughout the rest of this article, we will employ $N\in\Z_{\geq0}$ as a scaling parameter that will eventually be taken to infinity to recover relevant scaling limits, such as KPZ scaling limits for the height function of the interacting particle system of interest in this paper, both of which will be introduced and defined precisely in this section. To introduce the interacting particle system, we view it as a Markov process on a state space encoding particle configurations and provide its infinitesimal generator.
%%%
\begin{definition}\label{definition:model}
\fsp Given any subset $\mathbb{I}\subseteq\Z$, let us define $\Omega_{\mathbb{I}}=\{\pm1\}^{\mathbb{I}}$, and let us denote elements by $\eta$ with interpretation that $\eta_{\x}=1$ if there is a particle at $\x$ in the $\eta$-configuration and otherwise $\eta_{\x}=-1$; note these $\eta$-variables encode fluctuations in particle density. For convenience, we will also set $\Omega=\Omega_{\Z}$ to be the ``total" state space. Lastly, we will not refer to exactly which $\Omega_{\mathbb{I}}$ space any given configuration $\eta$ belongs to; equivalently, if we take $\eta\in\Omega$, then we realize $\eta\in\Omega_{\mathbb{I}}$ for any $\mathbb{I}\subseteq\Z$ by forgetting $\eta_{\x}$ for all $\x\not\in\mathbb{I}$.

Let us now fix a parameter $\beta_{\star}\in[0,2^{-1})$ to be interpreted as the strength of the slow bond in the forthcoming interacting particle system generator. We also define $\mathscr{L}_{\x,\y}$ as the generator for a speed-$1$ symmetric simple exclusion process on the bond $\{\x,\y\}\subseteq\Z$, and we additionally define $\mathscr{L}_{\x}=\mathscr{L}_{\x,\x+1}$ for convenience. Given any local function $\mathfrak{f}:\Omega\to\R$, where ``local" refers to depending only on a finite number of $\eta$-values per configuration $\eta\in\Omega$, we define $\mathscr{L}^{N}\mathfrak{f}=\mathscr{L}^{N,\mathrm{S}}\mathfrak{f}+\mathscr{L}^{N,\mathrm{A}}{\mathfrak{f}}$, where $\mathscr{L}^{N,\mathrm{S}}$ is defined by
\begin{align}
\mathscr{L}^{N,\mathrm{S}}\mathfrak{f}(\eta) \ = \ 2^{-1}N^{2}{\sum}_{\x\neq0}\mathscr{L}_{\x}\mathfrak{f}(\eta) + 2^{-1}N^{2-\beta_{\star}}\mathscr{L}_{0}\mathfrak{f}(\eta),
\end{align}
and $\mathscr{L}^{N,\mathrm{A}}$ is the following asymmetric part of the interacting particle system with no slow bond effect:
\begin{align}
\mathscr{L}^{N,\mathrm{A}}\mathfrak{f}(\eta) \ = \ 2^{-1}N^{3/2}{\sum}_{\x\in\Z}\mathbf{1}_{\eta_{\x}=-1}\mathbf{1}_{\eta_{\x+1}=1}\mathscr{L}_{\x}\mathfrak{f}(\eta) - 2^{-1}N^{3/2}{\sum}_{\x\in\Z}\mathbf{1}_{\eta_{\x}=1}\mathbf{1}_{\eta_{\x+1}=-1}\mathscr{L}_{\x}\mathfrak{f}(\eta).
\end{align}
To conclude, we define the \emph{support} of a function $\mathfrak{f}:\Omega\to\R$ as the subset $\mathbb{I}(\mathfrak{f})\subseteq\Z$ such that $\mathfrak{f}$ depends only on $\eta_{\x}$ for $\x\in\mathbb{I}$. Lastly, for any $\t\geq0$, we define $\eta_{\t}\in\Omega$ as the time-$\t$ configuration under the $\mathscr{L}^{N}$ Markov process with initial configuration $\eta_{0}$.
\end{definition}
%%%
%%%
\begin{remark}\label{remark:model}
\fsp We could have included a mesoscopic number of slow bonds instead of just one slow bond and this paper would only change in technical power-counting in the scaling parameter $N$. Indeed, estimates in this paper are explicit in powers of the scaling parameter whenever we deal with slow bonds, and multiplying negative powers of $N$ by sufficiently small but still positive powers of $N$ does not change asymptotic vanishing in the large-$N$ limit. However, we focus on one slow bond to make the article clearer. Lastly, we note the condition $\beta_{\star}\in[0,2^{-1})$, like in \cite{FGS}, is to guarantee the generator is an honest infinitesimal generator without negative rates. We also note that if $\beta_{\star}=0$, then we have the ASEP model in \cite{BG} without any slow bond effects. In this paper, we will eventually take $\beta_{\star}$ very small but still uniformly positive and independent of $N$. As we discuss after the statement of Theorem \ref{theorem:kpz}, which is the main result of this paper, even for such small $\beta_{\star}$, our work is not perturbative in $\beta_{\star}$.
\end{remark}
%%%
We now introduce the corresponding height function for the slow bond model of Definition \ref{definition:model}. This will serve as a discretization of the KPZ equation, as universality of the KPZ equation and Theorem \ref{theorem:kpz} both tell us. We will also introduce a microscopic version of the SHE \eqref{eq:she}, as the SHE is how we will make sense of the KPZ equation in this paper, similar to \cite{BG}.
%%%
\begin{definition}\label{definition:hf}
\fsp For any $\t\geq0$, we define $\mathbf{h}^{N}_{\t,0}$ as $2N^{-1/2}$ times the net flux of the bond $\{0,1\}$ up to time $\t$, which is the number of particles that went from $1$ to $0$ minus the number of particles that have gone from $0$ to $1$, both until time $\t$. Given $\x\in\Z$, we define the microscopic rescaled \emph{height function}, which is a function on $\R_{\geq0}\times\Z$, as the following ``integration" of $\eta$-variables; see \cite{BG}:
\begin{align}
\mathbf{h}^{N}_{\t,\x} \ = \ \mathbf{h}^{N}_{\t,0} + \mathbf{1}_{\x\geq1}N^{-1/2}{\sum}_{1\leq\y\leq\x}\eta_{\t,\y} - \mathbf{1}_{\x\leq-1}N^{-1/2}{\sum}_{\x\leq\y\leq0}\eta_{\t,\y}.
\end{align}
We also define the \emph{Gartner transform} $\mathbf{G}^{N}_{\t,\x}=\exp(-\mathbf{h}^{N}_{\t,\x}+\mathrm{R}_{N}\t)$, where $\mathrm{R}_{N}=2^{-1}N-(24)^{-1}$. We will also linearly interpolate values of $\mathbf{G}^{N}$ on $\R_{\geq0}\times\Z$ to obtain a continuous and spatially-piecewise-linear function on $\R_{\geq0}\times\R$, just for convenience. Let us note that we used $\mathbf{G}^{N}$ for the microscopic version of SHE instead of $\mathbf{Z}^{N}$, because fairly early on in Section \ref{section:mshe} we will modify $\mathbf{G}^{N}$ slightly into a different process, denoted by $\mathbf{Z}^{N}$, towards which basically the entire remainder of the article is dedicated. 
\end{definition}
%%%
%%%
\subsection{Main Results}
%%%
We will now present the main KPZ universality result of the paper. To start, we will need to introduce a class of allowable initial data for the interacting particle system, which is given in Definition \ref{definition:ns} below. Roughly speaking, initial data satisfying constraints for the ``near-stationary" initial data in Definition \ref{definition:ns} translate into analytically reasonable initial data for the proposed SHE limit of $\mathbf{G}^{N}$. First, however, we introduce a pair of rescaling operators used throughout this section.
%%%
\begin{definition}\label{definition:gammas}
\fsp Provided any function $\phi:\R_{\geq0}\times\Z\to\R$, we define $\Gamma^{N}\phi_{\t,\x}=\phi_{\t,N\x}$. Note that time is not scaled, as we have built the relevant time-scaling into the particle system already. Let us also define $\Gamma^{N,\mathbf{X}}\phi_{\x}=\phi_{0,N\x}$ for rescaling only the initial data.
\end{definition}
%%%
%%%
\begin{definition}\label{definition:ns}
\fsp A probability measure $\mu$ on the state space $\Omega$ is \emph{near-stationary} if for any $p\geq1$, for any $\alpha\in(0,2^{-1})$, and for some $\kappa\geq0$, we have the following moment estimates in which $\mathbf{G}^{N}_{0,\cdot}$, viewed as a function on $\Omega$, is distributed via $\mu$:
\begin{align}
\exp(-2p\kappa N^{-1}|\x|)\E|\mathbf{G}^{N}_{0,\x}|^{2p} \ \lesssim_{p} \ 1 \and \exp(-2p\kappa N^{-1}(|\x|+|\y|))\E|\mathbf{G}^{N}_{0,\x}-\mathbf{G}^{N}_{0,\y}|^{2p} { \ \lesssim_{\alpha,p} \ N^{-2p\alpha}|\x-\y|^{2p\alpha}.}
\end{align}
We also require that $\Gamma^{N,\mathbf{X}}\mathbf{G}^{N}$ converges locally uniformly and in probability to initial data denoted by $\Gamma^{\infty,\mathbf{X}}\mathbf{G}^{\infty}\in\mathscr{C}(\R)$. 
\end{definition}
%%%
Lastly, before we can state the main result of this paper, we must introduce the Skorokhod space $\mathscr{D}(\R_{\geq0},\mathscr{B})$ to be the space of cadlag paths valued in the Frechet space $\mathscr{B}$, which we will take to be $\mathscr{C}(\R)$ in Theorem \ref{theorem:kpz}. We equip this space with the Skorokhod topology that is detailed in Chapter 14 of \cite{Bil}; although \cite{Bil} only treats $\mathscr{B}=\R$, the same story holds for any Frechet space $\mathscr{B}$ if we replace the metric on $\R$ with the metric on $\mathscr{B}$ in general.
%%%
\begin{theorem}\label{theorem:kpz}
\fsp Assume the interacting particle system starts with near-stationary initial data {with $\kappa=0$}. We can pick $\bar{\beta}_{\star}>0$ independently of $N$ such that for any $\beta_{\star}\leq\bar{\beta}_{\star}$, the rescaled process $\Gamma^{N}\mathbf{G}^{N}$ converges to the solution of $\mathrm{SHE}$ in the large-$N$ limit with the initial data $\Gamma^{\infty,\mathbf{X}}\mathbf{G}^{\infty}$ from \emph{Definition \ref{definition:ns}}. This convergence is in the Skorokhod space $\mathscr{D}(\R_{\geq0},\mathscr{C}(\R))$.
\end{theorem}
%%%
%%%
\begin{remark}\label{remark:ns}
\fsp {In Theorem \ref{theorem:kpz}, we assumed $\kappa=0$ for our near stationary data; see Definition \ref{definition:ns}. (From this assumption, we get moment bounds uniform in space for our initial data of $\mathbf{G}^{N}$.) This assumption is just for convenience. Indeed, to remove it, one can just include sub-exponential weights to every space-norm of $\mathbf{G}^{N}$ and of any solution to any microscopic SPDE. This is exactly what is done in Section 3 of \cite{DT}, for example. All it does is formally add notation. (The math is otherwise identical.)}
\end{remark}
%%%
%%%
\begin{remark}\label{remark:kpzkpz}
\fsp Another initial data of interest for the KPZ equation is the so-called narrow-wedge initial data, which is given by a Dirac point mass at the origin and uncovers connections between KPZ and integrable systems; see \cite{ACQ,DT}. In \cite{ACQ,DT}, the interacting particle systems considered therein were shown to have height function converging to KPZ with narrow-wedge initial data as well, not just near-stationary data. The analogous result for slow bond models would be the following. Assume the interacting particle system starts with narrow-wedge initial data. We can, again, pick $\bar{\beta}_{\star}$ positive such that for any $\beta_{\star}\leq\bar{\beta}_{\star}$, the further rescaled process $2^{-1}N^{1/2}\Gamma^{N}\mathbf{G}^{N}$ converges to the solution of $\mathrm{SHE}$ in the large-$N$ limit with initial data the Dirac delta function supported at the origin $\{0\}\subseteq\R$. This convergence is in the Skorokhod space $\mathscr{D}(\R_{>0},\mathscr{C}(\R))$, which we will interpret to be the topological limit of $\mathscr{D}(\R_{\geq\e},\mathscr{C}(\R))$ for every $\e>0$, so that convergence is, by definition, convergence in each of these $\e$-spaces. To be precise, the narrow-wedge data is the measure concentrated on the deterministic configuration $\eta^{\mathrm{nw}}\in\Omega$ such that $\eta^{\mathrm{nw}}_{\x}=\mathrm{sgn}(\x)$ for all $\x\neq0$ and $\eta_{0}^{\mathrm{nw}}=-1$. The proof of the result that we recorded in this remark for narrow-wedge data is not difficult given the work of this paper, though it requires a lengthy argument. In a nutshell, we will require analysis of a discrete-type SPDE for $\mathbf{G}^{N}$. For regular near-stationary data, our analysis of the slow-bond input in said $\mathbf{G}^{N}$-equation amounts to estimates in a type of space-time uniform norm, under which $\mathbf{G}^{N}$ is well-behaved. For narrow-wedge data, we expect $\mathbf{G}^{N}$ to exhibit a heat kernel singularity at short times, so our analysis of the $\mathbf{G}^{N}$-equation would be with respect to a time-weighted space-time norm to control singular behavior of $\mathbf{G}^{N}$ for short times. Ultimately, this does not change our work beyond a handful of cosmetic changes. More technically, our analysis of the slow bond input in the $\mathbf{G}^{N}$ already addresses {heat}-kernel-type singularities, so the change of weights does not affect the success of the method that we develop herein. However, we refrain from adding details to avoid significantly lengthening this paper.
\end{remark}
%%%
%%%
\subsection{Additional Commentary}
%%%
The constraint of $\beta_{\star}\leq\bar{\beta}_{\star}$ for some positive $\bar{\beta}_{\star}$ is simply the statement that we must choose $\beta_{\star}$ sufficiently small, but we can still choose it to be uniformly positive and independent of $N$. We will discuss shortly in this subsection why this does not make Theorem \ref{theorem:kpz} a perturbative result, so that new ideas beyond the current literature \cite{BG,DT,Y,Y20} are needed to prove Theorem \ref{theorem:kpz}. In particular, the reader is invited to always think of $\beta_{\star}$ to be sufficiently small so that certain inequalities are true, such as $999\beta_{\star}\leq999^{-1}$. We definitely do not obtain an optimal value of $\bar{\beta}_{\star}$, but the methods, as written, do not let us access $\bar{\beta}_{\star}\approx2^{-1}$ that would be optimal, so we do not worry too much about the optimal $\bar{\beta}_{\star}$ obtained via the methods of this paper.

Next, as noted in \cite{ACQ,DT}, under a certain double limit, Remark \ref{remark:kpzkpz} shows that the height function $\mathbf{h}^{N}$ of the non-integrable and non-stationary interacting particle systems at hand converge under appropriate scaling to the so-called KPZ fixed point; we refer the reader to \cite{QS}. In a somewhat different direction, unlike \cite{FGS} our analysis towards proving Theorem \ref{theorem:kpz} will not require such a refined understanding of the invariant measures for slow bond models herein. This is one main advantage of our analysis, and we hope that it provides a way to explore KPZ limits and also possibly other large-scale behaviors for more general particle systems, including ``slow-boundary" generalizations of models studied in \cite{BDGN,EGN1,EGN2,GJNM} whose KPZ limits we believe are within reach based on the methods that are developed herein. We will pursue these directions in future research.

Let us now return to a deeper discussion on the history of KPZ universality and slow bond exclusion processes, starting with the former. Provided the success of Hairer's regularity structures in proving universality of the KPZ equation for SPDE growth models \cite{HQ}, it is quite natural to ask if these methods can analyze interacting particle system growth models/height functions. As we will soon illustrate, the height function $\mathbf{h}^{N}$ indeed solves a discretization of an SPDE that can be accessed with regularity structures as in \cite{HQ}. However, the details of this discretization are quite crucial and contain a major challenge when trying to employ regularity structures. In particular, the noise term in this SPDE discretization satisfied by $\mathbf{h}^{N}$ actually interacts with the both solution $\mathbf{h}^{N}$ and interacting particle system in a highly nontrivial way. This obstructs accessing refined information, such as microscopic correlation statistics for $\mathbf{h}^{N}$, that is necessary when employing the theory of regularity structures. This is one key difficulty in using regularity structures to attack Theorem \ref{theorem:kpz}, and there certainly may be others as well. As for the history of slow bond exclusion processes, these models happen to be prevalent in both statistical physics and probability communities, serving as a prototype for observing the changes in large-scale behaviors of highly complex systems in response to a microscopic perturbation; see \cite{BSaS,BSiS,MLM,JL1,JL2} for just a handful of relevant papers, of which at least the first two keep an eye towards KPZ statistics. Roughly speaking, a major lesson from these papers is that even microscopic perturbations such as slow bonds can induce drastic change in large-scale statistics; this issue was actually a topic of intense debate and speculation between physicists and mathematics until \cite{BSiS}. In particular, one view of Theorem \ref{theorem:kpz} and the main result of the current work is a homogenization result that adds more fuel to the aforementioned debate, as Theorem \ref{theorem:kpz} asserts the slow bond defect, for general slow bond particle systems, does not change large-scale behaviors of the height function. We conclude this extended discussion on slow bond models by citing \cite{FGN1,FGN2,FGN3,FGN4,FGN5,FN,FNV}, which study other non-KPZ aspects of homogenization theory and therefore naturally opens the door towards studying KPZ-related homogenization phenomena in these models that would only add more to the very modest literature on universality of the KPZ equation.

Let us now return to why Theorem \ref{theorem:kpz} is not perturbative in $\beta_{\star}$, despite allowing for choosing sufficiently small but still uniformly bounded from below and $N$-independent $\bar{\beta}_{\star}$. More generally, let us illustrate briefly the starting point for proving Theorem \ref{theorem:kpz} in the context of previous KPZ work \cite{BG,CST,CT,DT,FGS}. Looking at previous work of \cite{FGS} for slow bond models and SPDEs, we have the following leading-order SPDE for the microscopic height function $\mathbf{h}^{N}$, in which $c_{\beta_{\star}}\neq0$ depends only on $\beta_{\star}>0$:
\begin{align}
\partial_{\t}\mathbf{h}_{\t,\x}^{N} \ \sim \ 2^{-1}\Delta\mathbf{h}_{\t,\x}^{N} - 2^{-1}|\grad\mathbf{h}^{N}_{\t,\x}|^{2} + \mathrm{R}_{N} + \zeta^{N} + \mathbf{1}_{\x=0}c_{\beta_{\star}}N^{3/2}\mathfrak{g}_{\t,\x}, \label{eq:intro1}
\end{align}
where $\zeta^{N}$ is an approximation to the space-time white noise in KPZ, where $\mathrm{R}_{N}$ is the constant in Definition \ref{definition:hf} that is cancelled when we subtract from $\mathbf{h}^{N}$ the drift $\mathrm{R}_{N}\t$ in the Gartner transform, and $\mathfrak{g}_{\t,\x}\sim \eta_{\t,\x+1}-\eta_{\t,\x}$. Because $\mathfrak{g}$ is locally fluctuating, which means oscillating as the particle system evolves and swaps $\eta_{\t,\x}$ and $\eta_{\t,\x+1}$ values, we expect $\mathfrak{g}$ to have asymptotically negligible effect in the large-$N$ limit after we average over large space-time sets and turn the $N^{3/2}$ factor into $N^{1/2}$ since $\mathfrak{g}$ is localized in \eqref{eq:intro1} to a single spatial point. This is a key input of \cite{FGS}, which crucially requires the model to be globally stationary. For non-stationary models, microscopic correlation statistics become sufficiently intractable so the \cite{FGS} method is well beyond current reach. Namely, \eqref{eq:intro1} is significantly difficult to study for general non-stationary models, especially those with intractable invariant measures. We emphasize the slow bond parameter $\beta_{\star}$ appears in the coefficient $c_{\beta_{\star}}$. Even if $\beta_{\star}>0$ is arbitrarily small, as long as it is independent of $N$, the coefficient $c_{\beta_{\star}}$ is uniformly bounded away from zero, forcing us to deal with the last difficult term from the RHS of \eqref{eq:intro1} that we just discussed in this paragraph. We would only be able to avoid it if we choose $\beta_{\star}$ vanishing with $N$, which is not our case.

To work around the aforementioned obstruction to \eqref{eq:intro1}, let us consider instead the following \emph{approximate} microscopic KPZ equation, in which the Laplacian in \eqref{eq:intro1} is replaced by a non-constant coefficient Laplacian $\mathsf{A}(\x)\Delta$, where $\mathsf{A}$ is constant outside $\x=0$ and asymptotically vanishing at $\x=0$ to reflect the slow bond in the model. The benefit behind this cost of asymptotically non-elliptic operator is ultimately in forgetting $\mathfrak{g}$ in \eqref{eq:intro1} and replacing it with $N$ times the functional $\mathfrak{w}_{\t,\x}\sim\eta_{\t,\x}\eta_{\t,\x+1}$ along with replacing $c_{\beta_{\star}}$ with another constant $c_{\beta_{\star}}'\neq0$ that also depends only on $\beta_{\star}>0$:
\begin{align}
\partial_{\t}\mathbf{h}_{\t,\x}^{N} \ \sim \ 2^{-1}\mathsf{A}(\x)\Delta\mathbf{h}_{\t,\x}^{N} - 2^{-1}|\grad\mathbf{h}^{N}_{\t,\x}|^{2} + \mathrm{R}_{N} + \zeta^{N} + \mathbf{1}_{\x=0}c_{\beta_{\star}}'N\mathfrak{w}_{\t,\x}. \label{eq:intro2}
\end{align}
We note that $\mathsf{A}(0)\gtrsim N^{-\beta_{\star}}$ again reflecting the slow bond effect. This perspective of \eqref{eq:intro2} instead of \eqref{eq:intro1} then forces us to perform two major steps instead of one, but each of these two steps are much more accessible. The first step amounts to estimates for the heat operator $\partial_{\t}-2^{-1}\mathsf{A}(\x)\Delta$ and showing that at large scales, the support of $\mathsf{A}(\x)-1$ being the origin makes the coefficient $\mathsf{A}(\x)$ negligible. We approach this problem via a two-step heat kernel estimate strategy, the first step being a priori estimates via Nash-ian heat kernel analysis that has not been used in the literature for asymptotically degenerate diffusive particle systems to the author's knowledge, and the second step being a perturbative approach that exploits the local nature of the support of $\mathsf{A}(\x)-1$. Thus, we are then left with analyzing $N\mathfrak{w}$ in \eqref{eq:intro2}. Again, after averaging over large space-time sections, we may drop the $N$-factor, at which point we need to analyze the leading-order behavior of $\eta_{\t,0}\eta_{\t,1}$. Such asymptotic behavior is a major point of hydrodynamic limits. However, $\eta_{\t,0}\eta_{\t,1}$ is spatially localized. This is what leads us to develop a novel method to derive \emph{local} hydrodynamics for general interacting particle systems, and it ultimately implies $\eta_{\t,0}\eta_{\t,1}\sim0$ at large scales on which $\eta_{\t,0}$ and $\eta_{\t,1}$ fluctuate; see Definition \ref{definition:model} regarding the fluctuating property of $\eta_{\t,0}$. One takeaway and/or interpretation of this observation to look at \eqref{eq:intro2} instead of \eqref{eq:intro1} is replacing a highly difficult probabilistic problem by a much more tractable probabilistic problem at the cost of an analytic problem, the last of which are usually more more accessible given the extensive history and literature on partial differential equations.
%%%
\subsection{Organization}
%%%
In Section \ref{section:mshe}, we make precise the above approximate microscopic KPZ at the level of the Gartner transform $\mathbf{G}^{N}$. In this section, we also modify this Gartner transform slightly and conclude with both an outline of the proof of Theorem \ref{theorem:kpz} and an outline for the rest of this paper. In Section \ref{section:hke}, we develop heat kernel estimates for a rough and degenerate heat operator that is relevant because of the slow bond in the system. In Section \ref{section:lpe}, we present a main probabilistic ``local hydrodynamics" estimate alluded to throughout this introduction, whose proof is delayed until Section \ref{section:lpeproof} in order to avoid disrupting the flow of this article. In Section \ref{section:pc}, we basically reduce the proof of Theorem \ref{theorem:kpz} to the setting of earlier work \cite{BG,CST,CT} for integrable interacting particle systems, thereby forgetting the slow bond in some sense, and this is concluded by wrapping up the proof of Theorem \ref{theorem:kpz} in Section \ref{section:proof}. The appendix sections of this paper provide technical ingredients that are effectively straightforward modifications of estimates that are already established in previous literature \cite{BG,CST,CT,DT,Y,Y20} but are, strictly speaking, necessary for our analysis. Lastly, the final appendix section of this paper contains notation that is used frequently and serves as a source for the reader to refer to.
%%%
\subsection{Acknowledgements}
%%%
The author thanks Amir Dembo for advice, guidance, and useful discussion. The author and this work were funded by a fellowship from the Northern California chapter of the ARCS Foundation. The author would also like to thank anonymous referees for their invaluable input.
%
%
%
%%%
\section{Approximate Microscopic SHE}\label{section:mshe}
%%%
As promised in the previous section, we start with the following result; it can be interpreted as the non-elliptic height function SPDE but at the level of the Gartner transform. This is the key starting point of our work, after which we may present a detailed outline of the methods in this paper along with a ``heat-operator" version of the following stochastic differential dynamics.
%%%
\begin{prop}\label{prop:mshe1}
\fsp We have the following stochastic differential equation with notation introduced afterwards:
\begin{align}
\d_{\t}\mathbf{G}^{N}_{\t,\x} \ = \ 2^{-1}\Delta^{N,\beta_{\star}}_{\x}\mathbf{G}^{N}_{\t,\x}\d\t + \mathbf{G}_{\t,\x}^{N}\d\zeta_{\t,\x}^{N} + 2^{-1}c_{N}\mathbf{1}_{\x=0}\mathfrak{q}_{\t,\x}^{\mathrm{tot}}\mathbf{G}_{\t,\x}^{N}\d\t. \label{eq:mshe1}
\end{align}
%
%%%
\begin{itemize}
\item We define the non-elliptic operator $\Delta^{N,\beta_{\star}}_{\x}=N^{2}\mathbf{A}^{N,\beta_{\star}}_{\x}\Delta_{\x}$, where the static field in front of the Laplacian is defined via $\mathbf{A}_{\x}^{N,\beta_{\star}}=\mathbf{1}_{\x\neq0}+N^{-\beta_{\star}}\mathbf{1}_{\x=0}$; in words, it is a constant coefficient except for an asymptotically non-elliptic contribution supported at the slow bond. To be completely clear, we define $\Delta_{\x}\phi_{\x}=\phi_{\x+1}+\phi_{\x-1}-2\phi_{\x}$. We also defined the constant $c_{N}=N-N^{1-\beta_{\star}}$.
\item The differential martingale $\d\zeta^{N}_{\t,\x}$ is the ``infinitesimal increment" defined by specifying $\zeta^{N}_{\t,\x}$ as the compensation of the Poisson process that counts jumps in both directions (when allowed) between $\x$ and $\x+1$ until time $\t$; see the beginning of \emph{Section 2} in \cite{DT} but with the possible slow bond defect additionally introduced and for maximal jump-length equal to 1.
\item We defined the uniformly bounded functional $\mathfrak{q}_{\t,\x}^{\mathrm{tot}}=\mathfrak{q}_{\t,\x}+N^{-1/2}\wt{\mathfrak{q}}_{\t,\x}$, where $\mathfrak{q}_{\t,\x}=\eta_{\t,\x}\eta_{\t,\x+1}$ and $\wt{\mathfrak{q}}$ is uniformly bounded.
\end{itemize}
%%%
\end{prop}
%%%
In short, the proof of Proposition \ref{prop:mshe1} is inspecting the locality of $\mathbf{G}^{N}$ dynamics; away from the slow bond, we can cite Section 2 of \cite{DT}, for example, and at the slow bond, it, basically, suffices to observe that the growth of $\mathbf{G}^{N}$ is asymptotically ``slowed-down" by a factor of $N^{-\beta_{\star}}$, and the additional last term in \eqref{eq:mshe1} is the quadratic factor that occurs in KPZ growth. However, before we give the proof, we introduce a heat-operator version of \eqref{eq:mshe1}. After all, we will analyze the SDE-type equation \eqref{eq:mshe1} not by directly viewing as a differential equation but rather in its integrated form, similar to $\mathrm{SHE}$.
%%%
\begin{definition}\label{definition:mshe2}
\fsp Let us first define $\mathbf{H}^{\beta_{\star}}$ as the heat kernel on $\R_{\geq0}^{2}\times\Z^{2}$, though we only consider the first time-coordinate to be at most the second time-coordinate, associated to the following semi-discrete parabolic equation, so that for any $(\s,\t,\x,\y)$, we have
\begin{align}
\partial_{\t}\mathbf{H}^{\beta_{\star}}_{\s,\t,\x,\y} \ = \ 2^{-1}\Delta_{\x}^{N,\beta_{\star}}\mathbf{H}^{\beta_{\star}}_{\s,\t,\x,\y} \quad\mathrm{and}\quad \mathbf{H}^{\beta_{\star}}_{\s,\s,\x,\y} \ = \ \mathbf{1}_{\x=\y}.
\end{align}
We additionally define $\mathbf{H}^{\beta_{\star},\mathbf{X}}$ as the associated spatial heat operator and $\mathbf{H}^{\beta_{\star},\mathbf{T}}$ as a time-integrated version, both of these defined below in which $\phi:\R_{\geq0}\times\Z\to\R$ is a generic function; below, we will establish a few different notations that all mean the same thing but may be more convenient or clear depending on the context, and let us also emphasize the spatial heat operator $\mathbf{H}^{\beta_{\star},\mathbf{X}}$ acts on space-time functions only through their time-zero data:
\begin{align}
\mathbf{H}^{\beta_{\star},\mathbf{X}}_{\t,\x}(\phi) \ &= \ \mathbf{H}^{\beta_{\star},\mathbf{X}}_{\t,\x}(\phi_{0,\y}) \ = \ \mathbf{H}^{\beta_{\star},\mathbf{X}}_{\t,\x}(\phi_{0,\bullet}) \ = \ {\sum}_{y\in\Z}\mathbf{H}^{\beta_{\star}}_{0,\t,\x,\y}\phi_{0,\y} \\
\mathbf{H}^{\beta_{\star},\mathbf{T}}_{\t,\x}(\phi) \ &= \ \mathbf{H}^{\beta_{\star},\mathbf{T}}_{\t,\x}(\phi_{\s,\y}) \ = \ \mathbf{H}^{\beta_{\star},\mathbf{T}}_{\t,\x}(\phi_{\bullet,\bullet}) \ = \ \int_{0}^{\t}{\sum}_{y\in\Z}\mathbf{H}^{\beta_{\star}}_{\s,\t,\x,\y}\phi_{\s,\y}\d\s.
\end{align}
Observe the previous constructions extend to all $\beta_{\star}$, not just the choice of slow-bond exponent we have made in this paper. In this paper, we will also use the construction for $\beta_{\star}=0$, and simply remove $\beta_{\star}$ from the notation for the heat kernel, whose associated Laplacian is constant-coefficient, in this case. We will re-clarify this when these $\beta_{\star}=0$ operators become relevant.
\end{definition}
%%%
The following result is a classical semi-discrete version of the Duhamel principle and occurs frequently in previous KPZ works \cite{BG,CT,DT}, for example, so we do not provide a proof. If of interest, however, we note it follows by standard linear ODE theory.
%%%
\begin{corollary}\label{corollary:mshe3}
\fsp Retaining the notation in \emph{Proposition \ref{prop:mshe1}}, we have the following integral equation:
\begin{align}
\mathbf{G}_{\t,\x}^{N} \ = \ \mathbf{H}^{\beta_{\star},\mathbf{X}}_{\t,\x}(\mathbf{G}^{N}) + \mathbf{H}^{\beta_{\star},\mathbf{T}}_{\t,\x}(\mathbf{G}^{N}\d\zeta^{N}) + \mathbf{H}^{\beta_{\star},\mathbf{T}}_{\t,\x}\left(2^{-1}c_{N}\mathbf{1}_{\y=0}\mathfrak{q}_{\s,\y}^{\mathrm{tot}}\mathbf{G}_{\s,\y}^{N}\right). \label{eq:mshe3}
\end{align}
\end{corollary}
%%%
%%%
\begin{proof}[Proof of \emph{Proposition \ref{prop:mshe1}}]
If $\x\neq0$, it suffices to then follow the calculations in Section 2 of \cite{DT} for single-length jumps, because there is no slow bond effect for such $\x$, thus we will assume $\x=0$ for the rest of this argument. An explicit calculation with Poisson clocks as with the proof of Proposition 2.2 in \cite{DT} provides us the following first step, which can be seen as basically Proposition 2.2 and Section 2 of \cite{DT} but we have slowed-down the symmetric-jump speed across $\{0,1\}$ by a factor of $N^{-\beta_{\star}}$; in case it is of any clarity, we emphasize that $\mathbf{G}^{N}_{\t,0}$ changes only by jumps in the particle system across $\{0,1\}$, and moreover, the $\zeta^{N}$-noise comes from separating the stochastic Poisson-driven growth of $\mathbf{G}^{N}$ into a drift/infinitesimal-expected growth and resulting martingale:
\begin{align}
N^{-2}\d_{\t}\mathbf{G}^{N}_{\t,0} \ = \ \Phi^{\mathrm{s}}_{\t,0}\mathbf{G}^{N}_{\t,0}\d\t + \Phi^{\mathrm{a}}_{\t,0}\mathbf{G}^{N}_{\t,0}\d\t + \mathrm{R}_{N}\mathbf{G}^{N}_{\t,0}\d\t + \mathbf{G}^{N}_{\t,0}\d\zeta_{\t,0}^{N}. \label{eq:mshe11}
\end{align}
Above, we recall the renormalization constant $\mathrm{R}_{N}$ from the definition of the Gartner transform $\mathbf{G}^{N}$. We have also introduced the following quantities, the first of which encodes the speed of particle jumps across $\{0,1\}$ in the symmetric component of the random walk and the second encodes the speed of the asymmetric component; we emphasize only the first symmetric component sees the slow bond effect/factor of $N^{-\beta_{\star}}$:
\begin{align*}
\Phi_{\t,0}^{\mathrm{s}} \ = \ 8^{-1}N^{-\beta_{\star}}(1+\eta_{\t,1})(1-\eta_{\t,0})(\exp(-2N^{-1/2})-1) + 8^{-1}N^{-\beta_{\star}}(1-\eta_{\t,1})(1+\eta_{\t,0})(\exp(2N^{-1/2})-1) \\
\Phi_{\t,0}^{\mathrm{a}} \ = \ 8^{-1}N^{-1/2}(1+\eta_{\t,1})(1-\eta_{\t,0})(\exp(-2N^{-1/2})-1) - 8^{-1}N^{-1/2}(1-\eta_{\t,1})(1+\eta_{\t,0})(\exp(2N^{-1/2})-1).
\end{align*}
For complete transparency, we also record below the following Laplacian calculation that follows by elementary considerations as in the proof of Proposition 2.2/Section 2 of \cite{DT}; the Laplacian acts on the spatial variable of $\mathbf{G}^{N}$ and is evaluated at the origin:
\begin{align}
2^{-1}N^{-\beta_{\star}}\Delta_{\x}\mathbf{G}^{N}_{\t,0} \ = \ 2^{-1}N^{-\beta_{\star}}(\exp(-N^{-1/2}\eta_{\t,1})-1)\mathbf{G}^{N}_{\t,0} + 2^{-1}N^{-\beta_{\star}}(\exp(N^{-1/2}\eta_{\t,0})-1)\mathbf{G}^{N}_{\t,0}. \label{eq:mshe12}
\end{align}
The proposed equation \eqref{eq:mshe1} then follows from \eqref{eq:mshe11}, \eqref{eq:mshe12}, and a host of Taylor expansions of exponentials and power-counting. We describe these expansions as follows. First, the contribution of $\Phi^{\mathrm{s}}\mathbf{G}^{N}$ in \eqref{eq:mshe11} will match that of the LHS of \eqref{eq:mshe12} times $\d\t$ up to order $N^{-3/2}$ terms; this is just because this was the case in the calculations of Section 2 (see Proposition 2.2 and Proposition 2.3 therein) of \cite{DT}, and we have only multiplied each by the same small factor $N^{-\beta_{\star}}$. Meanwhile, looking at the same calculations, the contribution of $\Phi^{\mathrm{a}}\mathbf{G}^{N}$ matches that of $\mathrm{R}_{N}\mathbf{G}^{N}$ up to the error term of $2^{-1}c_{N}\mathfrak{q}^{\mathrm{tot}}$ from the statement of Proposition \ref{prop:mshe1}. Indeed, the renormalization constant $\mathrm{R}_{N}\mathbf{G}^{N}$ is designed to cancel the order $N^{-1}$ constant in $\Phi^{\mathrm{a}}\mathbf{G}^{N}$ in \cite{BG,DT}; what is left is the quadratic term $\mathfrak{q}_{\t,0}=\eta_{\t,0}\eta_{\t,1}$ times an order $N^{-1}$ constant plus order $N^{-3/2}$ terms. We emphasize that this last calculation for $\Phi^{\mathrm{a}}\mathbf{G}^{N}$ and $\mathrm{R}_{N}\mathbf{G}^{N}$ is blind to the slow bond effect and can be found in \cite{BG,DT} as written. This completes the proof after multiplying everything by $N^{2}$ to account for the $N^{2}$ speed of the particle system that we removed in \eqref{eq:mshe11}.
\end{proof}
%%%
%%%
\subsection{Technical Modification}
%%%
Based on the stochastic equations of Proposition \ref{prop:mshe1} and Corollary \ref{corollary:mshe3}, we will now consider a modification of $\mathbf{G}^{N}$ along the following lines. Technically, the $\d\zeta^{N}$-field from Proposition \ref{prop:mshe1} and Corollary \ref{corollary:mshe3} is not a \emph{globally} analytically-tractable field, since it is not ``uniformly bounded" in the following sense. In any positive time, the number of jumps seen in the $\d\zeta^{N}$-field is infinite, even if the time-scale is sub-microscopic, because $\d\zeta^{N}$-field is driven by infinitely many Poisson clocks of non-vanishing speeds. This introduces only few technical issues in later analysis when we want to control the total number of Poisson clock rings that are relevant on macroscopic length-scales in a very sub-microscopic time-scale; although we could find ways to avoid having to control such a total number of clock rings, doing so would make our writing severely less clear and require many more technical constructions. We now modify the $\d\zeta^{N}$-field by cutting off its contribution outside a very super-macroscopic length-scale, so that as far as Theorem \ref{theorem:kpz} is concerned we have barely changed things, and consider a solution to the stochastic equation driven by this modification, denoted by $\d\xi^{N}$. Below, the exponent 30 is one fixed choice much larger than the length-scale exponent of 1, and it could be replaced by anything larger or slightly smaller without affecting the ideas of this paper.
%%%
\begin{definition}\label{definition:mshe4}
\fsp Define $\d\xi_{\t,\x}^{N}=\mathbf{1}_{|\x|\leq N^{30}}\d\zeta_{\t,\x}^{N}$. Define $\mathbf{Z}^{N}$ as the solution to the following stochastic equation on $\R_{\geq0}\times\Z$:
\begin{align}
\d_{\t}\mathbf{Z}^{N}_{\t,\x} \ = \ 2^{-1}\Delta^{N,\beta_{\star}}_{\x}\mathbf{Z}^{N}_{\t,\x}\d\t + \mathbf{Z}_{\t,\x}^{N}\d\xi_{\t,\x}^{N} + 2^{-1}c_{N}\mathbf{1}_{\x=0}\mathfrak{q}_{\t,\x}^{\mathrm{tot}}\mathbf{Z}_{\t,\x}^{N}\d\t \and \mathbf{Z}_{0,\x}^{N}=\mathbf{G}_{0,\x}^{N}.
\end{align}
In particular, we have the following stochastic integral equation for $\mathbf{Z}^{N}$ via Duhamel expansion:
\begin{align}
\mathbf{Z}_{\t,\x}^{N} \ = \ \mathbf{H}^{\beta_{\star},\mathbf{X}}_{\t,\x}(\mathbf{G}^{N}) + \mathbf{H}^{\beta_{\star},\mathbf{T}}_{\t,\x}(\mathbf{Z}^{N}\d\xi^{N}) + \mathbf{H}^{\beta_{\star},\mathbf{T}}_{\t,\x}\left(2^{-1}c_{N}\mathbf{1}_{\y=0}\mathfrak{q}_{\s,\y}^{\mathrm{tot}}\mathbf{Z}_{\s,\y}^{N}\right).
\end{align}
Note that for the spatial heat operator in the previous stochastic integral equation, we can also put $\mathbf{Z}^{N}$ in place of $\mathbf{G}^{N}$, because we have stipulated in the stochastic differential equation prior that their initial data agree. We will sometimes refer to $\mathbf{Z}^{N}$ as the Gartner transform because it is basically the Gartner transform $\mathbf{G}^{N}$ modulo forgetting noise that is too far from the scales of interest to make a difference with enough probability; equivalently, calling $\mathbf{Z}^{N}$ the Gartner transform is justified by Lemma \ref{lemma:mshe5} below.
\end{definition}
%%%
Before we state the following result, it will be convenient to introduce a family of norms that we will often use when analyzing stochastic equations of $\mathbf{G}^{N}$ and $\mathbf{Z}^{N}$ type. We will recall this norm later in the paper when we turn to using it again.
%%%
\begin{definition}\label{definition:mshe4b}
\fsp Provided any function $\phi:\R_{\geq0}\times\Z\to\R$ along with any subsets $[0,\t_{+}]\subseteq\R_{\geq0}$ and $\mathbb{K}\subseteq\R$, which could be a discrete subset of $\Z$, for example, we define the norm $\|\phi\|_{\t_{+};\mathbb{K}}=\sup_{0\leq\t\leq\t_{+}}\sup_{\x\in\mathbb{K}}\exp(-|\x|/N)|\phi_{\t,\x}|$.
\end{definition}
%%%
%%%
\begin{lemma}\label{lemma:mshe5}
\fsp Take any compact set $\mathbb{K}\subseteq\R$ and define the microscopic coordinate set $\mathbb{K}^{N}=N\mathbb{K}\cap\Z$. We have the following outside an event of probability at most of order $N^{-100}$, in which the implied constant depends only on $\mathbb{K}$ and is independent of $N$:
\begin{align}
\|\mathbf{G}^{N}_{\t,\x}-\mathbf{Z}^{N}_{\t,\x}\|_{1;\mathbb{K}^{N}} \ \lesssim_{\mathbb{K}} \ N^{-100}.
\end{align}
\end{lemma}
%%%
%%%
\begin{remark}\label{remark:mshe6}
\fsp In view of Lemma \ref{lemma:mshe5}, because of the locally uniform topology we have put on the Skorokhod space in the statement of Theorem \ref{theorem:kpz}, it will suffice to prove Theorem \ref{theorem:kpz} after replacing $\mathbf{G}^{N}$ therein with $\mathbf{Z}^{N}$ from Definition \ref{definition:mshe4}. This is what we do.
\end{remark}
%%%
%%%
\begin{remark}\label{remark:mshe7}
\fsp The proof of Lemma \ref{lemma:mshe5} that we provide below is basically the Gronwall argument of Proposition 3.2 and Corollary 3.3 in \cite{DT}. We briefly explain what it is saying. The difference between the $\mathbf{Z}^{N}$ and $\mathbf{G}^{N}$ equations is $\d\zeta^{N}$ versus $\d\xi^{N}$, which only differ on the set of $|\x|>N^{30}$. Because $N^{30}$ is much farther than a simple random walk described by the heat kernel would travel in time $N^{2}$, we expect $\mathbf{G}^{N}\approx\mathbf{Z}^{N}$ to fail with exponentially low probability in a high power of $N$, and we expect that the difference $\mathbf{G}^{N}-\mathbf{Z}^{N}$ is similarly exponentially small in a high power of $N$. However, and this is a technical point, because the $\mathbf{G}^{N}$ and $\mathbf{Z}^{N}$ equations in Corollary \ref{corollary:mshe3} and Definition \ref{definition:mshe4} are both linear and multiplicative in their respective solutions, the difference $\mathbf{G}^{N}-\mathbf{Z}^{N}$ is also going to scale linearly in $\mathbf{G}^{N}$. However, the proof of Proposition 3.2 and Corollary 3.3 in \cite{DT} ensures that $\mathbf{G}^{N}$ is at most exponential in a bounded power of $N$ that comes from the fact that $\mathbf{G}^{N}$ grows at a speed of order at most $N^{2}$ due to our choice of scaling for the particle system. This bound for $\mathbf{G}^{N}$ is then easily beaten out by the exponentially low in a high power of $N$ factor from before. Given this explanation, the reader is invited to skip the proof of Lemma \ref{lemma:mshe5} because it is fairly beyond the main ideas of this paper and it is not saying anything too deep; it just gives a technical convenience explained prior to Definition \ref{definition:mshe4}.
\end{remark}
%%%
%%%
\begin{proof}[Proof of \emph{Lemma \ref{lemma:mshe5}}]
We first observe that the original Gartner transform $\mathbf{G}^{N}$ solves the following stochastic differential equation on $\R_{\geq0}\times\Z$, in which we replace the operator $\Delta^{N,\beta_{\star}}$ from Proposition \ref{prop:mshe1} with the constant-coefficient Laplacian $N^{2}\Delta$:
\begin{align}
\d_{\t}\mathbf{G}^{N}_{\t,\x} \ = \ 2^{-1}N^{2}\Delta\mathbf{G}^{N}_{\t,\x}\d\t + 2^{-1}(\Delta_{\x}^{N,\beta_{\star}}-N^{2}\Delta)\mathbf{G}^{N}_{\t,\x}\d\t + \mathbf{G}_{\t,\x}^{N}\d\zeta^{N}_{\t,\x} + 2^{-1}c_{N}\mathbf{1}_{\x=0}\mathfrak{q}_{\t,\x}^{\mathrm{tot}}\mathbf{G}^{N}_{\t,\x}\d\t. \label{eq:mshe51}
\end{align}
Observe the second term on the RHS of this stochastic equation is supported at $\x=0$ and can be realized as a linear combination of gradients of $\mathbf{G}^{N}$ scaled by $N^{2}$. By definition of the Gartner transform $\mathbf{G}^{N}$ and Taylor expansion, this can be written as an order $N^{3/2}$ coefficient times $\mathbf{G}^{N}$ itself at one point. Therefore, we may run the Gronwall inequality argument in the proof of Proposition 3.2 in \cite{DT} in order to deduce that uniformly in $(\t,\x)\in[0,1]\times\Z$, for any $p\geq1$ we have the estimate below for some $\kappa\geq0$, since the initial data of the Gartner transform $\mathbf{G}^{N}$ has all moments bounded uniformly over $\Z$ by assumption:
\begin{align}
\E|\mathbf{G}_{\t,\x}^{N}|^{2p} \ \lesssim_{p} \ \exp(\kappa p N^{3/2}). \label{eq:mshe52}
\end{align}
We additionally note that the difference $\mathbf{K}^{N}=\mathbf{G}^{N}-\mathbf{Z}^{N}$ solves the stochastic equation
\begin{align}
\d_{\t}\mathbf{K}^{N}_{\t,\x} \ = \ 2^{-1}N^{2}\Delta\mathbf{K}^{N}_{\t,\x}\d\t + 2^{-1}(\Delta_{\x}^{N,\beta_{\star}}-N^{2}\Delta)\mathbf{K}^{N}_{\t,\x}\d\t + \mathbf{K}_{\t,\x}^{N}\d\xi^{N}_{\t,\x} + 2^{-1}c_{N}\mathbf{1}_{\x=0}\mathfrak{q}_{\t,\x}^{\mathrm{tot}}\mathbf{K}^{N}_{\t,\x}\d\t + \mathbf{1}_{|\x|>N^{30}}\mathbf{G}^{N}_{\t,\x}\d\zeta_{\t,\x}^{N}. \nonumber
\end{align}
We may again run the Gronwall argument of Proposition 3.2 in \cite{DT}; this time, we observe that said Gronwall argument allows us to multiply every $\mathbf{K}^{N}$ and $\mathbf{G}^{N}$ term with an exponentially-decaying factor $\exp(-N^{-1}|\x|)$, because the heat kernel $\mathbf{H}$ defined by $\partial_{\t}-N^{2}\Delta$ can absorb factors of the form $\exp(N^{-1}|\x|)$. We again observe that the second term in the $\mathbf{K}^{N}$ equation can be realized as a linear combination of gradients of $\mathbf{K}^{N}$ at $\x=0$ scaled by $N^{2}$. Moreover, we now observe that because of the indicator function in the last $\mathbf{G}^{N}$-term in the $\mathbf{K}^{N}$ equation is supported outside of $|\x|\leq N^{30}$, when we multiply $\mathbf{G}^{N}$ by $\exp(-N^{-1}|\x|)$ for $|\x|>N^{30}$ the resulting $2p$-moments are exponentially small in $N$ courtesy of \eqref{eq:mshe52}. Therefore, the Gronwall inequality argument in the proof of Proposition 3.2 in \cite{DT} provides the following estimate for $\mathbf{K}^{N}$ given any $p\geq1$ and for some $\kappa\geq0$, in which $\d(\x)$ denotes the distance between $\x\in\Z$ and the set of $|\x|>N^{30}$. We emphasize the following is uniform over $(\t,\x)\in[0,1]\times\Z$. Moreover, on the RHS of \eqref{eq:mshe53}, the first factor is the exponential growth that comes from application of the Gronwall inequality. Such factor will hit both the initial data of $\mathbf{K}^{N}$ and the error term or ``forcing" given by the last term in the $\mathbf{K}^{N}$ equation above. The initial data of $\mathbf{K}^{N}$ vanishes by construction of $\mathbf{K}^{N}$ and $\mathbf{Z}^{N}$. The aforementioned forcing/last term in the $\mathbf{K}^{N}$ equation is controlled by $\mathbf{G}^{N}$ moment estimates but restricted to the set of $|\x|>N^{30}$. Because of the exponentially decaying weight in \eqref{eq:mshe53} below that we remarked on earlier in this paper, the contribution of this forcing is penalized by a sub-exponential factor that comes from tails of the $\mathbf{H}$ heat kernel; it depends on how far we look from the set of $|\x|>N^{30}$. Ultimately, the proof of Proposition 3.2 in \cite{DT} gives
\begin{align}
\exp(-pN^{-1}|\x|)\E|\mathbf{K}_{\t,\x}^{N}|^{2p} \ \lesssim_{p} \ \exp(2\kappa pN^{2})\exp(-\kappa^{-1}p N^{-1}\d(\x))\left({\sup}_{(\t,\x)\in[0,1]\times\Z}\E|\mathbf{G}_{\t,\x}^{N}|^{2p}\right). \label{eq:mshe53}
\end{align}
Observe that for $|\x|\lesssim_{\mathbb{K}}N$, the RHS of \eqref{eq:mshe53} is exponentially small in $N$ once we apply \eqref{eq:mshe52} as $\d(\x)\geq 2^{-1}N^{30}$ in this case.

The same Gronwall argument with the modifications above to account for the support of the $\mathbf{G}^{N}$-term in the $\mathbf{K}^{N}$ equation also implies that $\mathbf{K}^{N}$ is Holder regular in the compact space-time set $[0,1]\times\mathbb{K}^{N}$ with Holder norm having all moments and depending on $\mathbb{K}$, up to jumps of order $N^{-1/2}\mathbf{K}^{N}$ coming from the speed $N^{2}$-Poisson clocks supported on $\mathbb{K}^{N}$ on microscopic time-scales. In particular, because we have pointwise control in \eqref{eq:mshe53} and Holder regularity of the same scale uniformly on $[0,1]\times\mathbb{K}^{N}$, we may deduce the desired estimate on an event of the required high probability.
\end{proof}
%%%
%%%
\subsection{Strategy}
%%%
We will describe the method behind the proof of Theorem \ref{theorem:kpz} provided the inputs of this section. First, as we note when we give the proof of Theorem \ref{theorem:kpz}, it is enough to prove Theorem \ref{theorem:kpz} upon replacing the original Gartner transform $\mathbf{G}^{N}$ by the slightly modified Gartner transform $\mathbf{Z}^{N}$ in Definition \ref{definition:mshe4}. The three major steps to prove Theorem \ref{theorem:kpz} for $\mathbf{Z}^{N}$ are then as follows.
%%%
\begin{itemize}
\item In the stochastic integral equation for $\mathbf{Z}^{N}$, remove the last term on the RHS.
\item After removing the last term on the RHS of the $\mathbf{Z}^{N}$ equation, change $\mathbf{H}^{\beta_{\star}}$ heat kernels to $\mathbf{H}$ heat kernels.
\item Show that Theorem \ref{theorem:kpz} holds for the solution to the $\mathbf{Z}^{N}$ equation after these two adjustments in the previous bullet points.
\end{itemize}
%%%
The third and final step in the above outline is basically the work of \cite{BG}, though we require an additional hydrodynamic limit input of \cite{DT} because adapting the method to identify the limit quadratic variation of $\d\xi^{N}$ and $\d\zeta^{N}$ martingale differentials would require much work. In any case, however, this third step is basically well-known in the literature. The second step in the three-step method requires a heavy analysis of heat kernels $\mathbf{H}^{\beta_{\star}}$. Roughly speaking, the defect in the heat kernel $\mathbf{H}^{\beta_{\star}}$ is \emph{localized}, whereas Theorem \ref{theorem:kpz} looks for behavior on large space-time scales, for which microscopic behavior should be irrelevant; moreover, the slow bond defect is not strong enough to lead to any changes in the behavior of the random walk whose transition probabilities are given by the $\mathbf{H}^{\beta_{\star}}$ heat kernel. We discuss the picture for heat kernel estimates in more detail later in this section. The first step in the previous outline is based on the idea that the leading-order quadratic term in $\mathfrak{q}^{\mathrm{tot}}$ has vanishing hydrodynamic limit behavior, and thus the contribution of the last term on the RHS of the $\mathbf{Z}^{N}$ equation, which is basically at hydrodynamic scale modulo powers of $N^{\beta_{\star}}$ that come from technical calculations, is asymptotically vanishing. But there are several obstructions to making this precise, however, which include the fact that hydrodynamic limit estimates are generally with respect to a topology that is too weak to help us control the $\mathbf{Z}^{N}$ equation in any normed sense that is usually required to make sense of multiplicative SPDEs; in particular, we must make the hydrodynamic limit estimates of \cite{GPV,KL} quantitative, for example. The most important hurdle, however, is that hydrodynamic limit behaviors are generally established in the literature \cite{GPV,KL} at macroscopic space-time scales, whereas the contribution of the leading-order quadratic term in $\mathfrak{q}^{\mathrm{tot}}$ is \emph{localized in space}. In particular, we must develop a method to compute hydrodynamic limit behavior at \emph{local} scales. This local hydrodynamic limit is our novel innovation to the literature on large-scale interacting particle systems, and it is likely to provide future impact on research on interacting particle systems. We will discuss more the idea behind establishing these so-called local hydrodynamics in more detail in this section; we only mention here that all we need from $\mathfrak{q}^{\mathrm{tot}}$ and its leading-order $\mathfrak{q}$ is that its hydrodynamic limit is vanishing if $\eta\approx0$ in a spatially averaged sense initially. In particular, the specific quadratic nature of $\mathfrak{q}$, which might seem important given the context of the KPZ equation here, is not actually important.
%%%
\subsubsection{Heat Kernel Estimates}
%%%
We discuss first the content of Section \ref{section:hke}, which focuses on a plethora of heat kernel estimates for $\mathbf{H}^{\beta_{\star}}$ that, for example, compare $\mathbf{H}^{\beta_{\star}}$ to $\mathbf{H}$, estimate $\mathbf{H}^{\beta_{\star}}$ pointwise, and estimate regularity of $\mathbf{H}^{\beta_{\star}}$ for purposes that are actually important for estimating the last term on the RHS of the $\mathbf{Z}^{N}$ equation in Definition \ref{definition:mshe4}. In any case, we recall from earlier that we view $\mathbf{H}^{\beta_{\star}}$  basically perturbatively off $\mathbf{H}$. We make this precise via the Duhamel expansion, which is a resolvent Green's function perturbative expansion that is actually the same thing used to deduce Corollary \ref{corollary:mshe3} from Proposition \ref{prop:mshe1}; see Lemma \ref{lemma:hke5}. However, this perturbative expansion is difficult to use by itself to prove estimates for $\mathbf{H}^{\beta_{\star}}$ as the difference between $\mathbf{H}^{\beta_{\star}}$ and $\mathbf{H}$ equations is both in the leading-order derivative and in non-smooth fashion. Thus, to treat the non-smooth coefficient in front of the Laplacian in $\Delta^{N,\beta_{\star}}$, we require the approach of Nash for heat kernel estimates with rough coefficients. Usually, this is used for heat kernels of divergence-form equations, but in spatial dimension one, we show here it can be employed for non-divergence form equations as well; we could not find any instances of this in the heat kernel literature, and in case there actually is no literature for Nash heat kernel estimates for non-divergence-form heat kernels, our work would shed light towards that direction. Once we employ a Nash strategy to establish sub-optimal global a priori heat kernel estimates for $\mathbf{H}^{\beta_{\star}}$ and address the aforementioned non-smoothness in $\Delta^{N,\beta_{\star}}$, we may perturb using the Duhamel expansion. Perturbative expansions will also provide us sub-optimal but still sufficient regularity estimates for $\mathbf{H}^{\beta_{\star}}$. All of this is done in and is the purpose of Section \ref{section:hke}, which is key for the rest of this paper.
%%%
\subsubsection{Local Hydrodynamics}
%%%
Section \ref{section:lpe} is dedicated to stating a precise and quite strongly quantitative version of the following (what will be) fact -- the last term from the RHS of the $\mathbf{Z}^{N}$ equation in Definition \ref{definition:mshe4} vanishes in the large-$N$ limit in an appropriate space-time uniform norm, \emph{contingent upon} estimates for $\mathbf{Z}^{N}$ at the origin. We defer the proof of this to Section \ref{section:lpeproof} in order to avoid disrupting the flow of this paper, because the proof of the aforementioned vanishing statement, which is Proposition \ref{prop:lpe3}, is fairly long and technically involved. We provide here a short description and discussion of the proof of Proposition \ref{prop:lpe3}, or equivalently the content of Section \ref{section:lpeproof}. Recall from earlier that estimating the last term on the RHS of the $\mathbf{Z}^{N}$ equation amounts to replacing $\mathfrak{q}^{\mathrm{tot}}$ by its hydrodynamic limit. In classical works on hydrodynamic limits, beginning with \cite{GPV}, this is achieved by replacing one of the $\eta$-factors in the leading-order quadratic term $\mathfrak{q}$ from Proposition \ref{prop:mshe1} with an average of $\eta$-values on a bigger block, and afterwards matching said average with the hydrodynamic limit. Specifically, in the classical approach of \cite{GPV}, this requires spatially averaging $\mathfrak{q}^{\mathrm{tot}}$ and then applying a phenomenon of local equilibrium. Since we have no spatial averaging in the last term of the $\mathbf{Z}^{N}$ equation, due to the $\y=0$ localization therein, we will instead average $\mathfrak{q}^{\mathrm{tot}}$, or more precisely $\mathfrak{q}$, in time. The upshot here is that we only have to perform said time-average on \emph{local} time-scales, at which we can employ a quantitative version of local equilibrium in \cite{GPV} based on a \emph{local}, and therefore more robust, log-Sobolev inequality of \cite{Yau}, to implement time-averaging estimates for stationary particle systems instead of spatial-averaging estimates used in \cite{GPV} that ultimately end up being classical probability. One interpretation of our approach to \emph{local} hydrodynamics of this kind is a \emph{dynamical} and quantitative version of the classical one-block scheme of \cite{GPV}. In particular, we develop an improvement of the one-block scheme in \cite{GPV} that analyzes local and dynamical path-space functionals of the particle system via local equilibrium and technology that is exclusively available for stationary interacting particle systems and absolutely crucial for energy solution theory approaches to universality of KPZ, for example. In any case, after this replacement of $\mathfrak{q}^{\mathrm{tot}}$ by hydrodynamic limit, we must estimate this hydrodynamic limit, which we accomplish by interpreting said hydrodynamic limit in terms of a regularity estimate for $\mathbf{G}^{N}$, or equivalently $\mathbf{Z}^{N}$ up to an error with high-probability by Lemma \ref{lemma:mshe5}. Such an idea was actually used in \cite{DT} to estimate hydrodynamic limits but strictly at the fluctuation length-scale of order $N^{1/2}$. In the current work, we illustrate that said analysis holds for all mesoscopic scales, thereby highlighting its potential utility for future applications that include establishing a Boltzmann-Gibbs principle \cite{GJ15} and furthering our understanding of universality of the KPZ equation. We finish our brief discussion of local hydrodynamics by noting the aforementioned replacement by local hydrodynamic limit will require a classical entropy production and log-Sobolev inequality that are standard for the local equilibrium step in hydrodynamic limits \cite{GPV,Yau}. The estimate for localized stationary particle systems for time-averages of local functionals that we will employ is the classical Kipnis-Varadhan inequality; see Lemma \ref{lemma:lp11} and see \cite{GJ15} for further discussion about its importance for interacting particle systems more generally as well as necessary restriction to stationary models. Lastly, for the replacement step, because we must estimate terms in the $\mathbf{Z}^{N}$ equation in appropriately strong space-time norms, we also include a few technical gymnastics. A more detailed outline of the proof of Proposition \ref{prop:lpe3} and analysis for local hydrodynamics is given at the beginning of Section \ref{section:lpeproof}. We briefly remark that our method of local hydrodynamics may be improved with an additional localization of the two-blocks scheme in \cite{GPV} as well, but this would only complicate further the presentation of this paper and the resulting boost in $\bar{\beta}_{\star}$ and allowable $\beta_{\star}$ does not bring us significantly closer to the optimal value $\bar{\beta}_{\star}\approx2^{-1}$.
%%%
\subsubsection{Rest of this Paper}
%%%
Given Proposition \ref{prop:lpe3}, whose proof we briefly discussed immediately above, we make aforementioned removal of the last term in the $\mathbf{Z}^{N}$ equation and replacement of $\mathbf{H}^{\beta_{\star}}$ heat operators by $\mathbf{H}$ heat operators; this is the point of Section \ref{section:pc}. In Section \ref{section:proof}, we basically follow \cite{BG,DT} to prove Theorem \ref{theorem:kpz} for the resulting stochastic equation after removal and replacement in Section \ref{section:pc}. This would prove Theorem \ref{theorem:kpz} for $\mathbf{Z}^{N}$, from which we deduce Theorem \ref{theorem:kpz} for $\mathbf{G}^{N}$ itself via Lemma \ref{lemma:mshe5}. Section \ref{section:lpeproof} provides a proof for the key probabilistic ``local hydrodynamics" ingredient of Proposition \ref{prop:lpe3}. The appendix sections are auxiliary estimates that are either quite clear and/or straightforward or are minor adaptations on previous estimates in the literature \cite{BG,DT}.
%
%
%
%%%
\section{Heat Kernel Estimates}\label{section:hke}
%%%
The primary goal of this section is to provide us heat kernel estimates for $\mathbf{H}^{\beta_{\star}}$ from Definition \ref{definition:mshe2} in order to study terms in the heat equation \eqref{eq:mshe3} effectively. This section is rather lengthy, so we provide explanations for each estimate in order to clarify the exposition and assist the reader for an easier and swifter reading. In particular, let us provide the following outline.
%%%
\begin{itemize}
\item The coefficient in front of the perturbed Laplacian $\Delta^{N,\beta_{\star}}$ is neither elliptic nor smooth, so we cannot directly use any classical parabolic theory. The first subsection in this section is to adapt the famous Nash-approach to heat kernel estimates; in case this is of any help, the Nash approach is usually used for divergence-form equations, whereas $\Delta^{N,\beta_{\star}}$ is non-divergence-form, but this is not a major issue in spatial dimension one, which is the current setting. In any case, in this subsection we ultimately establish \emph{global} estimates for the heat kernel $\mathbf{H}^{\beta_{\star}}$ that are likely sub-optimal in view of an explicit invariant measure for the random walk described by $\mathbf{H}^{\beta_{\star}}$, but they are nonetheless enough for us if we assume $\beta_{\star}$ is not too small though still positive and fixed.
\item These aforementioned global estimates are \emph{certainly not} optimal if we look at the heat kernel away from the slow bond because of locality of the random walk described by $\mathbf{H}^{\beta_{\star}}$, or equivalently locality of the $\Delta^{N,\beta_{\star}}$ operator. Ultimately, outside a strictly-smaller-than-macroscopic neighborhood of the slow bond support, we have $\mathbf{H}^{\beta_{\star}}\approx\mathbf{H}$. This means the slow bond is effectively negligible in the large-$N$ limit at the level of heat kernels if $\beta_{\star}$ is positive and fixed but still small, which is a key reason why the Gartner transform limit does not see the slow bond, though we must make this sufficiently quantitative in this paper. Technically, we will use above global estimates as \emph{a priori} estimates in a local perturbative scheme to obtain optimal estimates for $\mathbf{H}^{\beta_{\star}}$ away from the slow bond. This perturbation is in $\beta_{\star}$, allowing us to gain something via $\beta_{\star}=0$ estimates in Appendix of \cite{DT}.
\item In a similar fashion, we perturb around $\beta_{\star}=0$ to obtain regularity estimates for $\mathbf{H}^{\beta_{\star}}$ as well. However, these will not be optimal.
\item Lastly, we obtain off-diagonal estimates for $\mathbf{H}^{\beta_{\star}}$; this will be done both perturbatively about $\beta_{\star}=0$, letting us use off-diagonal estimates in \cite{DT}, and with a Chapman-Kolmogorov equation for $\mathbf{H}^{\beta_{\star}}$ that lets us concretely realize off-diagonal estimates via elementary and standard random walk tail estimates.
\end{itemize}
%%%
We now establish some notation that we will use only in this section, unless otherwise mentioned.
%%%
\begin{definition}\label{definition:hke1}
\fsp Provided any function $\phi:\Z\to\R$ and any possibly infinite exponent $p\geq1$, we define $\|\phi\|_{\x;p}^{p}=\sum_{\x\in\Z}|\phi_{\x}|^{p}$ if $p$ is finite and $\|\phi\|_{\x;\infty}=\sup_{\x\in\Z}|\phi_{\x}|$ otherwise/if $p=\infty$.
\end{definition}
%%%
We will also take for granted the following PDE satisfied by the $\mathbf{H}^{\beta_{\star}}$ heat kernel; it is the classical \emph{forward Kolmogorov equation}, whose proof follows by the usual integration-by-parts/summation-by-parts calculation, which is equivalent to the elementary fact about Markov processes that the backwards evolution has generator given by the adjoint of the generator for the forwards evolution:
\begin{align}
\partial_{\t}\mathbf{H}_{\s,\t,\x,\y}^{\beta_{\star}} \ = \ 2^{-1}N^{2}\Delta_{\y}(\mathbf{H}_{\s,\t,\x,\y}^{\beta_{\star}}\mathbf{A}^{N,\beta_{\star}}_{\y}). \label{eq:kfe}
\end{align}
%
%%%
\subsection{Nash Inequalities}
%%%
We start with the following classical Nash inequality, which is a discrete version of the Nash inequality on the line. Its proof is a discretization of the dynamical proof of the continuum Nash inequality.
%%%
\begin{lemma}\label{lemma:hke2}
\fsp For any function $\phi:\Z\to\R$, we have $\|\phi\|_{\x;2}^{2}\lesssim\|\phi\|_{\x;1}^{4/3}\|\grad_{1}^{\mathbf{X}}\phi\|_{\x;2}^{2/3}$, where $\grad_{\mathfrak{l}}^{\mathbf{X}}\phi_{\x}=\phi_{\x+\mathfrak{l}}-\phi_{\x}$ for any $\mathfrak{l}\in\Z$.
\end{lemma}
%%%
%%%
\begin{proof}
For notational convenience, provided any function $\phi:\Z\to\R$, we define $\phi_{\t,\x}=\mathbf{H}^{\mathbf{X}}_{\t,\x}(\phi)$ as the action of the $\beta_{\star}=0$ heat flow on $\phi$. Differentiation with respect to $\t$ and the observation that $\Delta=-\grad_{1}^{\mathbf{X}}(\grad_{1}^{\mathbf{X}})^{\ast}$, with the superscript $\ast$ denoting adjoint with respect to $\ell^{2}(\Z)$ with constant or flat measure, provides the following inequality. Indeed, the operator identity $\Delta=-\grad_{1}^{\mathbf{X}}(\grad_{1}^{\mathbf{X}})^{\ast}$ implies that the lower bound in \eqref{eq:hke21} below would be an identity if we replace $|\grad_{1}^{\mathbf{X}}\phi_{\s,\x}|^{2}$ by the product $\grad_{1}^{\mathbf{X}}\phi_{\s,\x}\grad_{1}^{\mathbf{X}}\phi_{\x}$. The sum of this product over $\x\in\Z$ is controlled by the $\ell^{2}(\Z)$-norms of each factor by the Cauchy-Schwarz inequality. Thus, we have
\begin{align}
{\sum}_{\x\in\Z}\phi_{\t,\x}\phi_{\x} \ &= \ {\sum}_{\x\in\Z}|\phi_{\x}|^{2} + 2^{-1}N^{2}\int_{0}^{\t}{\sum}_{\x\in\Z}\Delta_{\x}\phi_{\s,\x}\phi_{\x}\d\s \\
&\geq \ {\sum}_{\x\in\Z}|\phi_{\x}|^{2}-2^{-1}N^{2}\int_{0}^{\t}\left({\sum}_{\x\in\Z}|\grad_{1}^{\mathbf{X}}\phi_{\s,\x}|^{2}\right)^{1/2}\left({\sum}_{\x\in\Z}|\grad_{1}^{\mathbf{X}}\phi_{\x}|^{2}\right)^{1/2}\d\s. \label{eq:hke21}
\end{align}
Observe that the first sum over $\x\in\Z$ appearing in \eqref{eq:hke21} is non-increasing in the time-variable $\s$. To see this, note that $\grad_{1}^{\mathbf{X}}\phi_{\s,\x}$ is the same as $(\grad_{1}^{\mathbf{X}}\phi)_{\s,\x}$, namely the action of the heat operator on $\grad_{1}^{\mathbf{X}}\phi$; this is because gradients commute with the Laplacian $\Delta$ and therefore with the heat operator $\mathbf{H}^{\mathbf{X}}$ as well since the heat operator is an analytic function/exponential of $\Delta$. Thus, we deduce that the sum over $\x\in\Z$ on the far RHS of \eqref{eq:hke21} is therefore the squared $\ell^{2}(\Z)$-norm of a heat operator acting on $\grad_{1}^{\mathbf{X}}\phi$, from which the non-increasing property follows by standard convexity of the heat operator $\mathbf{H}^{\mathbf{X}}$, or equivalently the Cauchy-Schwarz inequality upon realizing $\mathbf{H}^{\mathbf{X}}$ as an expectation with respect to the law of a simple random walk. We therefore deduce from \eqref{eq:hke21} that
\begin{align}
\|\phi\|_{\x;2}^{2} \ \leq \ \|\phi_{\t,\x}\|_{\x;\infty}\|\phi\|_{\x;1} + 2^{-1}N^{2}\t\|\grad_{1}^{\mathbf{X}}\phi\|_{\x;2}^{2} \ \lesssim \ \sup_{\x,\y\in\Z}|\mathbf{H}_{0,\t,\x,\y}|\|\phi\|_{\x;1}^{2} + 2^{-1}N^{2}\t\|\grad_{1}^{\mathbf{X}}\phi\|_{\x;2}^{2}. \label{eq:hke22}
\end{align}
We may now employ Proposition \ref{prop:hke3} with $\beta_{\star}=0$, which is proven in Proposition A.1 in \cite{DT} but we cite an upcoming result in this paper for the reader's convenience, for $\s=0$ and optimize in $\t$ to finish the proof.
\end{proof}
%%%
Using Lemma \ref{lemma:hke2} along with a twist of the usual Nash approach for global heat kernel estimates, we establish the following estimate in which $5/2$ could be replaced by any large constant if we are willing to take $\beta_{\star}$ to be smaller by some fixed positive factor in Theorem \ref{theorem:kpz}; we provide $5/2$ because it is what falls out from the proof of the following, in case this is of any interest.
%%%
\begin{prop}\label{prop:hke3}
\fsp For any $\s\leq\t$ and $\x,\y\in\Z$, we have the uniform estimates $0\leq\mathbf{H}^{\beta_{\star}}_{\s,\t,\x,\y}\lesssim N^{-1+\frac52\beta_{\star}}|\t-\s|^{-1/2}$.
\end{prop}
%%%
%%%
\begin{remark}\label{remark:hke4}
\fsp Let us note $\mathbf{H}^{\beta_{\star}}$ is always non-negative and at most 1 if ${\mathrm{s}}\leq\t$. This can be proven either by realizing it as a transition probability for a random walk or by maximum principle. We will often use this elementary fact without explicit proof or reference. Let us also note that if $\beta_{\star}=0$, then Proposition \ref{prop:hke3} recovers for us the classical on-diagonal pointwise upper bound for the classical semi-discrete heat kernel; see Appendix in \cite{DT}. In particular, the statement $\mathbf{H}^{\beta_{\star}}\approx\mathbf{H}$ is not directly suggested by Proposition \ref{prop:hke3} because of the $\beta_{\star}$-dependent power of $N$. We emphasize, again, this is because Proposition \ref{prop:hke3} is \emph{global}, whereas $\mathbf{H}^{\beta_{\star}}\approx\mathbf{H}$ is an approximation that holds in the large-$N$ limit in a slightly weaker sense; see Proposition \ref{prop:hke7} and Proposition \ref{prop:hke11}.
\end{remark}
%%%
%%%
\begin{proof}
We will assume $\s=0$; indeed, the heat kernel $\mathbf{H}^{\beta_{\star}}$ is time-homogeneous as the parabolic operator defining it does not have time-dependence in its coefficients. We first propose the following identity, which is known as the Chapman-Kolmogorov equation or the semigroup property, that we will use throughout this paper:
\begin{align}
\mathbf{H}_{0,\t,\x,\y}^{\beta_{\star}} = {\sum}_{\z\in\Z}\mathbf{H}_{2^{-1}\t,\t,\x,\z}^{\beta_{\star}}\mathbf{H}^{\beta_{\star}}_{0,2^{-1}\t,\z,\y}. \label{eq:cke}
\end{align}
Indeed, \eqref{eq:cke} follows by conditioning on the random walk position, for the simple random walk whose transition probabilities are given by $\mathbf{H}^{\beta_{\star}}$, at time $\t/2$. We now apply the Cauchy-Schwarz inequality to get the following estimate in which we recall the slow bond coefficient $\mathbf{A}^{N,\beta_{\star}}$ from Proposition \ref{prop:mshe1}:
\begin{align}
|\mathbf{H}_{0,\t,\x,\y}^{\beta_{\star}}|^{2} \ \leq \ \left({\sum}_{\z\in\Z}|\mathbf{H}_{2^{-1}\t,\t,\x,\z}^{\beta_{\star}}|^{2}\mathbf{A}_{\z}^{N,\beta_{\star}}\right)\left({\sum}_{\z\in\Z}|\mathbf{H}_{0,2^{-1}\t,\z,\y}^{\beta_{\star}}|^{2}(\mathbf{A}_{\z}^{N,\beta_{\star}})^{-1}\right) \ \overset{\bullet}= \ \mathbf{H}(\t,\x;1,2^{-1}\t)\mathbf{H}(2^{-1}\t,\y;2). \label{eq:hke31}
\end{align}
We will now apply the usual Nash argument for on-diagonal heat kernel estimates. First, we will differentiate the first term on the RHS of \eqref{eq:hke31} with respect to $\t$, treating the backwards time-variable $2^{-1}\t$ as fixed and independent of this differentiation. For this, we will employ the Kolmogorov forward PDE \eqref{eq:kfe}; this ultimately gives, upon defining $\mathbf{HA}_{\t,\x,\z}^{N,\beta_{\star}}=\mathbf{H}_{2^{-1}\t,\t,\x,\z}^{\beta_{\star}}\mathbf{A}_{\z}^{N,\beta_{\star}}$, that
\begin{align}
{\partial_{\t}\mathbf{H}(\t,\x;1,\s)|_{\s=2^{-1}\t}} \ = \ N^{2}{\sum}_{\z\in\Z}\mathbf{HA}_{\t,\x,\z}^{N,\beta_{\star}}\Delta_{\z}(\mathbf{HA}_{\t,\x,\z}^{N,\beta_{\star}}) \ = \ -N^{2}{\sum}_{\z\in\Z}|\grad_{1}^{\mathbf{X}}\mathbf{HA}_{\t,\x,\z}^{N,\beta_{\star}}|^{2}, \label{eq:hke32}
\end{align}
where the gradient on the far RHS of \eqref{eq:hke32} is with respect to $\z$; the far RHS of \eqref{eq:hke32} follows from summation-by-parts similar to \eqref{eq:hke21}. We now apply the Nash inequality of Lemma \ref{lemma:hke2} to the far RHS of \eqref{eq:hke32}; this gives
\begin{align}
{\partial_{\t}\mathbf{H}(\t,\x;1,\s)|_{\s=2^{-1}\t}} \ \lesssim \ -N^{2}\|\mathbf{HA}_{\t,\x,\z}^{N,\beta_{\star}}\|_{\z;2}^{6}\|\mathbf{HA}_{\t,\x,\z}^{N,\beta_{\star}}\|_{\z;1}^{-4}. \label{eq:hke33}
\end{align}
We now control both factors on the RHS of \eqref{eq:hke33}. To this end, let us first observe that $1\geq\mathbf{A}^{N,\beta_{\star}}\geq N^{-\beta_{\star}}$ uniformly by construction in Proposition \ref{prop:mshe1}. Therefore, we may replace $\mathbf{HA}^{N,\beta_{\star}}$ with $\mathbf{H}^{\beta_{\star}}(\mathbf{A}^{N,\beta_{\star}})^{1/2}$ in the first norm on the RHS of \eqref{eq:hke33} if we are willing to give up a factor of $N^{-3\beta_{\star}}$ and in the second norm on the RHS of \eqref{eq:hke33} for free because of the negative exponent for this second norm. Because the heat kernel $\mathbf{H}^{\beta_{\star}}$ is a probability measure in its forward spatial variable, we deduce
\begin{align}
\|\mathbf{HA}_{\t,\x,\z}^{N,\beta_{\star}}\|_{\z;2}^{6}\|\mathbf{HA}_{\t,\x,\z}^{N,\beta_{\star}}\|_{\z;1}^{-4} \ {\gtrsim} \ N^{-3\beta_{\star}}\|\mathbf{H}_{2^{-1}\t,\t,\x,\z}^{\beta_{\star}}(\mathbf{A}^{N,\beta_{\star}}_{\z})^{1/2}\|_{\z;2}^{6}\|\mathbf{H}_{2^{-1}\t,\t,\x,\z}^{\beta_{\star}}\|_{\z;1}^{-4} \ = \ N^{-3\beta_{\star}}\mathbf{H}(\t,\x;1,2^{-1}\t)^{3}. \label{eq:hke34}
\end{align}
For clarity, we now summarize the discussion after \eqref{eq:hke31}; from \eqref{eq:hke32}, \eqref{eq:hke33}, and \eqref{eq:hke34}, we deduce
\begin{align}
{\partial_{\t}\mathbf{H}(\t,\x;1,\s)|_{\s=2^{-1}\t}} \ \lesssim \ {-}N^{2}N^{-3\beta_{\star}}\mathbf{H}(\t,\x;1,2^{-1}\t)^{3}, \label{eq:hke35}
\end{align}
from which, by the comparison principle for ODEs and explicitly solving the equation $\partial_{\t}\phi(\t)=-\alpha\phi(\t)^{3}$, we deduce 
\begin{align}
\mathbf{H}(\t,\x;1,2^{-1}\t) \ \lesssim \ N^{-1+3\beta_{\star}/2}\t^{-1/2}. \label{eq:hke36}
\end{align}
We will now estimate the second factor $\mathbf{H}(2^{-1}\t,\y;2)$ in \eqref{eq:hke31}. This will be done using a similar procedure. First, we differentiate $\mathbf{H}(2^{-1}\t,\y;2)$ with respect to $2^{-1}\t$ by using the PDE for the $\mathbf{H}^{\beta_{\star}}$ heat kernel from Definition \ref{definition:mshe2}. We then apply the Nash inequality in Lemma \ref{lemma:hke2} and a summation-by-parts as we did for \eqref{eq:hke32} and ultimately deduce
\begin{align}
\partial_{2^{-1}\t}\mathbf{H}(2^{-1}\t,\y;2) \ = \ N^{2}{\sum}_{\z\in\Z}\mathbf{H}_{0,2^{-1}\t,\z,\y}^{\beta_{\star}}\Delta_{\z}\mathbf{H}_{0,2^{-1}\t,\z,\y}^{\beta_{\star}} \ \lesssim \ -N^{2}\|\mathbf{H}^{\beta_{\star}}_{0,2^{-1}\t,\z,\y}\|_{\z;2}^{6}
\|\mathbf{H}^{\beta_{\star}}_{0,2^{-1}\t,\z,\y}\|_{\z;1}^{-4}. \label{eq:hke37}
\end{align}
We note the square of the first norm on the far RHS of \eqref{eq:hke37} is equal to $\mathbf{H}(2^{-1}\t,\y;2)$ itself if we replace $\mathbf{H}^{\beta_{\star}}$ by $\mathbf{H}^{\beta_{\star}}(\mathbf{A}^{N,\beta})^{-1/2}$. Making such a replacement costs a factor of $N^{-3\beta_{\star}}$ after we take sixth-powers because $\mathbf{A}^{N,\beta_{\star}}\geq N^{-\beta_{\star}}$. We also note that the last $\|\|_{\z;1}$-norm of the heat kernel $\mathbf{H}^{\beta_{\star}}$ is not just equal to 1 as the heat kernel is not necessarily a probability measure with respect to its backwards spatial variable; the $\|\|_{\z;1}$-norm of the heat kernel is not conserved. However, we do have
\begin{align}
\|\mathbf{H}^{\beta_{\star}}_{0,2^{-1}\t,\z,\y}\|_{\z;1} \ \leq \ \|\mathbf{H}^{\beta_{\star}}_{0,2^{-1}\t,\z,\y}(\mathbf{A}^{N,\beta_{\star}}_{\z})^{-1}\|_{\z;1} \ = \ \|\mathbf{H}^{\beta_{\star}}_{0,0,\z,\y}(\mathbf{A}^{N,\beta_{\star}}_{\z})^{-1}\|_{\z;1} \ \leq \ N^{\beta_{\star}}. \label{eq:hke38}
\end{align}
The middle identity in \eqref{eq:hke38} follows from differentiating the sum of what is inside the second $\|\|_{\z;1}$-norm in $\t$; the $(\mathbf{A}^{N,\beta})^{-1}$ factor cancels with the $\mathbf{A}^{N,\beta}$ factor in front of the Laplacian defining the operator $\Delta^{N,\beta_{\star}}$ that appears when we employ the PDE for $\mathbf{H}^{\beta_{\star}}$ in Definition \ref{definition:mshe2}. The final estimate in \eqref{eq:hke38} follows by noting $\mathbf{H}^{\beta_{\star}}$ is a probability measure with respect to both spatial variables when the time-variables are equal by construction in Definition \ref{definition:mshe2}, as well as the lower bound $\mathbf{A}^{N,\beta_{\star}}\geq N^{-\beta_{\star}}$. Combining \eqref{eq:hke38} with \eqref{eq:hke37} and the paragraph in between these two estimates, we deduce
\begin{align}
\partial_{2^{-1}\t}\mathbf{H}(2^{-1}\t,\y;2) \ \lesssim \ -N^{2-7\beta_{\star}}\partial_{2^{-1}\t}\mathbf{H}(2^{-1}\t,\y;2)^{3}, \label{eq:hke39}
\end{align}
from which we deduce the following similar to how we deduced \eqref{eq:hke36}:
\begin{align}
\mathbf{H}(2^{-1}\t,\y;2) \ \lesssim \ N^{-1+7\beta_{\star}/2}\t^{-1/2}. \label{eq:hke39}
\end{align}
Combining \eqref{eq:hke36} and \eqref{eq:hke39} with the initial estimate \eqref{eq:hke31} completes the proof.
\end{proof}
%%%
%%%
\subsection{Perturbative Estimates}
%%%
The backbone behind obtaining sharper local estimates for $\mathbf{H}^{\beta_{\star}}$ is the following perturbative identity, which allows us to borrow $\beta_{\star}=0$ estimates in Appendix of \cite{DT}. The proof is just the classical Duhamel formula.
%%%
\begin{lemma}\label{lemma:hke5}
\fsp We have $\mathbf{H}^{\beta_{\star}}_{\s,\t,\x,\y}=\mathbf{H}_{\s,\t,\x,\y}-N^{2}\bar{\mathbf{A}}^{N,\beta_{\star}}\int_{\s}^{\t}\mathbf{H}_{\r,\t,\x,0}^{\beta_{\star}}\Delta\mathbf{H}_{\s,\r,0,\y}\d\r$, where $\bar{\mathbf{A}}^{N}=2^{-1}-2^{-1}N^{-\beta_{\star}}$.
\end{lemma}
%%%
%%%
\begin{proof}
Differentiate $\sum_{\z\in\Z}\mathbf{H}_{\r,\t,\x,\z}^{\beta_{\star}}\mathbf{H}_{\s,\r,\z,\y}$ in $\r$, use the fundamental theorem of calculus, and use the defining PDEs in Definition \ref{definition:mshe2}; because the heat kernels in Definition \ref{definition:mshe2} are time-homogeneous, differentiating in the backwards time-variable $\s$ is the same as differentiating in the forwards time-variable $\t$ but in the negative direction/with an additional negative sign.
\end{proof}
%%%
Lemma \ref{lemma:hke5} will be useful for obtaining estimates for $\mathbf{H}^{\beta_{\star}}$ in terms of those for $\mathbf{H}$; this leads naturally to consider the following.
%%%
\begin{definition}\label{definition:hke6}
\fsp For notational convenience, we define $\wt{\mathbf{H}}^{\beta_{\star}}=\mathbf{H}^{\beta_{\star}}-\mathbf{H}$ as functions in space-time.
\end{definition}
%%%
Let us clarify before actually stating the following result that the estimate therein effectively confirms $\mathbf{H}^{\beta_{\star}}\approx\mathbf{H}$ in the following weaker and non-global sense. Away from the slow bond, the following result says $\mathbf{H}^{\beta_{\star}}-\mathbf{H}$ behaves lower-order in $N$ than what we expect either heat kernel to behave like in terms of $N$-dependence; the time-dependent factor is worse for times less than 1, but it is integrable at $\s\approx\t$, which is all we need. In particular, if $|\t-\s|$ is order 1, then summing the following estimate for $\mathbf{H}^{\beta_{\star}}-\mathbf{H}$ over a macroscopic-size set whose length is order $N$ gives something that vanishes quantitatively in the large-$N$ limit. Briefly, the proof of the following Proposition \ref{prop:hke7} will be a fixed-point argument based on the perturbative formula in Lemma \ref{lemma:hke5} as in classical parabolic theory, but it also requires an additional input from Proposition \ref{prop:hke3}.
%%%
\begin{prop}\label{prop:hke7}
\fsp There exist uniformly positive $\e_{1},\e_{2}\gtrsim\beta_{\star}$ such that $|\wt{\mathbf{H}}^{\beta_{\star}}_{\s,\t,\x,\y}|\lesssim N^{-1-\e_{1}}|\t-\s|^{-1+\e_{2}}$ for all $\t\geq\s$ satisfying $\s,\t\leq10$ and for all $|\y|\geq N^{999\beta_{\star}}$; the constant $999$ can be replaced by any large constant if we change $\e_{1},\e_{2}$ by positive factors.
\end{prop}
%%%
%%%
\begin{proof}
Define a ``regularizing time" $\t^{\mathrm{cut}}=N^{-2+99\beta_{\star}}$ and $\mathbb{I}^{\beta_{\star}}$ to be the set of $|\y|\geq N^{999\beta_{\star}}$. We also define, for $\alpha$ positive,
\begin{align}
\Upsilon_{\s,\t}^{\beta_{\star}} \ = \ |\t-\s|^{1-\alpha}\sup_{\x\in\Z}\sup_{\y\in\mathbb{I}^{\beta_{\star}}}|\mathbf{H}^{\beta_{\star}}_{\s,\t,\x,\y}| \and \Gamma_{\s,\t}^{\beta_{\star}} \ = \ |\t-\s|^{1-\alpha}\sup_{\x\in\Z}\sup_{\y\in\mathbb{I}^{\beta_{\star}}}|\wt{\mathbf{H}}_{\s,\t,\x,\y}^{\beta_{\star}}|. 
\end{align}
We will first estimate $\Upsilon^{\beta_{\star}}$ and then afterwards $\Gamma^{\beta_{\star}}$. To analyze the former, let us consider first the following decomposition, which takes the perturbative identity in Lemma \ref{lemma:hke5} and cuts off the time-integral therein with respect to $\t^{\mathrm{cut}}$. To be precise, if $\s+\t^{\mathrm{cut}}\geq\t$, then last integral on the RHS of \eqref{eq:hke71} below is omitted entirely; we also recall notation in Lemma \ref{lemma:hke5}:
\begin{align}
\mathbf{H}_{\s,\t,\x,\y}^{\beta_{\star}} \ = \ \mathbf{H}_{\s,\t,\x,\y} - N^{2}\bar{\mathbf{A}}^{N,\beta_{\star}}\int_{\s+\t^{\mathrm{cut}}}^{\t}\mathbf{H}_{\r,\t,\x,0}^{\beta_{\star}}\Delta\mathbf{H}_{\s,\r,0,\y}\d\r - N^{2}\bar{\mathbf{A}}^{N,\beta_{\star}}\int_{\s}^{\s+\t^{\mathrm{cut}}}\mathbf{H}_{\r,\t,\x,0}^{\beta_{\star}}\Delta\mathbf{H}_{\s,\r,0,\y}\d\r. \label{eq:hke71}
\end{align}
Let us define the first integral on the RHS of \eqref{eq:hke71} as $\mathbf{H}^{\beta_{\star}}_{\s,\t,\x,\y}(1)$ and the second by $\mathbf{H}^{\beta_{\star}}_{\s,\t,\x,\y}(2)$. We first control $\mathbf{H}^{\beta_{\star}}_{\s,\t,\x,\y}(1)$ via
\begin{align}
\mathbf{H}^{\beta_{\star}}_{\s,\t,\x,\y}(1) \ &= \ \int_{\s+\t^{\mathrm{cut}}}^{\t}{\sum}_{\w\in\mathbb{I}^{\beta_{\star}}}\mathbf{H}_{\frac{\t+\r}{2},\t,\x,\w}^{\beta_{\star}}\mathbf{H}_{\r,\frac{\t+\r}{2},\w,0}^{\beta_{\star}}\Delta\mathbf{H}_{\s,\r,0,\y}\d\r +\int_{\s+\t^{\mathrm{cut}}}^{\t}{\sum}_{\w\not\in\mathbb{I}^{\beta_{\star}}}\mathbf{H}_{\frac{\t+\r}{2},\t,\x,\w}^{\beta_{\star}}\mathbf{H}_{\r,\frac{\t+\r}{2},\w,0}^{\beta_{\star}}\Delta\mathbf{H}_{\s,\r,0,\y}\d\r \nonumber \\
&\overset{\bullet}= \ \mathbf{H}^{\beta_{\star}}_{\s,\t,\x,\y}(1,1)+\mathbf{H}^{\beta_{\star}}_{\s,\t,\x,\y}(1,2). \label{eq:hke72} 
\end{align}
Let us make a clarifying remark; the $(1,1)$-term in \eqref{eq:hke72} will be analyzed via a fixed-point-type estimate that is linear in $\Upsilon^{\beta_{\star}}$, and the $(1,2)$-term will be analyzed as an error term via the a priori global estimate in Proposition \ref{prop:hke3} combined with the smallness of the complement of $\mathbb{I}^{\beta_{\star}}$, which is relevant in said $(1,2)$-term. The second integral on the RHS of \eqref{eq:hke71} will also be treated similarly as an error, namely by Proposition \ref{prop:hke3}, but instead of smallness of the complement of $\mathbb{I}^{\beta_{\star}}$, we use smallness of $\mathfrak{t}^{\mathrm{cut}}$ and eventually take $\y\in\mathbb{I}^{\beta_{\star}}$ in order to control $\Upsilon^{\beta_{\star}}$, which we emphasize only looks at forward variables in $\mathbb{I}^{\beta_{\star}}$.

Let us now recall the definition of $\Upsilon^{\beta_{\star}}$ given at the beginning of this proof. We also recall from the proof of Proposition \ref{prop:hke3} that the sum over $\Z$ of $\mathbf{H}^{\beta_{\star}}$ in the backwards spatial variable is bounded by $N^{\beta_{\star}}$; we used this to get \eqref{eq:hke38}, for example. This will be used to estimate the $(1,1)$-term in \eqref{eq:hke72} in the following fashion. For the sum over $\w\in\mathbb{I}^{\beta_{\star}}$, we control the first $\mathbf{H}^{\beta_{\star}}$ factor in the $(1,1)$-term by $\Upsilon^{\beta_{\star}}$ if we give up a factor of $|\t-\r|^{-1+\alpha}$ times some constants; this is because we have restricted to $\w\in\mathbb{I}^{\beta_{\star}}$ in the sum. The sum over $\w\in\mathbb{I}^{\beta_{\star}}$ is then left with just the second $\mathbf{H}^{\beta_{\star}}$ factor; we control this sum from the \eqref{eq:hke38} estimate mentioned in this paragraph. We are then left to estimate the Laplacian of $\mathbf{H}$, namely for $\beta_{\star}=0$, via Proposition A.1 in \cite{DT}. We deduce
\begin{align}
|\mathbf{H}_{\s,\t,\x,\y}^{\beta_{\star}}(1,1)| \ \leq \ N^{\beta_{\star}}\int_{\s+\t^{\mathrm{cut}}}^{\t}\Upsilon_{\frac{\t+\r}{2},\t}^{\beta_{\star}}|\t-\r|^{-1+\alpha}|\Delta\mathbf{H}_{\s,\r,0,\y}|\d\r \ \lesssim \ N^{-3+\beta_{\star}}\int_{\s+\t^{\mathrm{cut}}}^{\t}\Upsilon_{\frac{\t+\r}{2},\t}^{\beta_{\star}}|\t-\r|^{-1+\alpha}|\r-\s|^{-\frac32}\d\r. \label{eq:hke73}
\end{align}
As for the $(1,2)$-term in \eqref{eq:hke72}, we note the product of two copies of $\mathbf{H}^{\beta_{\star}}$ in the $(1,2)$-term is bounded above by $N^{-2+5\beta_{\star}}|\t-\r|^{-1}$ times uniformly bounded constants and by 1, courtesy of Proposition \ref{prop:hke3} and Remark \ref{remark:hke4}, respectively. Interpolating these bounds and employing the second-order gradient estimate for $\mathbf{H}$ in Proposition A.1 in \cite{DT}, which looks like a Gaussian heat kernel with time scaled by $N^{2}$, we deduce the estimate below where $\alpha'$ is positive and satisfies $\alpha'\lesssim\alpha$, with $\alpha$ that we picked at the beginning; let us clarify that by ``interpolation", we mean changing the aforementioned $|\t-\r|^{-1}$ factor to $|\t-\r|^{-1+\alpha}$ up to giving up $N^{\alpha'}$:
\begin{align}
|\mathbf{H}_{\s,\t,\x,\y}^{\beta_{\star}}(1,2)| \ \lesssim \ N^{-5+5\beta_{\star}+\alpha'}\int_{\s+\t^{\mathrm{cut}}}^{\t}|\Z\setminus\mathbb{I}^{\beta_{\star}}||\t-\r|^{-1+\alpha}|\r-\s|^{-\frac32}\d\r \ \lesssim \ N^{-5+9999\beta_{\star}+\alpha'}(\t^{\mathrm{cut}})^{-\frac12}|\t-\s|^{-1+\alpha}. \label{eq:hke74}
\end{align}
The last estimate in \eqref{eq:hke74} follows by recalling $\Z\setminus\mathbb{I}^{\beta_{\star}}$ is a discrete neighborhood of the origin of length at most $3N^{999\beta_{\star}}$ combined with a fairly elementary integral estimate; for this integral estimate, we note the non-integrable singularity $|\r-\s|^{-3/2}$ is regularized by looking only past time $\s+\t^{\mathrm{cut}}$, and because $|\r-\s|^{-3/2}$ is only a square root away from being basically integrable, we get the $-1/2$-power of the regularizing time $\t^{\mathrm{cut}}$. Recalling $\t^{\mathrm{cut}}=N^{-2+99\beta_{\star}}$, we deduce from \eqref{eq:hke74} the bound
\begin{align}
|\mathbf{H}_{\s,\t,\x,\y}^{\beta_{\star}}(1,2)| \ \lesssim \ N^{-4+9999\beta_{\star}+\alpha'}|\t-\s|^{-1+\alpha}. \label{eq:hke75}
\end{align}
The estimates \eqref{eq:hke73} and \eqref{eq:hke75} control the first integral on the RHS of \eqref{eq:hke71}. To control the second integral therein, we will first restrict to $|\y|\geq N^{999\beta_{\star}}$, which is okay if we are controlling $\Upsilon^{\beta_{\star}}$, since this quantity restricts to such forward variables $\y$. Because $|\y|\geq N^{999\beta_{\star}}$, the off-diagonal bound for the $\mathbf{H}$ heat kernel in Proposition A.1 in \cite{DT} implies the $\Delta\mathbf{H}$-term in this second integral of \eqref{eq:hke71} is exponentially small in $N$. Indeed, for times shorter than $\t^{\mathrm{cut}}=N^{-2+99\beta_{\star}}$, the probability that the simple random walk described by $\mathbf{H}$ propagates a distance at least $N^{999\beta_{\star}}$ is exponentially small in $N$ by standard concentration estimates for random walks/martingales. Therefore, because all other terms in the second integral in \eqref{eq:hke71} are uniformly bounded per Remark \ref{remark:hke4}, and because exponentially-small factors in $N$ beat any fixed positive power of $N$, for any positive $D$ we deduce
\begin{align}
|N^{2}\bar{\mathbf{A}}^{N,\beta_{\star}}\mathbf{H}_{\s,\t,\x,\y}^{\beta_{\star}}(2)| \ \lesssim_{D} \ N^{-D}. \label{eq:hke76}
\end{align}
Let us now combine \eqref{eq:hke71}, \eqref{eq:hke72}, \eqref{eq:hke73}, \eqref{eq:hke75}, and \eqref{eq:hke76} together with the on-diagonal estimate for $\mathbf{H}$ from Proposition A.1 in \cite{DT}, or equivalently Proposition \ref{prop:hke3} for $\beta_{\star}=0$; this provides, because $\bar{\mathbf{A}}^{N,\beta_{\star}}$ is uniformly bounded as constructed in Lemma \ref{lemma:hke5}, the following estimate, which we emphasize is uniform in all space-time variables $(\s,\t,\x,\y)$ with $\s,\t\leq10$.
\begin{align}
|\mathbf{H}_{\s,\t,\x,\y}^{\beta_{\star}}||\t-\s|^{1-\alpha} \ \lesssim \ N^{-1} + N^{-1+\beta_{\star}}|\t-\s|^{1-\alpha}\int_{\s+\t^{\mathrm{cut}}}^{\t}\Upsilon_{\frac{\t+\r}{2},\t}^{\beta_{\star}}|\t-\r|^{-1+\alpha}|\r-\s|^{-\frac32}\d\r + N^{-2+9999\beta_{\star}+\alpha'}. \label{eq:hke77}
\end{align}
We observe that the RHS of \eqref{eq:hke77} is uniform in allowable space-time coordinates $(\s,\t,\x,\y)$ with $\y\in\mathbb{I}^{\beta_{\star}}$, so \eqref{eq:hke77} extends to an estimate on $\Upsilon^{\beta_{\star}}$. From \eqref{eq:hke77} combined with the Gronwall inequality, which is still applicable because the non-$\Upsilon^{\beta_{\star}}$ factors in the integral on the RHS of \eqref{eq:hke77} are integrable on $[\s+\t^{\mathrm{cut}},\t]$ given the additional regularization $\t^{\mathrm{cut}}$, we deduce that, if $\beta_{\star}\leq999^{-1}$, for example and $\alpha$ is sufficiently small but fixed, we have the Gronwall-type upper bound for some fixed positive $\kappa$:
\begin{align}
\sup_{\x\in\Z}\sup_{{\y}\in\mathbb{I}^{\beta_{\star}}}|\mathbf{H}_{\s,\t,\x,\y}^{\beta_{\star}}||\t-\s|^{1-\alpha} \ = \ \Upsilon_{\s,\t}^{\beta_{\star}} \ \lesssim \ N^{-1}\exp\left(\kappa N^{-1+\beta_{\star}}|\t-\s|^{1-\alpha}\int_{\s+\t^{\mathrm{cut}}}^{\t}|\t-\r|^{-1+\alpha}|\r-\s|^{-\frac32}\d\r\right) \ \lesssim \ N^{-1}, \nonumber
\end{align}
where the last estimate follows by integrating. We now estimate $\Gamma^{\beta_{\star}}$. For this, by \eqref{eq:hke71}, \eqref{eq:hke72}, \eqref{eq:hke73}, \eqref{eq:hke75}, and \eqref{eq:hke76}, we get
\begin{align}
|\wt{\mathbf{H}}_{\s,\t,\x,\y}^{\beta_{\star}}||\t-\s|^{1-\alpha} \ \lesssim \ N^{-1+\beta_{\star}}|\t-\s|^{1-\alpha}\int_{\s+\t^{\mathrm{cut}}}^{\t}\Upsilon_{\frac{\t+\r}{2},\t}^{\beta_{\star}}|\t-\r|^{-1+\alpha}|\r-\s|^{-\frac32}\d\r + N^{-2+9999\beta_{\star}+\alpha'}. \label{eq:hke78}
\end{align}
Using the bound $\Upsilon^{\beta_{\star}}\lesssim N^{-1}$ for times before $10$ that we obtained via \eqref{eq:hke77} in the paragraph prior to \eqref{eq:hke78} and then a straightforward integration of the RHS of \eqref{eq:hke78} then implies $\Gamma^{\beta_{\star}}_{\s,\t}\lesssim N^{-1-\e_{1}}$ for all $\s\leq\t\leq10$ for some positive $\e_{1}$. Choosing $\alpha=\e_{2}$ then completes the proof of the proposed $\wt{\mathbf{H}}^{\beta_{\star}}$ estimate after rearranging the $|\t-\s|^{-1+\alpha}=|\t-\s|^{-1+\e_{2}}$ factor.
\end{proof}
%%%
%%%
\subsection{Regularity Estimates}
%%%
We could obtain space-time Holder regularity for $\mathbf{H}^{\beta_{\star}}$ through the Nash method, but this gives quite sub-optimal regularity and would not provide the necessary explicit and quantitative space-time gradient estimates. Our strategy is again by perturbing in $\beta_{\star}$ via Lemma \ref{lemma:hke5} and borrowing regularity estimates for $\beta_{\star}=0$ obtained in Appendix 1 of \cite{DT}. Let us note the following regularity estimates for $\mathbf{H}^{\beta_{\star}}$ are rather opaque, but roughly speaking they may be interpreted as a regularity estimate for $\beta_{\star}=0$ in Appendix 1 of \cite{DT} times small powers of $N$ plus an error term, independent of the length-scale or time-scale of the gradient. Note the estimates in Proposition \ref{prop:hke8} are uniform in space-time; this is unlike Proposition \ref{prop:hke7}.
%%%
\begin{prop}\label{prop:hke8}
\fsp Given any $\mathfrak{l}\in\Z$, we have the following estimate uniformly in $\t\geq\s$ satisfying $\s,\t\leq10$ and $\x,\y\in\Z$ in which $\beta_{\star}\lesssim\e_{1},\e_{2}\lesssim\beta_{\star}$ are uniformly positive in $N$:
\begin{align}
|\grad_{\mathfrak{l}}^{\mathbf{X}}\mathbf{H}_{\s,\t,\x,\y}^{\beta_{\star}}| \ \lesssim_{\e_{1},\e_{2}} \ N^{-2+2\beta_{\star}}|\t-\s|^{-1+\e_{1}}|\mathfrak{l}| + N^{-2+99999\beta_{\star}}|\t-\s|^{-1+\e_{2}}. \label{eq:hke8I}
\end{align}
For any positive $\t_{N}$, we have the following time-gradient bound in which $\e_{3},\e_{4}\gtrsim\beta_{\star}$ are uniformly, in $N$, positive, and in which $\grad^{\mathbf{T}}$ below acts on the $\t$-variable; here, given any $\phi:\R_{\geq0}\times\Z\to\R$ and time-scale $\r\in\R$, we define $\grad^{\mathbf{T}}_{\r}\phi_{\t,\x}=\phi_{\t+\r,\x}-\phi_{\t,\x}$ if $\t+\r\geq0$ and  $\grad^{\mathbf{T}}_{\r}\phi_{\t,\x}=\phi_{0,\x}-\phi_{\t,\x}$ if $\t+\r\leq0$, since we only have $\phi$-values at non-negative times:
\begin{align}
|\grad_{\t_{N}}^{\mathbf{T}}\mathbf{H}_{\s,\t,\x,\y}^{\beta_{\star}}| \ \lesssim_{\e_{3},\e_{4}} \ N^{-1+5\beta_{\star}}|\t-\s|^{-1+\e_{3}}\t_{N}^{1/2}+N^{-2+99999\beta_{\star}}|\t-\s|^{-1+\e_{4}}. \label{eq:hke8II}
\end{align}
\end{prop}
%%%
%%%
\begin{proof}
Let us first prove the spatial gradient estimate \eqref{eq:hke8I}, for which we first define $\mathbb{I}^{\beta_{\star}}$ as the set of $|\y|\geq N^{999\beta_{\star}}$ as well as the following $\Upsilon^{\beta_{\star}}$ quantity not to be confused with that in the proof of Proposition \ref{prop:hke7}; again, we take $\alpha$ an arbitrary positive number:
\begin{align}
\Upsilon_{\s,\t}^{\beta_{\star}} \ = \ |\t-\s|^{1-\alpha}\sup_{\x\in\Z}\sup_{\y\in\mathbb{I}^{\beta_{\star}}}|\grad_{\mathfrak{l}}^{\mathbf{X}}\mathbf{H}_{\s,\t,\x,\y}^{\beta_{\star}}|
\end{align}
Now, provided any $\x,\y\in\Z$, we employ the Chapman-Kolmogorov equation/semigroup property of $\mathbf{H}^{\beta_{\star}}$ and then decompose the resulting sum according to membership in $\mathbb{I}^{\beta_{\star}}$. This, along with the triangle inequality, yields the estimate below for any $\x,\y\in\Z$:
\begin{align}
|\grad_{\mathfrak{l}}^{\mathbf{X}}\mathbf{H}_{\s,\t,\x,\y}^{\beta_{\star}}| \ \leq \ {\sum}_{\w\in\mathbb{I}^{\beta_{\star}}}|\grad_{\mathfrak{l}}^{\mathbf{X}}\mathbf{H}_{\frac{\t+\s}{2},\t,\x,\w}^{\beta_{\star}}|\mathbf{H}_{\s,\frac{\t+\s}{2},\w,\y}^{\beta_{\star}} + {\sum}_{\w\not\in\mathbb{I}^{\beta_{\star}}}|\grad_{\mathfrak{l}}^{\mathbf{X}}\mathbf{H}_{\frac{\t+\s}{2},\t,\x,\w}^{\beta_{\star}}|\mathbf{H}_{\s,\frac{\t+\s}{2},\w,\y}^{\beta_{\star}}. \label{eq:hke8I1}
\end{align}
Observe that the spatial gradient of $\mathbf{H}^{\beta_{\star}}$ is controlled by $\mathbf{H}^{\beta_{\star}}$ at two different spatial coordinates and the same time coordinates. Thus, to estimate the second sum on the RHS of \eqref{eq:hke8I1}, we may employ the on-diagonal estimate for $\mathbf{H}^{\beta_{\star}}$ in Proposition \ref{prop:hke3}  and Remark \ref{remark:hke4} along with the estimate $|\Z\setminus\mathbb{I}^{\beta_{\star}}|\lesssim N^{999\beta_{\star}}$ that follows by construction at the beginning of this proof. We ultimately deduce the following, which we explain a little more afterwards and in which $\alpha'\lesssim\alpha$:
\begin{align}
|\grad_{\mathfrak{l}}^{\mathbf{X}}\mathbf{H}_{\s,\t,\x,\y}^{\beta_{\star}}| \ \leq \ {\sum}_{\w\in\mathbb{I}^{\beta_{\star}}}|\grad_{\mathfrak{l}}^{\mathbf{X}}\mathbf{H}_{\frac{\t+\s}{2},\t,\x,\w}^{\beta_{\star}}|\mathbf{H}_{\s,\frac{\t+\s}{2},\w,\y}^{\beta_{\star}} + N^{999\beta_{\star}}N^{-2+5\beta_{\star}+\alpha'}|\t-\s|^{-1+\alpha}. \label{eq:hke8I2}
\end{align}
If we employed the on-diagonal estimate for each heat kernel factor in the second sum in \eqref{eq:hke8I1}, then we would obtain a factor of $|\t-\s|^{-1}$. Via interpolation with the upper bound of 1 on the $\mathbf{H}^{\beta_{\star}}$ heat kernel, as in the proof of Proposition \ref{prop:hke7}, we trade $|\t-\s|^{-1}$ for $|\t-\s|^{-1+\alpha}$ if we give up a small power $N^{\alpha'}$. This explains \eqref{eq:hke8I2} from \eqref{eq:hke8I1} and the paragraph following \eqref{eq:hke8I1}. Recalling the definition of $\Upsilon^{\beta_{\star}}$ at the beginning of this proof along with an elementary rearrangement/collection of terms and \eqref{eq:hke38} used to estimate the 1-norm of the heat kernel in the backwards spatial variable, from \eqref{eq:hke8I2}, we get
\begin{align}
|\grad_{\mathfrak{l}}^{\mathbf{X}}\mathbf{H}_{\s,\t,\x,\y}^{\beta_{\star}}| \ \leq \ N^{\beta_{\star}}|\t-\s|^{-1+\alpha}\Upsilon_{\frac{\t+\s}{2},\t}^{\beta_{\star}} + N^{-2+9999\beta_{\star}+\alpha'}|\t-\s|^{-1+\alpha}. \label{eq:hke8I3}
\end{align}
We will now estimate $\Upsilon^{\beta_{\star}}$ using basically the same perturbative fixed-point estimates in the proof of Proposition \ref{prop:hke7}. In particular, we will first point out that we may take spatial gradients of the Duhamel formula in Lemma \ref{lemma:hke5} in the $\x$-variable to deduce spatial gradients of $\mathbf{H}^{\beta_{\star}}$ satisfy the same identity in Lemma \ref{lemma:hke5} that $\mathbf{H}^{\beta_{\star}}$ itself satisfies if we replace the first $\mathbf{H}$ term on the RHS therein with its own spatial gradient. We may now follow the proof of Proposition \ref{prop:hke7}, formally replacing $\Upsilon^{\beta_{\star}}$ therein by $\Upsilon^{\beta_{\star}}$ in this proof. Indeed, the upper bound for $\Upsilon^{\beta_{\star}}$ therein holds for $\Upsilon^{\beta_{\star}}$ herein because $\Upsilon^{\beta_{\star}}$ herein is controlled by the heat kernel at different spatial coordinates but same time-coordinates, and the only estimate on $\Upsilon^{\beta_{\star}}$ in the proof of Proposition \ref{prop:hke7} is the on-diagonal estimate of Proposition \ref{prop:hke3} and Remark \ref{remark:hke4} for the heat kernel $\mathbf{H}^{\beta_{\star}}$. Thus, we get the following analog of \eqref{eq:hke77} except the $N^{-1}$-term on the RHS of \eqref{eq:hke77}, which came from an on-diagonal bound for the $\mathbf{H}$ heat kernel, is replaced by a gradient of said heat kernel; we also recall the ```regularization time" $\t^{\mathrm{cut}}=N^{-2+99\beta_{\star}}$ from the proof of Proposition \ref{prop:hke7}:
\begin{align}
|\grad_{\mathfrak{l}}^{\mathbf{X}}\mathbf{H}_{\s,\t,\x,\y}^{\beta_{\star}}||\t-\s|^{1-\alpha} \ \lesssim \ |\grad_{\mathfrak{l}}^{\mathbf{X}}\mathbf{H}_{\s,\t,\x,\y}||\t-\s|^{1-\alpha} + N^{-1+\beta_{\star}}|\t-\s|^{1-\alpha}\int_{\s+\t^{\mathrm{cut}}}^{\t}\Upsilon_{\frac{\t+\r}{2},\t}^{\beta_{\star}}|\t-\r|^{-1+\alpha}|\r-\s|^{-\frac32}\d\r + N^{-2+9999\beta_{\star}+\alpha'}. \nonumber
\end{align}
The first term on the RHS of the previous estimate is controlled by $N^{-2+2\alpha}|\mathfrak{l}|$ courtesy of Proposition A.1 in \cite{DT}. In particular, this term can be controlled by pretending $\mathbf{H}$ is the limit Gaussian heat kernel of time scaled by $N^{2}$; this would give $N^{-2}|\mathfrak{l}||\t-\s|^{-1}$, but we may forget the singular time factor if we give up $N^{2\alpha}$ upon interpolating with $|\grad^{\mathbf{X}}\mathbf{H}|\lesssim\|\mathbf{H}\|_{\x;\infty}$ as time is scaled by $N^{2}$. So
\begin{align}
|\grad_{\mathfrak{l}}^{\mathbf{X}}\mathbf{H}_{\s,\t,\x,\y}^{\beta_{\star}}||\t-\s|^{1-\alpha} \ \lesssim \ N^{-2+2\alpha}|\mathfrak{l}| + N^{-1+\beta_{\star}}|\t-\s|^{1-\alpha}\int_{\s+\t^{\mathrm{cut}}}^{\t}\Upsilon_{\frac{\t+\r}{2},\t}^{\beta_{\star}}|\t-\r|^{-1+\alpha}|\r-\s|^{-\frac32}\d\r + N^{-2+9999\beta_{\star}+\alpha'}, \nonumber
\end{align}
from which, observing the RHS of this last display is uniform in $\x\in\Z$ and $\y\in\mathbb{I}^{\beta_{\star}}$, the RHS of this last estimate is also an estimate for $\Upsilon^{\beta_{\star}}$, and thus we obtain a linear recursive estimate for $\Upsilon^{\beta_{\star}}$; this ultimately gives, similar to the proof of Proposition \ref{prop:hke7}, that
\begin{align}
\Upsilon_{\s,\t}^{\beta_{\star}} \ \lesssim_{\alpha} \ N^{-2+2\alpha}|\mathfrak{l}| + N^{-2+9999\beta_{\star}+\alpha'}. \label{eq:hke8I5}
\end{align}
Plugging \eqref{eq:hke8I5} into \eqref{eq:hke8I3} and choosing $\alpha$ small enough multiple of $\beta_{\star}$ so that $-5\beta_{\star}+\alpha'\leq-3\beta_{\star}$ and $-2+7\beta_{\star}+\alpha'\leq-1-3\beta_{\star}$, we deduce \eqref{eq:hke8I}. We now move to proving \eqref{eq:hke8II}, which follows by basically the same perturbative mechanism. In particular, we first record the following analog of \eqref{eq:hke8I3} but for time-gradients in the forward time-variable, which follows by the proof of \eqref{eq:hke8I3} as all estimates for the spatial gradient of $\mathbf{H}^{\beta_{\star}}$ used to prove \eqref{eq:hke8I3} also hold for forwards time-gradients of $\mathbf{H}^{\beta_{\star}}$; we emphasize the utility of $\t_{N}\geq0$ is that time-dependent factors appearing in these time-gradient estimates for $\mathbf{H}^{\beta_{\star}}$ only improve with $\t_{N}\geq0$, and we therefore do not have to introduce $\t,\s$-dependent factors with additional $\t_{N}$-dependence for convenience:
\begin{align}
|\grad_{\t_{N}}^{\mathbf{T}}\mathbf{H}_{\s,\t,\x,\y}^{\beta_{\star}}| \ \leq \ N^{\beta_{\star}}|\t-\s|^{-1+\alpha}\Gamma_{\frac{\t+\s}{2},\t}^{\beta_{\star}} + N^{-2+9999\beta_{\star}+\alpha'}|\t-\s|^{-1+\alpha} \quad \mathrm{where} \quad \Gamma^{\beta_{\star}}_{\s,\t} \ = \ |\t-\s|^{1-\alpha}\sup_{\x\in\Z}\sup_{\y\in\mathbb{I}^{\beta_{\star}}}|\grad_{\t_{N}}^{\mathbf{T}}\mathbf{H}_{\s,\t,\x,\y}^{\beta_{\star}}|. \nonumber
\end{align}
Given this last bound, it suffices to prove the following $\Gamma^{\beta_{\star}}$ estimate for $\alpha$ a small multiple of $\beta_{\star}$ to obtain \eqref{eq:hke8II}:
\begin{align}
\Gamma_{\s,\t}^{\beta_{\star}} \ \lesssim \ N^{-1+3\beta_{\star}}\t_{N}^{1/2} + N^{-1-4\beta_{\star}}. \label{eq:hke8II1}
\end{align}
To prove \eqref{eq:hke8II1}, we take time-gradients with respect to $\t$ on both sides of the identity in Lemma \ref{lemma:hke5}; this gives
\begin{align}
\grad_{\t_{N}}^{\mathbf{T}}\mathbf{H}_{\s,\t,\x,\y}^{\beta_{\star}} \ = \ \grad_{\t_{N}}^{\mathbf{T}}\mathbf{H}_{\s,\t,\x,\y}-N^{2}\bar{\mathbf{A}}^{N,\beta_{\star}}\int_{\s}^{\t}\grad_{\t_{N}}^{\mathbf{T}}\mathbf{H}_{\r,\t,\x,0}^{\beta_{\star}}\Delta\mathbf{H}_{\s,\r,0,\y}\d\r-N^{2}\bar{\mathbf{A}}^{N,\beta_{\star}}\int_{\t}^{\t+\t_{N}}\mathbf{H}_{\r,\t+\t_{N},\x,0}^{\beta_{\star}}\Delta\mathbf{H}_{\s,\r,0,\y}\d\r. \label{eq:hke8II2}
\end{align}
At this point, we observe that we are almost in the same situation as the proof of Proposition \ref{prop:hke7}. In particular, the time-gradient of $\mathbf{H}^{\beta_{\star}}$ solves the Duhamel-type equation for $\mathbf{H}^{\beta_{\star}}$ in Lemma \ref{lemma:hke5} except the first term on the RHS therein is now replaced by the time-gradient of the $\mathbf{H}$ heat kernel, and we also have a short-time integral given by the last term on the RHS of \eqref{eq:hke8II2}. To estimate this term, we provide the following first step with explanation given afterwards:
\begin{align}
|N^{2}\bar{\mathbf{A}}^{N,\beta_{\star}}\int_{\t}^{\t+\t_{N}}\mathbf{H}_{\r,\t+\t_{N},\x,0}^{\beta_{\star}}\Delta\mathbf{H}_{\s,\r,0,\y}\d\r| \ \lesssim \ \int_{\t}^{\t+\t_{N}}N^{-1+5\beta_{\star}/2}|\t+\t_{N}-\r|^{-1/2} \cdot N^{2\alpha}|\r-\s|^{-1+\alpha}\d\r. \label{eq:hke8II3}
\end{align}
The bound \eqref{eq:hke8II3} follows by the on-diagonal heat kernel bound for $\mathbf{H}^{\beta_{\star}}$ in Proposition \ref{prop:hke8} combined with the bound $|\Delta\mathbf{H}_{\s,\r,0,\y}|\lesssim N^{2\alpha}|\r-\s|^{-1+\alpha}$ for any positive $\alpha$, which is recorded in Proposition A.1 of \cite{DT} and follows by interpolating a maximum principle bound $\mathbf{H}\leq1$ with the classical Gaussian derivative estimate after time-$N^{2}$ speeding-up. We now observe that the $|\r-\s|$ factor is controlled uniformly from below by picking $\r=\t$ in the integral on the RHS of \eqref{eq:hke8II3}; this gives
\begin{align}
\int_{\t}^{\t+\t_{N}}N^{-1+5\beta_{\star}/2}|\t+\t_{N}-\r|^{-1/2} \cdot N^{2\alpha}|\r-\s|^{-1+\alpha}\d\r \ &\lesssim \ N^{-1+5\beta_{\star}/2+2\alpha}|\t-\s|^{-1+\alpha}\int_{\t}^{\t+\t_{N}}|\t+\t_{N}-\r|^{-1/2}. \label{eq:hke8II4}
\end{align}
Observe the integral on the RHS of \eqref{eq:hke8II4} is controlled by $\t_{N}^{1/2}$ by straightforward integration. Combining this with \eqref{eq:hke8II3}, \eqref{eq:hke8II4}, and the remark given after \eqref{eq:hke8II2} about following the proof of Proposition \ref{prop:hke7}, upon replacing $\mathbf{H}^{\beta_{\star}}$ by its time-gradient, the latter of which satisfies the same necessary estimates as $\mathbf{H}^{\beta_{\star}}$ itself, we ultimately deduce, like with the estimate after \eqref{eq:hke8I3} for spatial gradients, the following for $\Gamma^{\beta_{\star}}$ in which $\lambda(N,\t_{N})=N^{-2+9999\beta_{\star}+\alpha'}+N^{-1+5\beta_{\star}/2+2\alpha}\t_{N}^{1/2}$:
\begin{align}
|\grad_{\t_{N}}^{\mathbf{T}}\mathbf{H}_{\s,\t,\x,\y}^{\beta_{\star}}|\ \lesssim \ |\grad_{\t_{N}}^{\mathbf{T}}\mathbf{H}_{\s,\t,\x,\y}| + N^{-1+\beta_{\star}}\int_{\s+\t^{\mathrm{cut}}}^{\t}\Gamma_{\frac{\t+\r}{2},\t}^{\beta_{\star}}|\t-\r|^{-1+\alpha}|\r-\s|^{-\frac32}\d\r + \lambda(N,\t_{N})|\t-\s|^{-1+\alpha}. \label{eq:hke8II5}
\end{align}
Multiplying by $|\t-\s|^{1-\alpha}$, taking a supremum over $\x\in\Z$ and $\y\in\mathbb{I}^{\beta_{\star}}$, and applying the Gronwall inequality, we then deduce the following like how we deduced the $\Upsilon^{\beta_{\star}}$ and $\Gamma^{\beta_{\star}}$ estimates in Proposition \ref{prop:hke7} and even \eqref{eq:hke8I5} from earlier in this proof:
\begin{align}
\Gamma_{\s,\t}^{\beta_{\star}} \ \lesssim \ |\t-\s|^{1-\alpha}\sup_{\x\in\Z}\sup_{\y\in\mathbb{I}^{\beta_{\star}}}|\grad_{\t_{N}}^{\mathbf{T}}\mathbf{H}_{\s,\t,\x,\y}| + \lambda(N,\t_{N}). \label{eq:hke8II6}
\end{align}
For $\alpha$ a small multiple of $\beta$, provided the time-gradient estimates in Proposition A.1 of \cite{DT}, it is a straightforward to deduce \eqref{eq:hke8II1} from \eqref{eq:hke8II6} after some elementary power-counting, so we are done.
\end{proof}
%%%
%%%
\subsection{Off-Diagonal Estimates}
%%%
For our proof of Theorem \ref{theorem:kpz}, similar to \cite{BG,DT} we will need off-diagonal estimates for $\mathbf{H}^{\beta_{\star}}$. We again establish a ``global" off-diagonal estimate along with a ``local" off-diagonal estimate away from the slow bond; this is similar to the difference between Proposition \ref{prop:hke3} and Proposition \ref{prop:hke7}. In particular, we again use a $\beta_{\star}$-perturbative argument below.
%%%
\begin{definition}\label{definition:hke9}
\fsp Provided any real number $\kappa$, we define the following exponential weight
\begin{align}
\exp_{\s,\t,\x,\y}^{N,\kappa} \ = \ \exp\left(\frac{\kappa|\x-\y|}{N^{2}|\t-\s|\vee1}\right).
\end{align}
\end{definition}
%%%
%%%
\begin{remark}\label{remark:hke10}
\fsp We briefly comment that the analysis in this section used to obtain off-diagonal estimates is actually able to lift our estimates with respect to weights in Definition \ref{definition:hke9} to Gaussian weights, though this will not be necessary.
\end{remark}
%%%
Before we state the following result, we comment that the clarifying remarks made before Proposition \ref{prop:hke7} are also in order for Proposition \ref{prop:hke11}; the only difference between these two results/estimates is the off-diagonal behavior.
%%%
\begin{prop}\label{prop:hke11}
\fsp For any positive $\kappa$, we have $\mathbf{H}_{\s,\t,\x,\y}^{\beta_{\star}}\lesssim_{\kappa}N^{-1+3\beta_{\star}}|\t-\s|^{-1/2}\exp_{\s,\t,\x,\y}^{N,-\kappa}$. We also have, for any $|\y|\geq N^{999\beta_{\star}}$ and $\t\geq\s$ satisfying $\s,\t\leq10$, the following estimate in which $\e_{1},\e_{2}\gtrsim\beta_{\star}$ are uniformly positive:
\begin{align}
|\wt{\mathbf{H}}^{\beta_{\star}}_{\s,\t,\x,\y}| \ \lesssim_{\kappa,\e_{1},\e_{2}} \ N^{-1-\e_{1}}|\t-\s|^{-1+\e_{2}}\exp_{\s,\t,\x,\y}^{N,-\kappa} \and |\wt{\mathbf{H}}^{\beta_{\star}}_{\s,\t,\x,\y}| \ \lesssim_{\kappa,\e_{1},\e_{2}} \ \exp_{\s,\t,\x,\y}^{N,-\kappa}.
\end{align}
For any such $|\y|\geq N^{999\beta_{\star}}$ and $\s,\t\leq10$, and for any $\e$ positive, we deduce the estimates
\begin{align}
\mathbf{H}_{\s,\t,\x,\y}^{\beta_{\star}} \ \lesssim_{\kappa,\e} \ N^{-1}|\t-\s|^{-1+\e}\exp_{\s,\t,\x,\y}^{N,-\kappa} \and \mathbf{H}_{\s,\t,\x,\y}^{\beta_{\star}}\ \lesssim_{\kappa,\e}\ \exp_{\s,\t,\x,\y}^{N,-\kappa}.
\end{align}
\end{prop}
%%%
%%%
\begin{proof}
Let us first prove the global estimate without any restrictions on the forwards spatial variable. To this end, we can first pick $\alpha$ sufficiently small but universal so that, by Proposition \ref{prop:hke3}, we have the following interpolation estimate, where the $5\beta_{\star}/2$-exponent in Proposition \ref{prop:hke3} gets boosted to something strictly larger, in this case $3\beta_{\star}$:
\begin{align}
|\mathbf{H}_{\s,\t,\x,\y}^{\beta_{\star}}| \ \leq \ |\mathbf{H}_{\s,\t,\x,\y}^{\beta_{\star}}|^{1-\alpha}|\mathbf{H}_{\s,\t,\x,\y}^{\beta_{\star}}|^{\alpha} \ \lesssim \ N^{-1+3\beta_{\star}}|\t-\s|^{-\frac12}|\mathbf{H}_{\s,\t,\x,\y}^{\beta_{\star}}|^{\alpha}. \label{eq:hke111}
\end{align}
It now suffices to prove the $\mathbf{H}^{\beta_{\star}}$ heat kernel is controlled by $\exp^{N,-\kappa}$ times $\kappa$-dependent constants for any positive $\kappa$; changing $\kappa$ by a factor of $\alpha$, this would prove the global $\mathbf{H}^{\beta_{\star}}$ estimate when combined with \eqref{eq:hke111}. In particular, we show, for any positive $\kappa$,
\begin{align}
|\mathbf{H}_{\s,\t,\x,\y}^{\beta_{\star}}| \ \lesssim_{\kappa} \ \exp^{N,-\kappa}_{\s,\t,\x,\y}. \label{eq:hke112}
\end{align}
We now observe that the LHS is the same as the probability a random walk, whose generator is $\mathbf{A}^{N,\beta_{\star}}\Delta$, goes from $\x$ to $\y$ in time $N^{2}|\t-\s|$. Because the aforementioned generator has no drift, this random walk is a continuous-time martingale with increments of length 1 and speed at most $N^{2}$; it is slower at the origin because of the slow bond. Therefore, the estimate \eqref{eq:hke112} follows from standard martingale inequalities, like the bounds that give the off-diagonal estimates/\eqref{eq:hke112} for $\beta_{\star}=0$ in Proposition A.1 of \cite{DT}. We note that we could even couple the $\mathbf{H}^{\beta_{\star}}$ random walk to the $\mathbf{H}$ random walk by coupling Poisson clocks, so that the maximal displacement of the $\mathbf{H}^{\beta_{\star}}$ random walk is controlled by that of the $\mathbf{H}$ random walk because the latter simply takes more jumps and covers the trajectory of the former, at which point we can employ maximal inequalities for the standard continuous-time symmetric simple random walk on $\Z$, which are classical and can be accessed with elementary/combinatorial methods.

We now restrict to $|\y|\geq N^{999\beta_{\star}}$. Using a similar interpolation, we can again choose $\alpha$ sufficiently small but universal so that, by Proposition \ref{prop:hke7}, we have the following in which $\e_{1},\e_{2}$ are the uniformly positive constants in the statement of Proposition \ref{prop:hke7}:
\begin{align}
|\wt{\mathbf{H}}_{\s,\t,\x,\y}^{\beta_{\star}}| \ = \ |\wt{\mathbf{H}}_{\s,\t,\x,\y}^{\beta_{\star}}|^{1-\alpha}|\wt{\mathbf{H}}_{\s,\t,\x,\y}^{\beta_{\star}}|^{\alpha} \ \lesssim \ N^{-1-\e_{1}/2}|\t-\s|^{-1+\e_{2}}|\wt{\mathbf{H}}_{\s,\t,\x,\y}^{\beta_{\star}}|^{\alpha}. \label{eq:hke113}
\end{align}
Upon changing $\e_{1}$ in the statement of Proposition \ref{prop:hke11} by a factor of $1/2$, it suffices to prove the last $\alpha$-power on the far RHS of is controlled by $\exp^{N,-\kappa}$ times $\kappa$-dependent factors for any positive $\kappa$; this is analogous to how we deduced the global $\mathbf{H}^{\beta_{\star}}$ estimate from \eqref{eq:hke111} earlier in this argument. In particular, we are left to establish \eqref{eq:hke112} but with $\wt{\mathbf{H}}^{\beta_{\star}}$ in place of $\mathbf{H}^{\beta_{\star}}$ on the LHS of \eqref{eq:hke112}. However, this follows by the triangle inequality, which controls $|\wt{\mathbf{H}}^{\beta_{\star}}|\leq\mathbf{H}^{\beta_{\star}}+\mathbf{H}$, and then realizing \eqref{eq:hke112} holds with $\beta_{\star}=0$ as well. Indeed, the simple random walk argument for \eqref{eq:hke112} holds uniformly in all non-negative $\beta_{\star}$. This completes the proof.
\end{proof}
%%%
Again similar to \cite{BG,DT}, for technical and analytic reasons, we will need a version of Proposition \ref{prop:hke11} but for gradients of $\mathbf{H}^{\beta_{\star}}$. We note that the estimates in the following result are not supposed to be optimal in any regard, but they will certainly be sufficient. In particular, we only give the gradient estimates for $\mathbf{H}^{\beta_{\star}}$ itself instead of for $\wt{\mathbf{H}}^{\beta_{\star}}$ as well, as was the case in Proposition \ref{prop:hke11}.
%%%
\begin{prop}\label{prop:hke12}
\fsp There exist $\beta_{\star}\lesssim\e_{1},\e_{2},\e_{3},\e_{4}\lesssim\beta_{\star}$ such that for any $|\mathfrak{l}|\leq N|\t-\s|^{1/2}\vee1$, for any positive $\mathfrak{t}_{N}$, and for any $\kappa$, we have the following gradient estimates for all $\t\geq\s$ satisfying $\s,\t\leq10$:
\begin{align}
|\grad_{\mathfrak{l}}^{\mathbf{X}}\mathbf{H}_{\s,\t,\x,\y}^{\beta_{\star}}| \ &\lesssim_{\kappa,\e_{1},\e_{2}} \ N^{-2+3\beta_{\star}}|\t-\s|^{-1+\e_{2}}|\mathfrak{l}|\exp_{\s,\t,\x,\y}^{N,-\kappa} + N^{-2+99999\beta_{\star}+\e_{2}}|\t-\s|^{-1+\e_{2}}\exp_{\s,\t,\x,\y}^{N,-\kappa} \\
|\grad_{\mathfrak{t}_{N}}^{\mathbf{T}}\mathbf{H}_{\s,\t,\x,\y}^{\beta_{\star}}| \ &\lesssim_{\kappa,\e_{3},\e_{4}} \ N^{-1+6\beta_{\star}+\e_{3}}|\t-\s|^{-1+\e_{3}}\t_{N}^{1/2-\e_{3}}\exp_{\s,\t,\x,\y}^{N,-\kappa}+N^{-2+99999\beta_{\star}+\e_{4}}|\t-\s|^{-1+\e_{4}}\exp_{\s,\t,\x,\y}^{N,-\kappa}.
\end{align}
\end{prop}
%%%
%%%
\begin{proof}
Following the proof for Proposition \ref{prop:hke11} via interpolation, because the proposed estimates herein are worse than those in Proposition \ref{prop:hke8} by small powers of $N$ and a slightly smaller power of $\t_{N}$, it suffices to prove that the LHS of both estimates/both gradients of $\mathbf{H}^{\beta_{\star}}$ are controlled by $\exp^{N,-\kappa}$. This, however, follows from controlling said gradients by $\mathbf{H}^{\beta_{\star}}$ itself at two different space-time coordinates and applying \eqref{eq:hke112} for each of these $\mathbf{H}^{\beta_{\star}}$ copies. Let us clarify that because we have $\mathbf{H}^{\beta_{\star}}$ at two different space-time coordinates, we actually get upper bounds in terms of $\exp^{N,-\kappa}$ at two different space-time coordinates. However, the assumption $|\mathfrak{l}|\leq N|\t-\s|^{1/2}$ and the positivity of $\t_{N}$ let us bound $\exp^{N,-\kappa}$ at said different space-time coordinates by $\exp^{N,-\kappa}$ at the same space-time coordinate if we allow for an extra uniformly bounded factor, so we are done.
\end{proof}
%%%
%
%
%
%%%
\section{Local Probabilistic Estimates}\label{section:lpe}
%%%
This section is dedicated to analyzing the last term in \eqref{eq:mshe3}. In particular, we provide the key estimate for this term that allows us to basically forget it when analyzing \eqref{eq:mshe3}. Before we state this estimate, however, let us introduce notation for the exponential-weighted space-time uniform norm we estimate this term with respect to. We emphasize that this norm is \emph{not} the one that is used in \cite{BG,DT} to analyze the Gartner transforms therein. In particular, we will only employ the following norm when making comparisons or estimating differences between the solution to \eqref{eq:mshe3} and the solution to \eqref{eq:mshe3} with some adjustments.
%%%
\begin{definition}\label{definition:lpe1}
\fsp Provided any function $\phi:\R_{\geq0}\times\Z\to\R$ and any subsets $[0,\t_{+}]\subseteq\R_{\geq0}$ and $\mathbb{K}\subseteq\R$, which could be a discrete subset of $\Z$, for example, we define the norm $\|\phi\|_{\t_{+};\mathbb{K};\max}=\sup_{0\leq\t\leq\t_{+}}\sup_{\x\in\mathbb{K}}|\phi_{\t,\x}|$. We also define the exponentially-weighted norm $\|\phi\|_{\t_{+};\mathbb{K}}=\sup_{0\leq\t\leq\t_{+}}\sup_{\x\in\mathbb{K}}\exp(-|\x|/N)|\phi_{\t,\x}|$, which we use in the next section.
\end{definition}
%%%
%%%
\begin{remark}\label{remark:lpe2}
\fsp We emphasize that the specific choice of $-1$ for the exponent in $\|\|_{\t_{+};\mathbb{K}}$-norms defined above is because we assumed the Gartner transform is bounded in moments uniformly in $\Z$, whereas if we had instead exponential-linear growth of the Gartner transform as in \cite{BG,DT}, we would need to increase this negative weight in the exponential. However, as long as this weight is fixed and independent of the scaling parameter $N$, this is not a big deal; as noted after Theorem \ref{theorem:kpz}, all we would need to do is standard insertion of exponential weights as in \cite{BG,DT}. We only make this comment in order to clarify notation/details of our work.
\end{remark}
%%%
%%%
\begin{prop}\label{prop:lpe3}
\fsp Recalling notation in \emph{Proposition \ref{prop:mshe1}} and \emph{Corollary \ref{corollary:mshe3}}, we have the following estimate for any possibly random time $\t_{\mathrm{st}}\in[0,1]$; we provide some clarification for the estimate below afterwards:
\begin{align}
\mathbf{P}\left(\|\mathbf{H}^{\beta_{\star},\mathbf{T}}_{\t,\x}\left(2^{-1}c_{N}\mathbf{1}_{\y=0}\mathfrak{q}_{\s,\y}^{\mathrm{tot}}\mathbf{Z}_{\s,\y}^{N}\right)\|_{\t_{\mathrm{st}};\Z;\max}(1+\|\mathbf{Z}^{N}\|_{\t_{\mathrm{st}};\Z}^{2})^{-1} \geq N^{-\gamma}\right) \ \lesssim \ N^{-\gamma}. \label{eq:lpe3I}
\end{align}
The constant $\gamma>0$ is positive and depends only on $\beta_{\star}$. The first norm in the expectation on the LHS of \eqref{eq:lpe3I} is taken with respect to $\t,\x$-variables that the space-time heat operator is evaluated at; the $\s,\y$-variables are for the integration-variables inside the heat operator definition, and in particular not for the aforementioned norm that is inside.
\end{prop}
%%%
The proof of Proposition \ref{prop:lpe3} will be quite lengthy and require several procedures, so we actually defer it to the final section of this paper; the interested reader is invited to skip directly there if desired. We conclude this section by briefly mentioning its utility, which is explained in more detail shortly in the upcoming section. Assuming $\mathbf{Z}^{N}$ is reasonably behaved and bounded from above, which is plausible given that it is supposed to look like the regular $\mathrm{SHE}$ according to Theorem \ref{theorem:kpz}, we can forget the $\|\mathbf{Z}^{N}\|^{2}$-factor inside the expectation in \eqref{eq:lpe3I}. Proposition \ref{prop:lpe3} then tells us the final term in \eqref{eq:mshe3} is negligibly small in probability with respect to the space-time uniform norm $\|\|_{1;\Z}$, at which point we can basically perform classical linear PDE analysis to forget the last term in \eqref{eq:mshe3}. The key difficulty is now controlling $\mathbf{Z}^{N}$, which we took for granted in this paragraph; this is addressed in the next section.
%
%
%
%%%
\section{Pathwise Comparisons}\label{section:pc}
%%%
The purpose of this section is basically to estimate the difference in the solution to \eqref{eq:mshe3} if we remove the last term therein and replace heat operators for our choice of $\beta_{\star}$ for this paper/Theorem \ref{theorem:kpz} with the integrable choice $\beta_{\star}=0$, which simply recovers ASEP of \cite{BG}. This estimation will be done exactly in this order, and it will be done with respect to the $\|\|_{1;\Z}$ norm in Definition \ref{definition:lpe1}.
%%%
\begin{definition}\label{definition:pc1}
\fsp Define $\mathbf{Y}^{N}$ and $\mathbf{W}^{N}$ as the solutions to the following two stochastic integral equations, which employ multiplicative noises that we clarify afterwards; we emphasize the following are posed on $\R_{\geq0}\times\Z$ similar to $\mathbf{Z}^{N}$ in Definition \ref{definition:mshe4}:
\begin{align}
\mathbf{Y}^{N}_{\t,\x} \ &= \ \mathbf{H}^{\beta_{\star},\mathbf{X}}_{\t,\x}(\mathbf{Z}^{N}) + \mathbf{H}^{\beta_{\star},\mathbf{T}}_{\t,\x}(\mathbf{Y}^{N}\d\xi^{N}) \and \mathbf{W}^{N}_{\t,\x} \ = \ \mathbf{H}^{\mathbf{X}}_{\t,\x}(\mathbf{Z}^{N}) + \mathbf{H}^{\mathbf{T}}_{\t,\x}(\mathbf{W}^{N}\d\xi^{N}).
\end{align}
The multiplicative noise term $\mathbf{Y}^{N}\d\xi^{N}$ is defined to be the jump in the $\d\xi^{N}$ martingale differential, minus the deterministic drift, and then times $\mathbf{Y}^{N}$ all at the same space-time point. The same is true for $\mathbf{W}^{N}\d\xi^{N}$ except $\mathbf{W}^{N}$ instead of $\mathbf{Y}^{N}$.
\end{definition}
%%%
%%%
\begin{prop}\label{prop:pc2}
\fsp Define the following two events in which $\gamma_{2}=999^{-999}\gamma$ with $\gamma$ from \emph{Proposition \ref{prop:lpe3}} and in which $\alpha>0$:
\begin{align}
\mathcal{E}(1) \ = \ \left\{\|\mathbf{Z}^{N}-\mathbf{Y}^{N}\|_{1;\Z}\geq N^{-\gamma_{2}}\right\} \and \mathcal{E}(2) \ = \ \left\{\|\mathbf{Y}^{N}-\mathbf{W}^{N}\|_{1;\Z}\geq N^{-\alpha}\right\}.
\end{align}
There exists a uniformly, in $N$, positive choice of $\alpha$ for $\mathcal{E}(2)$ and a positive constant $\gamma_{3}$ depending only on $\beta_{\star}$ such that
\begin{align}
\mathbf{P}(\mathcal{E}(1))+\mathbf{P}(\mathcal{E}(2)) \ \lesssim \ N^{-\gamma_{3}}. \label{eq:pc2I}
\end{align}
\end{prop}
%%%
At a technical level, the proof of Proposition \ref{prop:pc2} has several components that we list below.
%%%
\begin{itemize}
\item We will first show that $\mathbf{Y}^{N}$ and $\mathbf{W}^{N}$ are ``somewhat reasonably well-behaved" in the sense of $\|\|_{1;\Z}$-a priori estimates. This will not be very difficult and mostly based on \cite{BG,DT} analysis of Gartner transforms. In particular, we will control $\|\|_{1;\Z}$-norms with a union bound and moment estimates for $\mathbf{Y}^{N}$ and $\mathbf{W}^{N}$. We note our $\|\|_{1;\Z}$-estimates will not have to be optimal and can grow in $N$ sufficiently slowly; see Lemma \ref{lemma:pc3} below and its proof for details on this bullet point.
\item Next, we prove the $\mathcal{E}(1)$-estimate. If we knew, a priori, that $\mathbf{Z}^{N}$ were reasonably bounded in the same sense as $\mathbf{Y}^{N}$ and $\mathbf{W}^{N}$ in the previous bullet point, then Proposition \ref{prop:lpe3} and standard linear theory as in \cite{BG,DT} would be enough to estimate $\mathbf{Z}^{N}-\mathbf{Y}^{N}$. However, such an a priori estimate for $\mathbf{Z}^{N}$ is not trivial to obtain; in principle, direct estimates for the third term on the RHS of \eqref{eq:mshe3} may blow up if we are not careful. This is a major reason why Theorem \ref{theorem:kpz} is difficult. To circumvent this issue, we employ a stochastic continuity argument. Roughly speaking, if we expect $\mathbf{Y}^{N}\approx\mathbf{Z}^{N}$ and $\mathbf{Y}^{N}$ is reasonable, then $\mathbf{Z}^{N}$ is also reasonable. In particular, it is more reasonable than we may need for Proposition \ref{prop:lpe3} to inform us that $\mathbf{Y}^{N}\approx\mathbf{Z}^{N}$ in the first place. Because $\mathbf{Z}^{N}$ cannot change too much in very short sub-microscopic times, we ultimately deduce that if $\mathbf{Y}^{N}\approx\mathbf{Z}^{N}$ until some time $\t_{+}$, then we know $\mathbf{Y}^{N}\approx\mathbf{Z}^{N}$ plus a small error for a small time after $\t_{+}$. But this last estimate after time $\t_{+}$ is enough to imply the last term in \eqref{eq:mshe3} is negligible by Proposition \ref{prop:lpe3}, thus $\mathbf{Y}^{N}\approx\mathbf{Z}^{N}$ not just until but after $\t_{+}$ without the aforementioned small error. This continuity argument, which is applied often in PDE analysis, depends heavily on the fact that our a priori estimates for $\mathbf{Y}^{N}$ and $\mathbf{W}^{N}$ are uniform in space-time/with respect to $\|\|_{1;\Z}$, and most importantly the same is true for the estimate in Proposition \ref{prop:lpe3}. Moreover, the precise estimates that we obtain to prove the $\mathcal{E}(1)$-estimate in Proposition \ref{prop:pc2} will be slightly worse than what Proposition \ref{prop:lpe3} suggests, for example, which is technically important to making this stochastic continuity work.
\item Finally, we now explain the $\mathcal{E}(2)$-estimate. This basically follows the previous bullet point, except that we already have good a priori estimates for both $\mathbf{Y}^{N}$ and $\mathbf{W}^{N}$, and instead of Proposition \ref{prop:lpe3}, we must now compare the heat kernels in terms of $\beta_{\star}$ via Proposition \ref{prop:hke8} and Proposition \ref{prop:hke11}. These latter estimates only provide a little-oh comparison of heat operators and heat kernels in some neighborhood away from the slow bond, but this does not introduce any real difficulties because said neighborhood of the slow bond is mesoscopic and of length much smaller than the macroscopic length-scale $N$.
\end{itemize}
%%%
We conclude this section by inviting the reader to take Proposition \ref{prop:pc2} for granted, at least upon a first reading, and skip to the next section. That is, the proof of Proposition \ref{prop:pc2} contains many technical details that will not be used in the rest of this paper, though it is not sufficiently long or opaque that we defer it entirely to another section.
%%%
\subsection{A Priori Estimates}
%%%
We will first address the first bullet point above. In what follows, the constant $\delta$ can be chosen to be any arbitrarily small positive constant, as long as we keep it fixed; the reader may take $\delta=999^{-999}\gamma_{2}$, with $\gamma_{2}$ from Proposition \ref{prop:pc2}.
%%%
\begin{lemma}\label{lemma:pc3}
\fsp Let $\mathcal{E}(\mathbf{Y}^{N};\delta)$ be the event where $\|\mathbf{Y}^{N}\|_{1;\Z}\geq N^{\delta}$, and let $\mathcal{E}(\mathbf{W}^{N};\delta)$ be the event where $\|\mathbf{W}^{N}\|_{1;\Z}\geq N^{\delta}$, where $\delta$ is positive and fixed and otherwise arbitrarily small. Provided any positive constant $D$, we have
\begin{align}
\mathbf{P}(\mathcal{E}(\mathbf{Y}^{N};\delta))+\mathbf{P}(\mathcal{E}(\mathbf{W}^{N};\delta)) \ \lesssim_{D} \ N^{-D}.
\end{align}
\end{lemma}
%%%
%%%
\begin{proof}
We will prove the $\mathbf{Y}^{N}$ and $\mathcal{E}(\mathbf{Y}^{N};\delta)$ estimates; the proof of $\mathbf{W}^{N}$ and $\mathcal{E}(\mathbf{W}^{N};\delta)$ estimates follows from the same argument, because all ingredients, such as heat kernel estimates for $\mathbf{H}^{\beta_{\star}}$, also hold when analyzing $\mathbf{W}^{N}$ and the $\beta_{\star}=0$ heat kernel $\mathbf{H}$. To this end, we will first employ Lemma \ref{lemma:ste}. Letting $\mathbb{I}^{\mathbf{T}}=\{\mathfrak{j}N^{-100}\}_{\mathfrak{j}=0,\ldots,N^{100}}$ denote a sub-microscopic-scale discretization of the continuum time-interval $[0,1]$, which is finite and therefore approachable via union bound, we have the following by Lemma \ref{lemma:ste}:
\begin{align}
\|\mathbf{Y}^{N}\|_{1;\Z} \ \lesssim \ \sup_{(\t,\x)\in\mathbb{I}^{\mathbf{T}}\times\Z}\exp(-|\x|/N)|\mathbf{Y}^{N}_{\t,\x}|. \label{eq:pc31}
\end{align}
We now decompose $\Z=\cup_{\mathfrak{k}\in\Z}\mathbb{I}^{\mathbf{X},\mathfrak{k}}$ into almost-disjoint blocks of length {$N$} given by {$\mathbb{I}^{\mathbf{X},\mathfrak{k}}$}$=\llbracket\mathfrak{k}N,\mathfrak{k}N+N\rrbracket$, so that
\begin{align}
\sup_{(\t,\x)\in\mathbb{I}^{\mathbf{T}}\times\Z}\exp(-|\x|/N)|\mathbf{Y}^{N}_{\t,\x}| \ \lesssim \ \sup_{\mathfrak{k}\in\Z}\sup_{(\t,\x)\in\mathbb{I}^{\mathbf{T}}\times\mathbb{I}^{\mathbf{X},\mathfrak{k}}}\exp(-|\mathfrak{k}|)|\mathbf{Y}^{N}_{\t,\x}|. \label{eq:pc32}
\end{align}
Courtesy of \eqref{eq:pc31} and \eqref{eq:pc32}, we deduce the probability of $\mathcal{E}(\mathbf{Y}^{N};\delta)$ is controlled by the probability the RHS of \eqref{eq:pc32} exceeds $N^{\delta/2}$, with the decrease exponent accounting for implied constants in \eqref{eq:pc31} and \eqref{eq:pc32}. By union bound over $\mathfrak{k}\in\Z$ on the RHS of \eqref{eq:pc32}, we additionally deduce the following probability estimate, in which $\mathcal{E}(\mathbf{Y}^{N};\delta;\mathfrak{k})$ denotes the event in which $|\mathbf{Y}^{N}|$ exceeds $N^{\delta/2}$ when restricted to the space-time discrete set $\mathbb{I}^{\mathbf{T}}\times\mathbb{I}^{\mathbf{X},\mathfrak{k}}$ after multiplying by the sub-exponential factor $\exp(-|\mathfrak{k}|)$ on the RHS of \eqref{eq:pc32}:
\begin{align}
\mathbf{P}(\mathcal{E}(\mathbf{Y}^{N};\delta)) \ \leq \ {\sum}_{\mathfrak{k}\in\Z}\mathbf{P}(\mathcal{E}(\mathbf{Y}^{N};\delta;\mathfrak{k})). \label{eq:pc33}
\end{align}
We will now estimate $\mathbf{P}(\mathcal{E}(\mathbf{Y}^{N};\delta;\mathfrak{k}))$ probabilities in summable fashion in $\mathfrak{k}\in\Z$; for this, the sub-exponential factor $\exp(-|\mathfrak{k}|)$ will be crucial. To start, we employ \emph{another} union bound per $\mathfrak{k}\in\Z$ over the space-time discrete set $\mathbb{I}^{\mathbf{T}}\times\mathbb{I}^{\mathbf{X},\mathfrak{k}}$; this gives the bound below, in which $\mathcal{E}(\mathbf{Y}^{N};\delta;\mathfrak{k};\t,\x)$ is the probability that $\exp(-|\mathfrak{k}|)|\mathbf{Y}^{N}_{\t,\x}|\geq N^{\delta/2}$:
\begin{align}
\mathbf{P}(\mathcal{E}(\mathbf{Y}^{N};\delta;\mathfrak{k})) \ \leq \ {\sum}_{(\t,\x)\in\mathbb{I}^{\mathbf{T}}\times\mathbb{I}^{\mathbf{X},\mathfrak{k}}}\mathbf{P}(\mathcal{E}(\mathbf{Y}^{N};\delta;\mathfrak{k};\t,\x)) \ \lesssim \ N^{101}\sup_{(\t,\x)\in\mathbb{I}^{\mathbf{T}}\times\mathbb{I}^{\mathbf{X},\mathfrak{k}}}\mathbf{P}(\mathcal{E}(\mathbf{Y}^{N};\delta;\mathfrak{k};\t,\x)), \label{eq:pc34}
\end{align}
where the final estimate in \eqref{eq:pc34} follows from observing $|\mathbf{I}^{\mathbf{T}}|\lesssim N^{100}$ and $\mathbf{I}^{\mathbf{X},\mathfrak{k}}$ is a block of length order $N$ uniformly in $\mathfrak{k}\in\Z$. Now, we will estimate each probability on the RHS of \eqref{eq:pc34} uniformly in $(\t,\x)$, starting with the Chebyshev inequality for $p\geq1$:
\begin{align}
\mathbf{P}(\mathcal{E}(\mathbf{Y}^{N};\delta;\mathfrak{k};\t,\x)) \ \leq \ N^{-p\delta}\exp(-2p|\mathfrak{k}|)\E|\mathbf{Y}^{N}_{\t,\x}|^{2p}. \label{eq:pc35}
\end{align}
Combining \eqref{eq:pc33}, \eqref{eq:pc34}, and \eqref{eq:pc35} and defining $\mathbf{Y}_{}^{N,p}$ as the supremum of $\E|\mathbf{Y}^{N}_{\t,\x}|^{2p}$ over $(\t,\x)\in[0,1]\times\Z$, for $p\geq1$ we get
\begin{align}
\mathbf{P}(\mathcal{E}(\mathbf{Y}^{N};\delta)) \ \lesssim \ N^{101-p\delta}\mathbf{Y}_{}^{N,p}{\sum}_{\mathfrak{k}\in\Z}\exp(-2p|\mathfrak{k}|) \ \lesssim \ N^{101-p\delta}\mathbf{Y}^{N,p}_{}. \label{eq:pc36}
\end{align}
According to \eqref{eq:pc36}, it suffices to bound $\E|\mathbf{Y}^{N}_{\t,\x}|^{2p}$ uniformly in $(\t,\x)$, which would estimate $\mathbf{Y}^{N,p}_{}$, by a $p$-dependent constant, in which $p\geq1$ is chosen sufficiently large so that $101-p\delta\leq -D$ for the exponent $D$ chosen in the statement of Lemma \ref{lemma:pc3}. We note that the point behind our estimates in the proof thus far are to reduce estimating $\|\mathbf{Y}^{N}\|_{1;\Z}$, which is uniform and simultaneous control of $\mathbf{Y}^{N}$ in space-time, in terms of one-point space-time moment bounds for $\mathbf{Y}^{N,p}$. To estimate $\mathbf{Y}^{N,p}_{}$, we will appeal to its defining stochastic integral equation; according to Definition \ref{definition:pc1}, we first have, in which $\|\cdot\|_{\omega;2p}=(\E|\cdot|^{2p})^{1/2p}$, the bound
\begin{align}
\|\mathbf{Y}^{N}_{\t,\x}\|_{\omega;2p}^{2} \ &\lesssim \ \|\mathbf{H}^{\beta_{\star},\mathbf{X}}_{\t,\x}(\mathbf{Z}^{N})\|_{\omega;2p}^{2} + \|\mathbf{H}^{\beta_{\star},\mathbf{T}}_{\t,\x}(\mathbf{Y}^{N}\d\xi^{N})\|_{\omega;2p}^{2}. \label{eq:pc37}
\end{align}
To estimate the first term on the RHS of \eqref{eq:pc37}, we follow the beginning of the proof of Proposition 3.2 in \cite{DT}. In particular, the spatial heat operator $\mathbf{H}^{\beta_{\star},\mathbf{X}}$ is still convex, positivity-preserving, and maps the constant function to itself, as this heat kernel still describes a probability measure in its forward spatial variable. Thus, because we assume near-stationary initial data, the first term on the RHS of \eqref{eq:pc37} is controlled by a $p$-dependent constant. It remains to control the second term on the RHS of \eqref{eq:pc37}. For this, we again follow the proof of Proposition 3.2 in \cite{DT} and its analysis of the martingale integral relevant therein. However, instead of the BDG-type inequality of Lemma 3.1 in \cite{DT}, we employ Lemma \ref{lemma:mge}; the only resulting difference is that we instead end up integrating, in time, the spatial supremum $\mathbf{Y}^{N,p,\mathbf{X}}_{\s}$ of $\E|\mathbf{Y}_{\s,\x}^{N}|^{2p}$ over $\x\in\Z$, not the pointwise moments of $\mathbf{Y}^{N}$ like in the proof of Proposition 3.2 in \cite{DT}, though this does not change the method in the proof of Porposition 3.2 in \cite{DT}. Thus, we have the following, in which the correct power of $\mathbf{Y}^{N,p,\mathbf{X}}$ on the RHS can be seen by noting the LHS of the estimate below scales quadratically in $\mathbf{Y}^{N}$, whereas $(\mathbf{Y}^{N,p,\mathbf{X}})^{1/2p}$ scales linearly in $\mathbf{Y}^{N}$:
\begin{align}
\|\mathbf{H}^{\beta_{\star},\mathbf{T}}_{\t,\x}(\mathbf{Y}^{N}\d\xi^{N})\|_{\omega;2p}^{2} \ \lesssim_{p} \ N\int_{0}^{\t}(\mathbf{Y}^{N,p,\mathbf{X}}_{\s})^{1/p}{\sum}_{\y\in\Z}|\mathbf{H}_{\s,\t,\x,\y}^{\beta_{\star}}|^{2}\d\s + N^{-1}(\mathbf{Y}^{N,p,\mathbf{X}}_{0})^{1/p}. \label{eq:pc38}
\end{align}
We note the second term on the RHS of \eqref{eq:pc38} is controlled by a $p$-dependent constant because it is the initial data and therefore controlled a priori for near-stationary data. We clarify that its relevance comes from considering time-scales $\t\leq N^{-2}$ on the LHS of \eqref{eq:pc38} where the evolution of the LHS of \eqref{eq:pc38} can basically be tracked in terms of jumps in $\mathbf{Y}^{N}\d\xi^{N}$; this is done implicitly in the proof of Proposition 3.2 of \cite{DT} when considering the time-factor $1\wedge|\t-\s|^{-1/2}$ therein. If we were to {follow the} analysis in the proof of Proposition 3.2 in \cite{DT}, we would employ on-diagonal estimates for $\mathbf{H}^{\beta_{\star}}$. However, here, the heat kernel estimates for $\mathbf{H}^{\beta_{\star}}$ uniformly over $\Z$ are not enough to beat the factor of $N$ on the RHS of \eqref{eq:pc38}. To remedy this technical issue, we decompose the sum over $\y\in\Z$ in the integral on the RHS of \eqref{eq:pc38} into two pieces. The first piece is over the small set of $|\y|\leq N^{999\beta_{\star}}$, and the second piece is over the complement of this set. To be precise, we have the inequality
\begin{align}
N\int_{0}^{\t}(\mathbf{Y}^{N,p,\mathbf{X}}_{\s})^{1/p}{\sum}_{\y\in\Z}|\mathbf{H}_{\s,\t,\x,\y}^{\beta_{\star}}|^{2}\d\s \ &= \ \Phi_{1,\t,\x} + \Phi_{2,\t,\x}, \label{eq:pc39}
\end{align}
where the pieces $\Phi_{1,\t,\x}$ and $\Phi_{2,\t,\x}$, defined below explicitly, were described prior to \eqref{eq:pc39}:
\begin{align*}
\Phi_{1,\t,\x} \ \overset{\bullet}= \ N\int_{0}^{\t}(\mathbf{Y}^{N,p,\mathbf{X}}_{\s})^{1/p}{\sum}_{|\y|\leq N^{999\beta_{\star}}}|\mathbf{H}_{\s,\t,\x,\y}^{\beta_{\star}}|^{2}\d\s \and \Phi_{2,\t,\x} \ \overset{\bullet}= \ N\int_{0}^{\t}(\mathbf{Y}^{N,p,\mathbf{X}}_{\s})^{1/p}{\sum}_{|\y|> N^{999\beta_{\star}}}|\mathbf{H}_{\s,\t,\x,\y}^{\beta_{\star}}|^{2}\d\s. 
\end{align*}
We will first control $\Phi_{1,\t,\x}$. To this end, we control the sum over $\y$ therein by $N^{999\beta_{\star}}$, which is the size of the set of summation, times the supremum of $|\mathbf{H}^{\beta_{\star}}|^{2}$. This latter supremum is controlled by interpolating the on-diagonal estimate of Proposition \ref{prop:hke3} and the maximum principle inequality $|\mathbf{H}^{\beta_{\star}}|\leq1$. Thus, for arbitrarily small but fixed and positive $\varrho$, we have a positive $\varrho'$ so that
\begin{align}
\Phi_{1,\t,\x} \ \lesssim \ N^{1+999\beta_{\star}-2+5\beta_{\star}+\varrho}\int_{0}^{\t}(\mathbf{Y}_{\s}^{N,p,\mathbf{X}})^{1/p}|\t-\s|^{-1+\varrho'}\d\s \ \lesssim \ N^{-1+9999\beta_{\star}+\varrho}\int_{0}^{\t}(\mathbf{Y}_{\s}^{N,p,\mathbf{X}})^{1/p}|\t-\s|^{-1+\varrho'}\d\s. \label{eq:pc310}
\end{align}
To estimate $\Phi_{2,\t,\x}$ on the RHS of \eqref{eq:pc39}, we note $\mathbf{H}^{\beta_{\star}}$ satisfies the same on-diagonal upper bound as the heat kernel in Proposition 3.2 in \cite{DT}; this follows from Proposition \ref{prop:hke11} with $\kappa=0$. Therefore, similar to the proof of Proposition 3.2 in \cite{DT}, we have the estimate below as we also have that $\mathbf{H}^{\beta_{\star}}$ is a probability measure on $\Z$ in its forward spatial variable:
\begin{align}
\Phi_{2,\t,\x} \ \lesssim \ \int_{0}^{\t}(\mathbf{Y}_{\s}^{N,p,\mathbf{X}})^{1/p}|\t-\s|^{-1/2}\d\s. \label{eq:pc311}
\end{align}
We now combine \eqref{eq:pc37}, \eqref{eq:pc38}, \eqref{eq:pc39}, \eqref{eq:pc310}, and \eqref{eq:pc311} along with our estimates, in terms of $p$-dependent constants, for the first term on the RHS of \eqref{eq:pc37} and second term on the RHS of \eqref{eq:pc38}. This provides
\begin{align}
\|\mathbf{Y}^{N}_{\t,\x}\|_{\omega;2p}^{2} \ \lesssim_{p} \ 1 + N^{-1+9999\beta_{\star}+\varrho}\int_{0}^{\t}(\mathbf{Y}_{\s}^{N,p,\mathbf{X}})^{1/p}|\t-\s|^{-1+\varrho'}\d\s + \int_{0}^{\t}(\mathbf{Y}_{\s}^{N,p,\mathbf{X}})^{1/p}|\t-\s|^{-1/2}\d\s. \label{eq:pc312}
\end{align}
Because the RHS of \eqref{eq:pc312} is uniform in the variable $\x\in\Z$ present on the LHS of \eqref{eq:pc312}, the previous estimate \eqref{eq:pc312} extends the same upper bound for $(\mathbf{Y}_{\t}^{N,p,\mathbf{X}})^{1/p}$, since this final term is given by the supremum over $\x\in\Z$ of the LHS of \eqref{eq:pc312}. To deduce a $p$-dependent upper bound for $\mathbf{Y}^{N,p,\mathbf{X}}_{\t}$ that is uniform in $\t\in[0,1]$, it suffices to then employ the Gronwall inequality, upon taking $\varrho$ in \eqref{eq:pc312} to be small and fixed. As such upper bound on $\mathbf{Y}^{N,p,\mathbf{X}}$ is uniform in $\t\in[0,1]$, it extends to an upper bound for $\mathbf{Y}^{N,p}$ itself as well; the term $\mathbf{Y}^{N,p}$ was defined prior to \eqref{eq:pc36}. This implies the desired $p$-dependent upper bound for $\mathbf{Y}^{N,p}$, which completes the proof courtesy of the sentence after \eqref{eq:pc36}, so we are done.
\end{proof}
%%%
%%%
\subsection{$\mathcal{E}(1)$ Estimate}
%%%
We start by constructing a set of random times that the aforementioned stochastic continuity argument will be based on. In particular, we propagate these random times to the final time 1 with high probability, which then ultimately yields the $\mathcal{E}(1)$ estimate claimed in Proposition \ref{prop:pc2}.
%%%
\begin{definition}\label{definition:pc4}
\fsp We will first define $\mathfrak{t}_{1}$ as the first time where $\mathbf{Y}^{N}$ exceeds $N^{\delta}$ like in Lemma \ref{lemma:pc3}. We also define $\mathfrak{t}_{2}$ as the first time where the difference $\mathbf{Z}^{N}-\mathbf{Y}^{N}$ exceeds $N^{-\gamma_{2}}$ with $\gamma_{2}$ from Proposition \ref{prop:pc2}. Precisely, recalling norms in Definition \ref{definition:lpe1}, we define
\begin{align}
{\t_{1}} \ = \ \inf\left\{\t\in[0,1]: \ \|\mathbf{Y}^{N}\|_{\t;\Z} \geq N^{\delta}\right\}\wedge1 \and \t_{2} \ = \ \inf\left\{\t\in[0,1]:\|\mathbf{Z}^{N}-\mathbf{Y}^{N}\|_{\t;\Z}\geq N^{-\gamma_{2}}\right\}\wedge1.
\end{align}
Define $\t_{3}$ as the first time where the estimate in Proposition \ref{prop:lpe3} fails; in particular, the constant $\gamma$ below is that in Proposition \ref{prop:lpe3}:
\begin{align}
\t_{3} \ = \ \inf\left\{\t\in[0,1]:\|\mathbf{H}^{\beta_{\star},\mathbf{T}}_{\t,\x}(2^{-1}c_{N}\mathbf{1}_{\y=0}\mathfrak{q}_{\s,\y}^{\mathrm{tot}}\mathbf{Z}^{N}_{\s,\y})\|_{\t;\Z;\max}\geq N^{-\gamma/2}\right\}.
\end{align}
Lastly, we glue these random times together and define $\t_{\mathbf{Y}}=\t_{1}\wedge\t_{2}\wedge\t_{3}$.
\end{definition}
%%%
%%%
\begin{remark}\label{remark:pc5}
{\fsp We briefly note that $\t_{3}$ is never infinite. Indeed, the space-time heat operator $\mathbf{H}^{\beta_{\star},\mathbf{T}}$ is always equal to zero at $\t=0$, so the infimum defining $\t_{3}$ is never of the empty set.}
\end{remark}
%%%
%%%
\begin{remark}\label{remark:pc6}
\fsp Observe that before the time $\t_{1}\wedge\t_{2}$, and thus before the time $\t_{\mathbf{Y}}$, we deduce $\|\mathbf{Z}^{N}\|\leq N^{3\delta}$ if $N$ is sufficiently large depending only on $\delta$ and $\gamma_{2}$. This follows just by the triangle inequality and constructions of $\t_{1}$ and $\t_{2}$ from Definition \ref{definition:pc4}.
\end{remark}
%%%
Proposition \ref{prop:pc2} follows from showing $\t_{2}=1$ with the required high probability. In particular, it will be enough to show $\t_{\mathbf{Y}}=1$ with the required high probability. To make this presentation clean, we introduce a few more pieces of notation.
%%%
\begin{definition}\label{definition:pc7}
\fsp Let us first define $\mathbf{D}^{N,1}=\mathbf{Z}^{N}-\mathbf{Y}^{N}$; we first observe that $\mathbf{D}^{N,1}$ solves the following stochastic integral equation, which can be readily verified by \eqref{eq:mshe3} for $\mathbf{Z}^{N}$ and construction of $\mathbf{Y}^{N}$ in Definition \ref{definition:pc1}:
\begin{align}
\mathbf{D}^{N,1}_{\t,\x} \ = \ \mathbf{H}^{\beta_{\star},\mathbf{T}}_{\t,\x}(\mathbf{D}^{N,1}\d\xi^{N}) + \mathbf{H}^{\beta_{\star},\mathbf{T}}_{\t,\x}(2^{-1}c_{N}\mathbf{1}_{\y=0}\mathfrak{q}_{\s,\y}^{\mathrm{tot}}\mathbf{Z}^{N}_{\s,\y}).
\end{align}
We now define a slight modification of $\mathbf{D}^{N,1}$, to be denoted by $\mathbf{D}^{N,1,1}$, through the following stochastic integral equation, whose solutions are certainly unique if we view the first term on the RHS below as a linear term in the solution $\mathbf{D}^{N,1,1}$ and if we view the second term therein as an external ``force" or ``source" term; we note that because $\t_{\mathbf{Y}}\leq\t_{3}$ by construction in Definition \ref{definition:pc4}, and because the second term on the RHS below may be viewed as a spatial heat operator smoothing, evaluated at time $\t-\t_{\mathbf{Y}}\wedge\t$, of a space-time heat operator evaluated at time $\t_{\mathbf{Y}}\wedge\t$, the second term on the RHS below is controlled by $N^{-\gamma/2}$ in absolute value:
\begin{align}
\mathbf{D}^{N,1,1}_{\t,\x} \ = \ \mathbf{H}^{\beta_{\star},\mathbf{T}}_{\t,\x}(\mathbf{D}^{N,1,1}\d\xi^{N}) + \mathbf{H}^{\beta_{\star},\mathbf{T}}_{\t,\x}(\mathbf{1}_{\s\leq\t_{\mathbf{Y}}}2^{-1}c_{N}\mathbf{1}_{\y=0}\mathfrak{q}_{\s,\y}^{\mathrm{tot}}\mathbf{Z}^{N}_{\s,\y}).
\end{align}
\end{definition}
%%%
In words, the modification $\mathbf{D}^{N,1,1}$ introduced in Definition \ref{definition:pc7} cuts off the second term in the RHS of the $\mathbf{D}^{N,1}$ equation so that said term is controlled by a vanishing-in-$N$ term. Because we expect $\t_{\mathbf{Y}}=1$ with high probability, we may that expect this cutoff does nothing with high probability. Indeed, the first step to confirming this prediction is the following identification between $\mathbf{D}^{N,1}$ and $\mathbf{D}^{N,1,1}$ until the time $\t_{\mathbf{Y}}$. The proof is standard uniqueness for linear equations; until $\t_{\mathbf{Y}}$, both $\mathbf{D}^{N,1}$ and $\mathbf{D}^{N,1,1}$ solve the same linear equation, whose solutions are unique even with the random martingale term.
%%%
\begin{lemma}\label{lemma:pc8}
\fsp We have $\mathbf{D}^{N,1}_{\t,\x}=\mathbf{D}^{N,1,1}_{\t,\x}$ for all $\x\in\Z$ and for all $\t\in[0,\t_{\mathbf{Y}}]$.
\end{lemma}
%%%
Lemma \ref{lemma:pc8} will be how we relate $\mathbf{D}^{N,1,1}$ estimates and $\mathbf{D}^{N,1}$ estimates. The following result provides the former estimates; its proof illustrates the utility in the cutoff via $\mathbf{1}(\t\leq\t_{\mathbf{Y}})$. In particular, the second term in the $\mathbf{D}^{N,1,1}$ equation is \emph{deterministically} small, {whereas the} second term in the $\mathbf{D}^{N,1}$ equation is only small in probability, and the estimates we require for $\mathbf{D}^{N,1,1}$ are achieved by higher moment estimates as in \cite{BG,DT}, for which ``high probability" statements are insufficient for technical reasons. We clarify that the first estimate in Lemma \ref{lemma:pc9} below is a moment estimate that holds uniformly in space-time, whereas the second estimate therein is a high probability bound for the exponentially-weighted space-time supremum/``space-time maximal process" of $\mathbf{D}^{N,1,1}$, which is much stronger than the moment estimate below and therefore requires an additional small power of $N$. We will eventually relate the latter ``maximal process" estimate to the former moment estimate via union bound.
%%%
\begin{lemma}\label{lemma:pc9}
\fsp Provided any $\t\in[0,1]$ and $\x\in\Z$ and $p\geq1$, we have $\E|\mathbf{D}^{N,1,1}_{\t,\x}|^{2p}\lesssim_{p}N^{-p\gamma}$, where $\gamma$ is the positive constant from \emph{Definition \ref{definition:pc4}}. Moreover, let us define the following event, in which $\delta$ is the constant in \emph{Lemma \ref{lemma:pc3}}, but we could also choose just any very small constant as long as it is only allowed to depend on and be much smaller than $\beta_{\star}$ and $\gamma/2$:
\begin{align}
\mathcal{E}(\mathbf{D}^{N,1,1}) \ = \ \left\{\|\mathbf{D}^{N,1,1}\|_{1;\Z}\geq N^{\delta}N^{-\gamma/2}\right\}.
\end{align}
Given any positive constant $D$, we have the probability estimate $\mathbf{P}(\mathcal{E}(\mathbf{D}^{N,1,1})) \lesssim_{\delta,D} N^{-D}$.
\end{lemma}
%%%
%%%
\begin{proof}
Similar to the proof of Lemma \ref{lemma:pc3}, it suffices to prove the pointwise moment estimate in the statement of Lemma \ref{lemma:pc9}; this follows by discretizing $\|\|_{1;\Z}$-norms via Lemma \ref{lemma:ste} and a union bound decomposition over $\mathbb{I}^{\mathbf{T}}\times\Z$ as in the proof of Lemma \ref{lemma:pc3} leading up to \eqref{eq:pc36}. To make our presentation more convenient and clear, we rewrite said moment bound below, in which $p\geq1$ and we recall the $2p$-norm $\|\cdot\|_{\omega;2p}=(\E|\cdot|^{2p})^{1/2p}$ with respect to all randomness in the particle system:
\begin{align}
\|\mathbf{D}^{N,1,1}_{\t,\x}\|_{\omega;2p}^{2p} \ \lesssim_{p} \ N^{-p\gamma}. \label{eq:pc91}
\end{align}
First let $\mathbf{D}^{N,1,1,p,\mathbf{X}}_{\t}$ denote the supremum over $\Z$ of the LHS of \eqref{eq:pc91}, namely at time $\t\geq0$. It trivially suffices to obtain \eqref{eq:pc91} uniformly over $\t\geq0$ after replacing the LHS with $\mathbf{D}^{N,1,1,p,\mathbf{X}}_{\t}$. To this end, we employ the stochastic integral equation defining $\mathbf{D}^{N,1,1}$ that was constructed in Definition \ref{definition:pc7}.  Similar to \eqref{eq:pc37} but without the initial data term therein, we have
\begin{align}
\|\mathbf{D}^{N,1,1}_{\t,\x}\|_{\omega;2p}^{2} \ \lesssim \ \|\mathbf{H}^{\beta_{\star},\mathbf{T}}_{\t,\x}(\mathbf{D}^{N,1,1}\d\xi^{N})\|_{\omega;2p}^{2} + \|\mathbf{H}^{\beta_{\star},\mathbf{T}}_{\t,\x}(\mathbf{1}_{\s\leq\t_{\mathbf{Y}}}2^{-1}c_{N}\mathbf{1}_{\y=0}\mathfrak{q}_{\s,\y}^{\mathrm{tot}}\mathbf{Z}^{N}_{\s,\y})\|_{\omega;2p}^{2}. \label{eq:pc92}
\end{align}
We treat the first term on the RHS of \eqref{eq:pc92} similar to the second term on the RHS of \eqref{eq:pc37}; indeed, our analysis for that term only required heat kernel estimates for $\mathbf{H}^{\beta_{\star}}$ and the general martingale inequality of Lemma \ref{lemma:mge}. We defer this momentarily. First, we note the second term on the RHS of \eqref{eq:pc92} is uniformly bounded above by $N^{-\gamma}$, as the heat operator itself is uniformly bounded above by $N^{-\gamma/2}$ by definition of $\t_{\mathbf{Y}}$ in Definition \ref{definition:pc4}. Therefore, using \eqref{eq:pc38}, \eqref{eq:pc39}, \eqref{eq:pc310}, and \eqref{eq:pc311} but for $\mathbf{D}^{N,1,1}$ in place of $\mathbf{Y}^{N}$ (noting $\mathbf{D}^{N,1,1}$ has zero initial data) along with the aforementioned estimate for the second term on the RHS of \eqref{eq:pc92}, we deduce, with any choice of small but fixed and positive $\varrho$, the following estimate:
\begin{align}
\|\mathbf{D}^{N,1,1}_{\t,\x}\|_{\omega;2p}^{2} \ \lesssim_{p} \ N^{-\gamma} + N^{-1+9999\beta_{\star}+\varrho}\int_{0}^{\t}(\mathbf{D}_{\s}^{N,1,1,p,\mathbf{X}})^{1/p}|\t-\s|^{-1+\varrho'}\d\s + \int_{0}^{\t}(\mathbf{D}_{\s}^{N,1,1,p,\mathbf{X}})^{1/p}|\t-\s|^{-1/2}\d\s. \label{eq:pc93}
\end{align}
Application of the Gronwall inequality, as in what we did in the proof of Lemma \ref{lemma:pc3} after \eqref{eq:pc312}, proves the desired moment bound \eqref{eq:pc91}, so this completes the proof entirely as we noted in the paragraph preceding \eqref{eq:pc91}.
\end{proof}
%%%
We are now in position to prove the $\mathcal{E}(1)$ estimate in Proposition \ref{prop:pc2} in spirit of the remark prior to Definition \ref{definition:pc7}.
%%%
\begin{lemma}\label{lemma:pc10}
\fsp There exists a positive constant $\gamma_{3}$ depending only on $\beta_{\star}$ such that $\mathbf{P}(\t_{\mathbf{Y}}\neq1)\lesssim N^{-\gamma_{3}}$.
\end{lemma}
%%%
%%%
\begin{proof}
Observe the event $\{\t_{\mathbf{Y}}\neq1\}$ implies that $\t_{1}\neq1$ or $\t_{2}\neq1$ or $\t_{3}\neq1$. A somewhat more refined version of this observation leads to the union bound below, which follows by conditioning on $\t_{1}\neq1$ or $\t_{1}=1$, and then conditioning on $\t_{\mathbf{Y}}=\t_{2}$ or $\t_{\mathbf{Y}}=\t_{3}$:
\begin{align}
\mathbf{P}(\t_{\mathbf{Y}}\neq1) \ \leq \ \mathbf{P}(\t_{1}\neq1)+\mathbf{P}(\t_{1}=1,\t_{2}\leq\t_{3},\t_{2}\neq1)+\mathbf{P}(\t_{1}=1,\t_{3}\leq\t_{2},\t_{3}\neq1). \label{eq:pc101}
\end{align}
By Lemma \ref{lemma:pc3} and the a priori $\mathbf{Y}^{N}$ estimates therein, the first probability on the RHS of \eqref{eq:pc101} is at most $N^{-100}$ times a uniformly bounded factor. Therefore, it suffices to estimate the last two probabilities on the RHS of \eqref{eq:pc101}. First, by Lemma \ref{lemma:pc9} and by Lemma \ref{lemma:ste}, we may condition on the event in which both of the estimates below hold; we note that the following short-time estimates are true uniformly over short-time scales $\t_{N}\in[0,N^{-100}]$ and over $\t\in[0,1]$ as well, and that $\delta\geq0$ is a very small multiple of $\gamma$:
\begin{align}
\|\mathbf{D}^{N,1,1}\|_{1;\Z} \ \leq \ N^{\delta-\gamma/2} \and \|\mathbf{Z}^{N}\|_{\t;\Z}^{-1}\|\grad_{\t_{N}}^{\mathbf{T}}\mathbf{Z}^{N}\|_{\t;\Z} + \|\mathbf{Y}^{N}\|_{\t;\Z}^{-1}\|\grad_{\t_{N}}^{\mathbf{T}}\mathbf{Y}^{N}\|_{\t;\Z} \ \lesssim \ N^{-1/3}. \label{eq:pc102}
\end{align}
Indeed, Lemma \ref{lemma:ste} guarantees that \eqref{eq:pc102} holds, uniformly over allowable $\t_{N}$ and $\t$, outside an event of probability at most $N^{-100}$. We may also condition on the following, whose failure holds with sufficiently small probability by Proposition \ref{prop:lpe3} for $\t_{\mathrm{st}}=\t_{3}$:
\begin{align}
\|\mathbf{H}^{\beta_{\star},\mathbf{T}}_{\t,\x}(2^{-1}c_{N}\mathbf{1}_{\y=0}\mathfrak{q}_{\s,\y}^{\mathrm{tot}}\mathbf{Z}^{N}_{\s,\y})\|_{\t_{3};\Z;\max} \ \leq \ N^{-3\gamma/4}\|\mathbf{Z}^{N}\|_{\t_{3};\Z}. \label{eq:pc102b}
\end{align}
We now show that conditioning on \eqref{eq:pc102} and \eqref{eq:pc102b}, simultaneously/uniformly over allowable $\t_{N}$ and $\t$, the events in the second and third probabilities on the RHS of \eqref{eq:pc101} are empty; this would finish the proof of Lemma \ref{lemma:pc10} provided \eqref{eq:pc101}. To this end, we consider each of these two events separately in the following pair of bullet points.
%%%
\begin{itemize}
\item Consider first the event $\t_{1}=1$ and $\t_{2}\leq\t_{3}$ and $\t_{2}\neq1$, conditioning on \eqref{eq:pc102} to hold simultaneously in both $\t_{N}\in[0,N^{-100}]$ and $\t\in[0,1]$. As $\t_{2}\in[0,1]$ with probability 1 by construction, we know $\t_{2}<1$. Now, we take $\t^{\mathrm{short}}\in(0,N^{-100}]$ such that $\t_{2}+\t^{\mathrm{short}}\leq1$. Recall $\mathbf{D}^{N,1}=\mathbf{Z}^{N}-\mathbf{Y}^{N}$. We first have the following elementary consequence of the triangle inequality:
\begin{align}
\|\mathbf{D}^{N,1}\|_{\t_{2}+\t^{\mathrm{short}};\Z} \ \leq \ \|\mathbf{D}^{N,1}\|_{\t_{2};\Z} + {\sup}_{0\leq\s\leq\t^{\mathrm{short}}}\|\grad_{\s}^{\mathbf{T}}\mathbf{D}^{N,1}\|_{\t_{2};\Z}. \label{eq:pc103}
\end{align}
The second term on the RHS is controlled by time-gradients of both $\mathbf{Z}^{N}$ and $\mathbf{Y}^{N}$ by definition of $\mathbf{D}^{N,1}$, the triangle inequality, and linearity of $\grad^{\mathbf{T}}$ operators. On the other hand, because $\t_{2}=\t_{\mathbf{Y}}$ on the event considered in this bullet point, we obtain the identification $\mathbf{D}^{N,1}=\mathbf{D}^{N,1,1}$ in the first norm on the RHS of \eqref{eq:pc103}; this follows from Lemma \ref{lemma:pc8}, and we therefore control the first term on the RHS via \eqref{eq:pc102}. Ultimately, we get the following estimate:
\begin{align}
\|\mathbf{D}^{N,1}\|_{\t_{2}+\t^{\mathrm{short}};\Z} \ \lesssim \ N^{\delta-\gamma/2} + N^{-1/3}\|\mathbf{Y}^{N}\|_{\t_{2};\Z} + N^{-1/3}\|\mathbf{Z}^{N}\|_{\t_{2};\Z}. \label{eq:pc104}
\end{align}
Because $\t_{1}=1$, the second term on the RHS of \eqref{eq:pc104} is controlled by $N^{-1/3+\delta}$, for example. On the other hand, we know the final norm on the RHS of \eqref{eq:pc104} is controlled by that of $\mathbf{Y}^{N}$ and that of $\mathbf{D}^{N,1}$. The former of these is again controlled by $N^{\delta}$, whereas the latter of these is controlled by 1, for example, courtesy of construction of $\t_{2}$ in Definition \ref{definition:pc4}. Thus, we also deduce a bound of $N^{-1/3+\delta}$ for this last term on the RHS of \eqref{eq:pc104}. Plugging the two bounds from this paragraph into \eqref{eq:pc104} implies that the LHS of \eqref{eq:pc104} is bounded above by a uniformly bounded factor times $N^{\delta-\gamma/2}$, since $\gamma$ is positive and very small. Because $\gamma_{2}$ is a very small multiple of $\gamma$, where $\gamma_{2}$ is from Proposition \ref{prop:pc2} and Definition \ref{definition:pc4}, we deduce that the LHS of \eqref{eq:pc104} is controlled by $N^{\delta-\gamma_{2}}$, with the extra $N^{-\gamma/2+\gamma_{2}}$ accounting for the aforementioned uniformly bounded factor provided that $N$ is sufficiently large depending only on the slow bond parameter $\beta_{\star}$. This contradicts the definition of $\t_{2}$; we just deduced that we can push $\t_{2}$ forward to any time until $\t_{2}+\t^{\mathrm{short}}$, and the lower bound defining $\t_{2}$ would not have been achieved as long as we condition on \eqref{eq:pc102}, because until $\t_{2}=\t_{\mathbf{Y}}$, said lower bound is far from true, and short times cannot change things much. Therefore, the event on which $\t_{1}=1$ and $\t_{\mathbf{Y}}=\t_{2}\leq\t_{3}$ in question in this bullet point is empty.
\item Let us now move to the event in the last probability on the RHS of \eqref{eq:pc101}, this time conditioning on \eqref{eq:pc102b}. Because $\t_{3}\leq\t_{2}$ in this bullet point, we know that until time $\t_{3}$, we have control on $\|\mathbf{Z}^{N}\|_{\t_{3};\Z}$ via control on $\|\mathbf{Y}^{N}\|_{\t_{3};\Z}$ and $\|\mathbf{D}^{N,1}\|_{\t_{3};\Z}$ via $\t_{1}$ and $\t_{2}$, respectively. In particular, using the triangle inequality and by $\t_{3}\leq\t_{1},\t_{2}$, we have the bound $\|\mathbf{Z}^{N}\|_{\t_{3};\Z}\lesssim N^{\delta}$. Thus, by \eqref{eq:pc102b},
\begin{align}
\|\mathbf{H}^{\beta_{\star},\mathbf{T}}_{\t,\x}(2^{-1}c_{N}\mathbf{1}_{\y=0}\mathfrak{q}_{\s,\y}^{\mathrm{tot}}\mathbf{Z}^{N}_{\s,\y})\|_{\t_{3};\Z;\max} \ \leq \ N^{-3\gamma/4}(1+\|\mathbf{Z}^{N}\|_{\t_{3};\Z}^{2}) \ \lesssim \ N^{-3\gamma/4+2\delta}. \label{eq:pc105}
\end{align}
With $\delta$ small so that $-3\gamma/4+2\delta\leq-2\gamma/3$, we deduce the LHS of \eqref{eq:pc105} is bounded above strictly by $N^{-\gamma/2}$, with the implied constants/factors being accounted for by $N^{-2\gamma/3+\gamma/2}$ if $N\geq0$ is sufficiently large but again depending only on the slow bond parameter $\beta_{\star}$, like with the previous bullet point. Because space-time heat operators $\mathbf{H}^{\beta_{\star},\mathbf{T}}$ are uniformly continuous in space-time, we may push $\t_{3}$ on the LHS of \eqref{eq:pc105} forward slightly such that the LHS of \eqref{eq:pc105} is still controlled above by $N^{-\gamma/2}$. Again, like the previous bullet point, the amount we push forward past $\t_{3}$ may depend heavily on $N$, though just any positive amount is sufficient. This shows that conditioning on \eqref{eq:pc102b}, the event in the last probability on the RHS of \eqref{eq:pc101} is empty.
\end{itemize}
%%%
Combining the previous two bullet points with \eqref{eq:pc101} and the paragraph after \eqref{eq:pc101} completes the proof.
\end{proof}
%%%
%%%
\subsection{$\mathcal{E}(2)$ Estimate}
%%%
The proof of the $\mathcal{E}(2)$ estimate in Proposition \ref{prop:pc2} will follow from a similar consideration as in the previous subsection to control $\mathcal{E}(1)$. However, at a technical level, many details are different; we must now replace heat kernels $\mathbf{H}^{\beta_{\star}}$ with $\mathbf{H}$, so in this subsection we account for the difference in power-counting and details in the heat operators after making this change in heat kernels. As noted after Proposition \ref{prop:pc2}, however, we have the benefit of a priori estimates for the two objects $\mathbf{Y}^{N}$ and $\mathbf{W}^{N}$ we are comparing in $\mathcal{E}(2)$, so no cutoffs of type $\mathbf{D}^{N,1,1}$ from Definition \ref{definition:pc7}, or even the random times in Definition \ref{definition:pc4}, are required here. We will again introduce convenient notation as before.
%%%
\begin{definition}\label{definition:pc11}
\fsp Define $\mathbf{D}^{N,2}=\mathbf{Y}^{N}-\mathbf{W}^{N}$ that, by construction in Definition \ref{definition:pc1}, solves the stochastic integral equation below where $\wt{\mathbf{H}}^{\beta_{\star},\mathbf{T}}=\mathbf{H}^{\beta_{\star},\mathbf{T}}-\mathbf{H}^{\mathbf{T}}$ is a difference of space-time heat operators with different heat kernels, extending the heat-kernel-difference construction in Definition \ref{definition:hke6} to heat operators:
\begin{align}
\mathbf{D}^{N,2}_{\t,\x} \ = \ {\mathbf{H}^{\beta_{\star},\mathbf{X}}_{\t,\x}(\mathbf{Z}^{N})-\mathbf{H}^{\mathbf{X}}_{\t,\x}(\mathbf{Z}^{N})}+\mathbf{H}_{\t,\x}^{\mathbf{T}}(\mathbf{D}^{N,2}\d\xi^{N}) + \wt{\mathbf{H}}_{\t,\x}^{\beta_{\star},\mathbf{T}}(\mathbf{Y}^{N}\d\xi^{N}).
\end{align}
\end{definition}
%%%
Although we could certainly apply the same strategy used in the previous section for the $\mathcal{E}(1)$ estimate given by cutting off the second term on the RHS of the $\mathbf{D}^{N,2}$ equation with high probability, it turns out the moment estimate strategy for Lemma \ref{lemma:pc9} can be applied directly to control $\mathbf{D}^{N,2}$. This provides, ultimately, the following analog of Lemma \ref{lemma:pc9} but for $\mathbf{D}^{N,2}$ in place of $\mathbf{D}^{N,1,1}$.
%%%
\begin{lemma}\label{lemma:pc12}
\fsp There exists a positive and universal constant $\alpha\geq999\delta$ such that for any $\t\in[0,1]$ and $\x\in\Z$ and $p\geq1$, we have the moment bound $\E|\mathbf{D}^{N,2}_{\t,\x}|^{2p}\lesssim_{p} N^{-2p\alpha}$; here, we pick the choice of $\delta$ from \emph{Lemma \ref{lemma:pc3}}. Similar to \emph{Lemma \ref{lemma:pc9}}, define the event
\begin{align}
\mathcal{E}(\mathbf{D}^{N,2})  = \ \left\{\|\mathbf{D}^{N,2}\|_{1;\Z}\geq N^{2\delta}N^{-\alpha}\right\}.
\end{align}
Given any positive constant $D$, we have the probability estimate $\mathbf{P}(\mathcal{E}(\mathbf{D}^{N,2}))\lesssim_{\delta,D}N^{-D}$.
\end{lemma}
%%%
%%%
\begin{proof}
We first note that Lemma \ref{lemma:ste} applies to both $\mathbf{Y}^{N}$ and $\mathbf{W}^{N}$. In particular, it suffices to bound $\mathbf{D}^{N,2}$ uniformly in space-time with the sub-exponential weight in $\|\|_{1;\Z}$ after we discretize time $[0,1]$ into $\mathbb{I}^{\mathbf{T}}=\{\mathfrak{j}N^{-100}\}_{\mathfrak{j}=0,\ldots,N^{100}}$. Indeed, the supremum of $\mathbf{D}^{N,2}$ over $[0,1]\times\Z$ with the sub-exponential weight in $\|\|_{1;\Z}$ is controlled by the supremum of the same object but over $\mathbb{I}^{\mathbf{T}}\times\Z$ plus the supremum of scale-$N^{-100}$ time-gradients of $\mathbf{D}^{N,2}$ times the same sub-exponential weight. However, the time-gradients of $\mathbf{D}^{N,2}$ are controlled by time-gradients of $\mathbf{Y}^{N}$ and $\mathbf{W}^{N}$, respectively, by construction in Definition \ref{definition:pc11}. These are controlled by Lemma \ref{lemma:ste} combined with the high probability bounds for $\mathbf{Y}^{N}$ and $\mathbf{W}^{N}$ of $N^{\delta}$ from Lemma \ref{lemma:pc3} in terms of $N^{2\delta-1/2}$, which is certainly lower-order than $N^{2\delta-\alpha}$ from $\mathcal{E}(\mathbf{D}^{N,2})$. In particular, we deduce that it suffices to establish the stated moment bound in Lemma \ref{lemma:pc12}, from which the $\mathcal{E}(\mathbf{D}^{N,2})$-probability estimate follows by the discretization strategy described in this paragraph. For this, we again appeal to the integral equation satisfied by $\mathbf{D}^{N,2}$ from Definition \ref{definition:pc11}. To start, we have
\begin{align}
\|\mathbf{D}^{N,2}_{\t,\x}\|_{\omega;2p}^{2} \ &\lesssim \ {\|\mathbf{H}^{\beta_{\star},\mathbf{X}}_{\t,\x}(\mathbf{Z}^{N})-\mathbf{H}^{\mathbf{X}}_{\t,\x}(\mathbf{Z}^{N})\|_{\omega;2p}^{2}}+\|\mathbf{H}_{\t,\x}^{\mathbf{T}}(\mathbf{D}^{N,2}\d\xi^{N})\|_{\omega;2p}^{2} + \|\wt{\mathbf{H}}_{\t,\x}^{\beta_{\star},\mathbf{T}}(\mathbf{Y}^{N}\d\xi^{N})\|_{\omega;2p}^{2}. \label{eq:pc121}
\end{align}
{Let us bound the first term on the RHS of \eqref{eq:pc121}. Note that $\|\mathbf{Z}^{N}_{0,\x}\|_{\omega;2p}\lesssim_{p}1$ uniformly in $\x$. Indeed, this bound holds for $\mathbf{G}^{N}$; see Definition \ref{definition:ns} {for $\kappa=0$ as assumed in Theorem \ref{theorem:kpz}}. Moreover, $\mathbf{G}^{N}$ and $\mathbf{Z}^{N}$ have the same initial data; see Definition \ref{definition:mshe4}. Using this with Proposition \ref{prop:hke11} for $\kappa=1$, Cauchy-Schwarz, the fact that $\sum_{\y\in\Z}|\wt{\mathbf{H}}^{\beta_{\star}}_{0,\t,\x,\y}|\leq\sum_{\y\in\Z}\mathbf{H}^{\beta_{\star}}_{0,\t,\x,\y}+\sum_{\y\in\Z}\mathbf{H}_{0,\t,\x,\y}\lesssim1$, and an elementary sub-exponential series estimate gives
\begin{align}
\|\mathbf{H}^{\beta_{\star},\mathbf{X}}_{\t,\x}(\mathbf{Z}^{N})-\mathbf{H}^{\mathbf{X}}_{\t,\x}(\mathbf{Z}^{N})\|_{\omega;2p}^{2} \ &= \ \|\sum_{\y\in\Z}\wt{\mathbf{H}}^{\beta_{\star}}_{0,\t,\x,\y}\mathbf{Z}^{N}_{0,\y}\|_{\omega;2p}^{2} \ \lesssim \ \sum_{\y\in\Z}|\wt{\mathbf{H}}^{\beta_{\star}}_{0,\t,\x,\y}| \label{eq:edit1}\\
&\lesssim \ \sum_{\y\in\Z}N^{-1-\e_{1}}|\t|^{-1+\e_{2}}\exp\left(-\frac{|\x-\y|}{(N\t^{1/2})\vee1}\right) \ \lesssim \ N^{-\e_{1}}|\t|^{-\frac12+\e_{2}}.\label{eq:edit2}
\end{align}
(Briefly, summing the sub-exponential factor takes a factor of $N^{-1}\t^{-1/2}$ per standard heat kernel scaling. This is what gives the last estimate; see the proof of Corollary A.2 in \cite{DT} for a similar bound.) On the other hand, by adding/subtracting $\mathbf{Z}^{N}_{0,\x}$ and using the triangle inequality, we have
\begin{align}
\|\mathbf{H}^{\beta_{\star},\mathbf{X}}_{\t,\x}(\mathbf{Z}^{N})-\mathbf{H}^{\mathbf{X}}_{\t,\x}(\mathbf{Z}^{N})\|_{\omega;2p}^{2} \ &\lesssim \ \|\mathbf{H}^{\beta_{\star},\mathbf{X}}_{\t,\x}(\mathbf{Z}^{N})-\mathbf{Z}_{0,\x}^{N}\|_{\omega;2p}^{2}+\|\mathbf{Z}_{0,\x}^{N}-\mathbf{H}^{\mathbf{X}}_{\t,\x}(\mathbf{Z}^{N})\|_{\omega;2p}^{2}.\label{eq:edit3}
\end{align}
Our analysis for the second term on the RHS of \eqref{eq:edit3} will be the same as that for the first term therein, so we focus on said first term. Because the heat kernel defining $\mathbf{H}^{\beta_{\star},\mathbf{X}}$ is a probability measure on $\Z$ with respect to the forward spatial variable, we have 
\begin{align}
\mathbf{H}^{\beta_{\star},\mathbf{X}}_{\t,\x}(\mathbf{Z}^{N})-\mathbf{Z}_{0,\x}^{N} \ = \ \sum_{\y\in\Z}\mathbf{H}^{\beta_{\star}}_{0,\t,\x,\y}(\mathbf{Z}^{N}_{0,\y}-\mathbf{Z}^{N}_{0,\x}).
\end{align}
Note $\|\mathbf{Z}^{N}_{0,\y}-\mathbf{Z}^{N}_{0,\x}\|_{\omega;2p}^{2p}\lesssim_{\upsilon,p} N^{-2p\upsilon}|\x-\y|^{2p\upsilon}$ for any $\upsilon\in(0,1/2)$. Again, this is true for $\mathbf{G}^{N}$, and the initial data for $\mathbf{Z}^{N}$ and $\mathbf{G}^{N}$ are the same. As with the proof of \eqref{eq:edit1}-\eqref{eq:edit2}, by using this with Proposition \ref{prop:hke11}, the previous display, Cauchy-Schwarz, and an elementary series bound, we get
\begin{align}
\|\mathbf{H}^{\beta_{\star},\mathbf{X}}_{\t,\x}(\mathbf{Z}^{N})-\mathbf{Z}_{0,\x}^{N}\|_{\omega;2p}^{2} \ \lesssim_{\upsilon,p} \ \sum_{\y\in\Z}|\mathbf{H}_{0,\t,\x,\y}^{\beta_{\star}}|\cdot N^{-2\upsilon}|\x-\y|^{2\upsilon} \ \lesssim \ \t^{\upsilon}. \label{eq:edit4}
\end{align}
(The last bound is just the usual $2\upsilon$-moment bound for a symmetric simple random walk of speed $\mathrm{O}(N^{2})$.) Again, the estimate \eqref{eq:edit4} also holds upon replacing its LHS by the last term in \eqref{eq:edit3}. We deduce the following from this, \eqref{eq:edit1}-\eqref{eq:edit2}, and \eqref{eq:edit3}:
\begin{align}
\|\mathbf{H}^{\beta_{\star},\mathbf{X}}_{\t,\x}(\mathbf{Z}^{N})-\mathbf{H}^{\mathbf{X}}_{\t,\x}(\mathbf{Z}^{N})\|_{\omega;2p}^{2} \ \lesssim \ \min\left(N^{-\e_{1}}|\t|^{-\frac12+\e_{2}},\t^{\upsilon}\right) \quad\mathrm{for} \ \mathrm{all} \ \t\geq0 \ \mathrm{and} \ \upsilon\in(0,1/2).
\end{align}
Fix $\upsilon=1/3$. Choosing $\t$ to be a sufficiently small power of $N^{-1}$ implies that the first term on the RHS of \eqref{eq:pc121} satisfies the moment bound that we stated for $\mathbf{D}^{N,2}$ in the statement of Lemma \ref{lemma:pc12}.
} 

{It remains to similarly bound the last two} terms on the RHS of \eqref{eq:pc121}, starting with the first {of these}. We again use \eqref{eq:pc38}, \eqref{eq:pc39}, \eqref{eq:pc310}, and \eqref{eq:pc311} for the first term on the RHS of \eqref{eq:pc121} but after replacing $\mathbf{Y}^{N}$ in \eqref{eq:pc38}, \eqref{eq:pc39}, \eqref{eq:pc310}, and \eqref{eq:pc311} by $\mathbf{D}^{N,2}$; we reemphasize that \eqref{eq:pc38}, \eqref{eq:pc39}, \eqref{eq:pc310}, and \eqref{eq:pc311} applies to $\mathbf{D}^{N,2}$ as well, because we only depend on heat kernel estimates for $\mathbf{H}^{\beta_{\star}}$ and Lemma \ref{lemma:mge}. Because the initial data of $\mathbf{D}^{N,2}$ vanishes, we get, with $\mathbf{D}^{N,2,p,\mathbf{X}}_{\t}$ the supremum over $\x\in\Z$ of $\|\mathbf{D}_{\t,\x}^{N,2}\|_{\omega;2p}^{2p}$ and with $\varrho>0$ small but fixed,
\begin{align}
\|\mathbf{H}_{\t,\x}^{\mathbf{T}}(\mathbf{D}^{N,2}\d\xi^{N})\|_{\omega;2p}^{2} \ \lesssim \ N^{-1+9999\beta_{\star}+\varrho}\int_{0}^{\t}(\mathbf{D}_{\s}^{N,2,p,\mathbf{X}})^{1/p}|\t-\s|^{-1+\varrho'}\d\s + \int_{0}^{\t}(\mathbf{D}_{\s}^{N,2,p,\mathbf{X}})^{1/p}|\t-\s|^{-1/2}\d\s. \label{eq:pc122}
\end{align}
For the second term on the RHS of \eqref{eq:pc121}, we again employ \eqref{eq:pc38}, \eqref{eq:pc39}, \eqref{eq:pc310}, and \eqref{eq:pc311} but with the following modifications. First, all of the $\mathbf{Y}^{N,p,\mathbf{X}}$-terms in these estimates are uniformly bounded by $p$-dependent constants, because we can now apply the conclusion of Lemma \ref{lemma:pc3} that \eqref{eq:pc38}, \eqref{eq:pc39}, \eqref{eq:pc310}, and \eqref{eq:pc311} were used to prove. Moreover, by appealing to Proposition \ref{prop:hke7}, for $|\y|\geq N^{999\beta_{\star}}$, we obtain extra $N^{-\e_{1}}$-factors for the difference $\wt{\mathbf{H}}^{\beta_{\star}}$ of heat kernels $\mathbf{H}^{\beta_{\star}}$ and $\mathbf{H}$, if we allow for a slightly more singular but still integrable singularity $|\t-\s|^{-1+\e_{2}}$; here, both $\e_{1}$ and $\e_{2}$ are positive. This observation changes only our analysis of $\Phi_{2,\t,\x}$ in \eqref{eq:pc39} because only this terms deals with $|\y|\geq N^{999\beta_{\star}}$; it provides an additional $N^{-\e_{1}}$ factor with the aforementioned change in time-integrable singularity. Ultimately, we deduce the following, again because $\mathbf{Y}^{N}$ is uniformly bounded in moments:
\begin{align}
 \|\wt{\mathbf{H}}_{\t,\x}^{\beta_{\star},\mathbf{T}}(\mathbf{Y}^{N}\d\xi^{N})\|_{\omega;2p}^{2} \ &\lesssim_{p} \ N^{-1+9999\beta_{\star}+\varrho}\int_{0}^{\t}|\t-\s|^{-1+\varrho'}\d\s + N^{-\e_{1}} \int_{0}^{\t}|\t-\s|^{-1+\e_{2}}\d\s \\
 &\lesssim_{\varrho',\e_{2}} \ N^{-1+9999\beta_{\star}+\varrho} + N^{-\e_{1}}. \label{eq:pc123}
\end{align}
We now combine \eqref{eq:pc121}, \eqref{eq:pc122}, and \eqref{eq:pc123} with $\varrho>0$ small but fixed  along with the Gronwall inequality that we used in the proofs of Lemma \ref{lemma:pc3} and Lemma \ref{lemma:pc9} in order to deduce the stated $\mathbf{D}^{N,2}$ pointwise moment estimate in the statement of Lemma \ref{lemma:pc12}. As noted prior to \eqref{eq:pc121}/in the first paragraph of this proof, this completes the proof of Lemma \ref{lemma:pc12}, so we are done.
\end{proof}
%%%
%%%
\subsection{Proof of Proposition \ref{prop:pc2}}
%%%
As noted prior to Definition \ref{definition:pc7}, the probability of $\mathcal{E}(1)$ in Proposition \ref{prop:pc2} is controlled by the probability of $\mathfrak{t}_{\mathbf{Y}}\neq1$, which is controlled by Lemma \ref{lemma:pc10}. Moreover, the probability of $\mathcal{E}(2)$ in Proposition \ref{prop:pc2} is controlled by that of $\mathcal{E}(\mathbf{D}^{N,2})$ in Lemma \ref{lemma:pc12} if we adjust $\alpha$ by subtracting the much smaller $2\delta$. This completes the proof. \qed
%
%
%
%%%
\section{Proof of Theorem \ref{theorem:kpz}}\label{section:proof}
%%%
We first refer to Remark \ref{remark:mshe6}, which says it is enough to prove Theorem \ref{theorem:kpz} for $\mathbf{Z}^{N}$ in place of the original Gartner transform $\mathbf{G}^{N}$. Because the space-time-locally-uniform topology is stronger than the Skorokhod topology on $\mathscr{D}(\R_{\geq0},\mathscr{C}(\R))$ \cite{Bil}, by Proposition \ref{prop:pc2} it suffices to prove Theorem \ref{theorem:kpz} except replacing $\mathbf{Z}^{N}$ therein with $\mathbf{W}^{N}$ from Definition \ref{definition:pc1}. Indeed, we write $\mathbf{Z}^{N}=\mathbf{W}^{N}+(\mathbf{Z}^{N}-\mathbf{Y}^{N})+(\mathbf{Y}^{N}-\mathbf{W}^{N})$. The second and third terms in this decomposition for $\mathbf{Z}^{N}$ vanish in the space-time-locally-uniform topology, even after $\Gamma^{N}$-rescaling operators, by Proposition \ref{prop:pc2}.  Therefore, if $\Gamma^{N}\mathbf{W}^{N}$ converges to the $\mathrm{SHE}$ in Theorem \ref{theorem:kpz}, then $\mathbf{Z}^{N}$ is equal to something that converges to $\mathrm{SHE}$ in Theorem \ref{theorem:kpz} plus terms that vanish in probability with respect to the topology the aforementioned convergences holds in, at which point standard probability on complete separable metric spaces allows us to deduce Theorem \ref{theorem:kpz} as written. In particular, it suffices to prove the following result, for which we can basically follow the proof of Theorem 1.1 in \cite{DT} except for a few remarks that we highlight afterwards.
%%%
\begin{prop}\label{prop:kpz1}
\fsp The sequence $\Gamma^{N}\mathbf{W}^{N}$ converges to the solution of $\mathrm{SHE}$ in the Skorokhod space of \emph{Theorem \ref{theorem:kpz}} with initial data given by $\Gamma^{\infty,\mathbf{X}}\mathbf{G}^{\infty}$ from \emph{Definition \ref{definition:ns}}.
\end{prop}
%%%
It is easy to see that by construction of $\mathbf{W}^{N}$ in Definition \ref{definition:pc1} that $\Gamma^{N,\mathbf{X}}\mathbf{W}^{N}$ converges to $\Gamma^{\infty,\mathbf{X}}\mathbf{G}^{\infty}$ since the $\mathbf{W}^{N}$ initial data is given by $\mathbf{Z}^{N}$/$\mathbf{G}^{N}$ initial data. Moreover, observe $\mathbf{W}^{N}$ solves an equation like the Gartner transform equation in \cite{BG} or Section 2 of \cite{DT} but without the additional ``gradient terms" and ``weakly vanishing" terms therein. In particular, we may now follow the proof of Theorem 1.1 in \cite{DT}. The tightness required in Proposition \ref{prop:kpz1} follows from the proof of Proposition 3.2 and Corollary 3.3 in \cite{DT} because this only requires the structure of the SPDE that $\mathbf{W}^{N}$ in this case solves. To identify limit points, we follow Section 3.2 of \cite{DT}. The key issue is deriving the limit of the multiplicative noise $\mathbf{W}^{N}\d\xi^{N}$, so to speak. Even though the particle system dynamics herein are different from \cite{BG} or simple-\cite{DT} models, the difference in particle system dynamics is given by a single slow bond. Because $\mathbf{W}^{N}\d\xi^{N}$ at a single point does not change the limit of $\mathbf{W}^{N}\d\xi^{N}$, which can be seen in the proof of Proposition 1.5 in Section 3.2 of \cite{DT}, for example, we are in the position of Section 3.2 of \cite{DT}. Following the analysis therein, we are left with proving the following hydrodynamic limit that identifies the bracket process of $\mathbf{W}^{N}\d\xi^{N}$ in a weak space-time integrated sense.
%%%
\begin{lemma}\label{lemma:kpz2}
\fsp Consider any function $\mathfrak{w}_{\s,\y}=\prod_{\mathfrak{j}=1}^{\mathfrak{m}}\eta_{\s,\y_{\mathfrak{j}}}$, where $\y_{1},\ldots,{\y_{\mathfrak{m}}}$ are of the form $\y+i_{\mathfrak{j}}$ for distinct $i_{1},\ldots,i_{\mathfrak{m}}$ independent of $\s,\y$. For any compactly supported $\phi\in\mathscr{C}^{\infty}_{\mathrm{c}}(\R)$ and any $0\leq\t\leq1$ fixed, we have the following for $\alpha=1,2$:
\begin{align}
\lim_{N\to\infty}\int_{0}^{\t}N^{-1}{\sum}_{\x\in\Z}\phi_{N^{-1}\x}\mathfrak{w}_{\s,\x}(\mathbf{W}_{\s,\x}^{N})^{\alpha}\d\s \ = \ 0.
\end{align}
\end{lemma}
%%%
%%%
\begin{proof}
Following the proof of Lemma 2.5 in \cite{DT}, it suffices to replace $\mathfrak{w}_{\s,\x}$ in the integral in question by a gradient of \emph{the Gartner transform} $\mathbf{G}^{N}$; relevance of $\mathbf{G}^{N}$, even for an estimate without $\mathbf{G}^{N}$ factors a priori, is the point that we want to emphasize in this argument. By Lemma \ref{lemma:mshe5} and Proposition \ref{prop:pc2}, we can estimate gradients of $\mathbf{G}^{N}$ by gradients of $\mathbf{W}^{N}$ up to errors that vanish in the large-$N$ limit. Inspecting the proof of Lemma 2.5 in \cite{DT}, this is enough to follow said proof upon replacing all Gartner transforms therein with $\mathbf{W}^{N}$, even for the gradient estimate for $\mathbf{G}^{N}$ that becomes relevant. Note the norm used in Lemma \ref{lemma:mshe5} and Proposition \ref{prop:pc2} gives locally uniform vanishing of $\mathbf{G}^{N}-\mathbf{W}^{N}$ in probability on macroscopic compact sets, so the aforementioned replacements are justifiable in the topology of vanishing claimed in the statement of the lemma.
\end{proof}
%%%
%
%
%
%%%
\section{Proof of Proposition \ref{prop:lpe3}}\label{section:lpeproof}
%%%
The rest of this paper is dedicated to proving the key estimate in Proposition \ref{prop:lpe3}. For better organization, we provide the main steps in the following list of bullet points, which the reader should feel invited to treat as an outline for this section to refer back to whenever desired. The main ideas were explained in Section \ref{section:mshe}, but we provide the actual steps needed below.
%%%
\begin{itemize}
\item The first step we take is a cutoff, or equivalently smoothing, procedure of the space-time heat operator. In particular, we smooth out the short-time singularity in the heat kernel $\mathbf{H}^{\beta_{\star}}$ from Proposition \ref{prop:hke3}. This is mostly for technical convenience; we estimate the error in this smoothing directly via Proposition \ref{prop:hke3} and elementary deterministic integral calculations. Again, by employing deterministic integral calculations and the heat kernel estimate in Proposition \ref{prop:hke3}, we will forget the order $N^{-1/2}$ term in $\mathfrak{q}^{\mathrm{tot}}$; see Proposition \ref{prop:mshe1} for the construction of $\mathfrak{q}^{\mathrm{tot}}$. Thus, we will only concern ourselves with $\mathfrak{q}$ from $\mathfrak{q}^{\mathrm{tot}}$ for the rest of this outline.
\item The next step that we take is an entropy production estimate and local equilibrium estimate. In short, this will reduce our analysis of local functionals of the particle system to calculations for stationary particle systems. It depends on standard entropy-energy-type duality \cite{GPV,KL}, which is somewhat weakened by the slow bond effect that increases the needed time for the particle system to relax to its invariant measure. It then leverages such an estimate with the local log-Sobolev inequality of \cite{Yau} to provide us a quantitative version of the following statement -- locally, statistics relax to stationary measure statistics very quickly.
\item With the local equilibrium input of the previous bullet point, we replace $\mathfrak{q}$ with its invariant measure expectation on length-scale $N^{20\beta_{\star}}$. This will be done through localizing to invariant measure calculations and then via energy solution theory ideas \cite{GJ15} at a local scale. This is, again, a precise version of the heuristic that locally, statistics relax to invariant measure expectations. Here, the choice $10\beta_{\star}$ is a sufficiently small exponent so that statistics to relax to invariant measure sufficiently fast, but it is also large enough to gain the benefit of larger-scale time-averages explained in the paragraph below.
\item The upshot of the length-scale $N^{20\beta_{\star}}$ is as follows. The previous bullet point replaces $\mathfrak{q}$ with a function that is controlled by the average density of $\eta$-variables on a block of length-scale $N^{20\beta_{\star}}$. If we believe the $\eta$-variables are sufficiently fluctuating, then this average of $\eta$-variables should be order $N^{-10\beta_{\star}}$ by central limit theorem considerations, for example. This, when combined with the heat kernel estimate for $\mathbf{H}^{\beta_{\star}}$ from Proposition \ref{prop:hke3}, will beat out the factor of $c_{N}\lesssim N$ that is present on the LHS of \eqref{eq:lpe3I} and give us the fixed negative power of $N$ on the RHS of \eqref{eq:lpe3I}. To show this fluctuating property of $\eta$-variables, which is entirely elementary if we are given a product invariant measure, we will interpret averages of $\eta$-variables in terms of properly rescaled gradients of the height function $\mathbf{h}^{N}$; by calculus, we will be able to control this in terms of gradients of its exponential, namely the Gartner transform $\mathbf{Z}^{N}$, at which point it suffices to establish regularity estimates for $\mathbf{Z}^{N}$ present in both \cite{BG,DT}. This scheme was employed in \cite{DT} at the mesoscopic length-scale of roughly $N^{1/2}$. One observation here is that it works at potentially much smaller length-scales of $N^{20\beta_{\star}}$ without much loss. Actually, we will not get the aforementioned $N^{-10\beta_{\star}}$ estimate; we will pick up an arbitrarily small but fixed and positive power of $N$ as well, but this is harmless, as it only forces us to adapt $\gamma$ on the RHS of \eqref{eq:lpe3I} by a small positive factor, for instance.
\end{itemize}
%%%
%%%
\begin{remark}
\fsp The previous outline suggests if we can replace $\mathfrak{q}$ by a spatial-average-type term on a longer length-scale for spatial averaging, then our $\eta$-fluctuation estimates improve, and this ultimately allows us to pick a bigger value of $\bar{\beta}_{\star}$ in Theorem \ref{theorem:kpz}. One way to progressively upgrade the length-scale for spatial averaging, after the initial replacement-by-spatial-average, is by using the two-blocks strategy of \cite{GJ15,GPV} but in a spatially local sense, similar to the one-block step in the outline above. A spatially localized two-blocks scheme would additionally contribute to the local hydrodynamics method herein. We do not include any of the details for this lengthy and involved procedure, because it would not push $\bar{\beta}_{\star}$ in Theorem \ref{theorem:kpz} much closer to $1/2$, but it may be of interest.
\end{remark}
%%%
Let us now make the previous outline precise. This starts with the following decomposition.
%%%
\begin{definition}\label{definition:lp1}
\fsp First, we define space-time averaging operators. Provided any function $\phi:\R_{\geq0}\times\Z\to\R$ and length/time-scales $\mathfrak{l}\in\Z_{\geq0}$ and $\t\geq0$, respectively, we define the following two space-time averaging operators:
\begin{align}
\mathsf{I}^{\mathbf{X},\mathfrak{l}}(\phi_{\s,\y}) \ \overset{\bullet}= \ \wt{\sum}_{\w=1,\ldots,\mathfrak{l}}\phi_{\s,\y+\w} \and \mathsf{I}^{\mathbf{T},\t}(\phi_{\s,\y}) \ \overset{\bullet}= \ \t^{-1}\int_{0}^{\t}\phi_{\s+\r,\y}\d\r
\end{align}
with the interpretation that $\mathbf{I}^{\mathbf{X},0}$ and $\mathbf{I}^{\mathbf{T},0}$ are identity operators.
\end{definition}
%%%
%%%
\begin{definition}\label{definition:lp2}
\fsp With notation to be defined after, we write
\begin{align}
\mathbf{H}^{\beta_{\star},\mathbf{T}}_{\t,\x}\left(2^{-1}c_{N}\mathbf{1}_{\y=0}\mathfrak{q}_{\s,\y}^{\mathrm{tot}}\mathbf{Z}_{\s,\y}^{N}\right) \ = \ \mathbf{H}_{\t,\x}^{\beta_{\star},\mathbf{T},1} + \mathbf{H}_{\t,\x}^{\beta_{\star},\mathbf{T},2} + \mathbf{H}_{\t,\x}^{\beta_{\star},\mathbf{T},3}, \label{eq:lp1}
\end{align}
where the first term on the RHS of \eqref{eq:lp1} cuts off the time-integral in the space-time heat operator a precise distance near the short-time singularity of the heat kernel, where the second term on the RHS is given by modifying $\mathfrak{q}^{\mathrm{tot}}$ upon replacing the leading-order quadratic term by introducing an $\eta$-average, and the last term is what is left after this replacement and cutoff. Precisely, we have
\begin{align}
\mathbf{H}_{\t,\x}^{\beta_{\star},\mathbf{T},1} \ \overset{\bullet}= \ \mathbf{H}^{\beta_{\star},\mathbf{T}}_{\t,\x}\left(\mathbf{1}_{\s\geq\t^{\star}}2^{-1}c_{N}\mathbf{1}_{\y=0}\mathfrak{q}_{\s,\y}^{\mathrm{tot}}\mathbf{Z}_{\s,\y}^{N}\right),
\end{align}
where $\t^{\star}=(\t-N^{-5\beta_{\star}-\e})_{+}$ for a fixed but small positive $\e$. Now write $\bar{\mathfrak{q}}_{\s,\y}=\eta_{\s,\y}(\eta_{\s,\y+1}-\mathsf{I}^{\mathbf{X},\mathfrak{l}_{N}}(\eta_{\s,\y+1}))$ for $\mathfrak{l}_{N}=N^{20\beta_{\star}}$ and
\begin{align}
\mathbf{H}_{\t,\x}^{\beta_{\star},\mathbf{T},2} \ \overset{\bullet}= \ \mathbf{H}^{\beta_{\star},\mathbf{T}}_{\t,\x}\left(\mathbf{1}_{\s\leq\t^{\star}}2^{-1}c_{N}\mathbf{1}_{\y=0}\bar{\mathfrak{q}}_{\s,\y}\mathbf{Z}_{\s,\y}^{N}\right). 
\end{align}
We are left with the following last term on the RHS of \eqref{eq:lp1} with the difference $\mathfrak{d}\overset{\bullet}=\mathfrak{q}^{\mathrm{tot}}-\bar{\mathfrak{q}}$ away from the heat kernel singularity:
\begin{align}
\mathbf{H}_{\t,\x}^{\beta_{\star},\mathbf{T},3} \ \overset{\bullet}= \ \mathbf{H}^{\beta_{\star},\mathbf{T}}_{\t,\x}\left(\mathbf{1}_{\s\leq\t^{\star}}2^{-1}c_{N}\mathbf{1}_{\y=0}\mathfrak{d}_{\s,\y}\mathbf{Z}_{\s,\y}^{N}\right).
\end{align}
\end{definition}
%%%
By the decomposition \eqref{eq:lp1} in Definition \ref{definition:lp2} along with the triangle inequality and union bound over three events, it suffices to prove that the estimate in Proposition \ref{prop:lpe3} holds upon replacing the LHS of \eqref{eq:lp1} with each term on the RHS of \eqref{eq:lp1}.
%%%
\begin{lemma}\label{lemma:lp3}
\fsp We have the following deterministic inequality where $\e$ is introduced in \emph{Definition \ref{definition:lp2}}:
\begin{align}
\|\mathbf{H}_{\t,\x}^{\beta_{\star},\mathbf{T},1}\|_{\t_{\mathrm{st}};\Z;\max} \ \lesssim \ N^{-\e/2}\|\mathbf{Z}^{N}\|_{\t_{\mathrm{st}};\Z}. 
\end{align}
\end{lemma}
%%%
%%%
\begin{prop}\label{prop:lp4}
\fsp There exists a positive constant $\gamma$ as in \emph{Proposition \ref{prop:lpe3}} such that
\begin{align}
\mathbf{P}\left(\|\mathbf{H}_{\t,\x}^{\beta_{\star},\mathbf{T},2}\|_{\t_{\mathrm{st}};\Z;\max}\|\mathbf{Z}^{N}\|_{\t_{\mathrm{st}};\Z}^{-1} \geq N^{-\gamma}\right) \ \lesssim \ N^{-\gamma}.
\end{align}
\end{prop}
%%%
%%%
\begin{prop}\label{prop:lp5}
\fsp There exists a positive constant $\gamma$ as in \emph{Proposition \ref{prop:lpe3}} such that
\begin{align}
\mathbf{P}\left(\|\mathbf{H}_{\t,\x}^{\beta_{\star},\mathbf{T},3}\|_{\t_{\mathrm{st}};\Z;\max}(1+\|\mathbf{Z}^{N}\|_{\t_{\mathrm{st}};\Z}^{2})^{-1} \geq N^{-\gamma}\right) \ \lesssim \ N^{-\gamma}.
\end{align}
\end{prop}
%%%
%%%
\begin{proof}[Proof of \emph{Proposition \ref{prop:lpe3}}]
With the appropriate probability, each term on the RHS of \eqref{eq:lp1} is bounded above in the $\|\|_{\t_{\mathrm{st}};\Z;\max}$-norm by $N^{-\gamma}\|\mathbf{Z}^{N}\|_{\t_{\mathrm{st}};\Z} + N^{-\gamma}(1+\|\mathbf{Z}^{N}\|_{\t_{\mathrm{st}};\Z}^{2})$ with $\gamma$ now depending on $\e$, which we will take fixed anyway. Note the intersection of the three events on which we have said estimate for each of the three terms on the RHS of \eqref{eq:lp1} also holds with the appropriate probability, namely with complement of probability at most $N^{-\gamma}$ times a uniformly bounded factor. On this intersection, by the triangle inequality, we deduce the LHS of \eqref{eq:lp1} is controlled by the same quantity. At this point, we complete the proof by noting $\|\mathbf{Z}^{N}\|_{\t_{\mathrm{st}};\Z}\leq1+\|\mathbf{Z}^{N}\|_{\t_{\mathrm{st}};\Z}^{2}$, which is an elementary convexity inequality, so we are done.
\end{proof}
%%%
The rest of this section is dedicated to proofs for Lemma \ref{lemma:lp3}, Proposition \ref{prop:lp4}, and Proposition \ref{prop:lp5}, which would complete the proof of Proposition \ref{prop:lpe3} provided the argument written above. The proof of Lemma \ref{lemma:lp3} will be fairly straightforward, since it is a deterministic argument. The proofs of Proposition \ref{prop:lp4} and Proposition \ref{prop:lp5} will be somewhat more complicated; they will require a set of preliminary ingredients/estimates that we list, use to prove Propositions \ref{prop:lp4} and \ref{prop:lp5}, respectively, and then provide proofs of afterwards in order to not disrupt the flow of the writing. However, in order to make the argument clearer, we will provide brief explanations for why the aforementioned preliminary ingredients are true as we list them. 
%%%
\subsection{Proof of Lemma \ref{lemma:lp3}}
%%%
Recall from Proposition \ref{prop:mshe1} that $c_{N}\lesssim N$ and $\mathfrak{q}^{\mathrm{tot}}$ is uniformly bounded. Therefore, by definition of $\mathbf{H}^{\beta_{\star},\mathbf{T},1}$ in Definition \ref{definition:lp2}, we may control everything inside of the space-time heat operator in $\mathbf{H}^{\beta_{\star}}$ uniformly in the space-time integration variables and be left with integrating the heat kernel $\mathbf{H}^{\beta_{\star}}$. For this, we emphasize that what is inside the space-time heat operator defining $\mathbf{H}^{\beta_{\star},\mathbf{T},1}$ is concentrated at spatial variable equal to the origin. Therefore, we have
\begin{align}
|\mathbf{H}_{\t,\x}^{\beta_{\star},\mathbf{T},1}| \ \lesssim \ N\|\mathbf{Z}^{N}\mathbf{1}_{\{0\}\subseteq\Z}\|_{\t;\Z;\max}\int_{\t^{\star}}^{\t}\mathbf{H}_{\s,\t,\x,0}^{\beta_{\star}}\d\s \ \lesssim \ N^{\frac52\beta_{\star}}\|\mathbf{Z}^{N}\mathbf{1}_{\{0\}\subseteq\Z}\|_{\t;\Z;\max}\int_{\t^{\star}}^{\t}|\t-\s|^{-1/2}\d\s, \label{eq:lp31}
\end{align}
where the last estimate in \eqref{eq:lp31} follows from the heat kernel estimate in Proposition \ref{prop:hke3}; to be totally transparent, the product of $\mathbf{Z}^{N}$ with the indicator function in \eqref{eq:lp31} refers to taking only the values of $\mathbf{Z}^{N}$ at space-time coordinates of the form $(\s,0)\subseteq[0,1]\times\Z$. We now observe that for functions supported at the origin, the $\|\|_{\t_{\mathrm{st}};\Z;\max}$ and $\|\|_{\t_{\mathrm{st}};\Z}$ norms are equivalent. Therefore, for $\t\leq\t_{\mathrm{st}}$, the norm on the far RHS of \eqref{eq:lp31} is controlled by $\|\mathbf{Z}^{N}\|_{\t_{\mathrm{st}};\Z}$. We now conclude the proof by integrating the integral on the far RHS of \eqref{eq:lp31} and recalling $\t^{\star}$, by construction in Definition \ref{definition:lp2}, satisfies $|\t-\t^{\star}|^{1/2}\leq N^{-5\beta_{\star}/2-\e/2}$. \qed
%%%
\subsection{Proof of Proposition \ref{prop:lp4}}
%%%
As discussed at the beginning of the section, the proof of Proposition \ref{prop:lp4} will be based on reducing our estimates of $\mathbf{H}^{\beta_{\star},\mathbf{T},2}$ to estimates at the stationary measure. Therefore, we will require a set of estimates concerning reduction to the stationary measure and estimates of statistics at the stationary measure. For the former, we introduce an entropy production inequality and reduction to equilibrium, all packaged in the same result. First, we introduce notation for these stationary measures, which are given by product ``grand-canonical" measures and ergodic ``canonical" measures, as well as important and useful functionals and maps on probability densities with respect to these stationary measures, some of which we introduce now because it is convenient but will not use until later in proofs of lemmas in these sections.
%%%
\begin{definition}\label{definition:lp6}
\fsp Provided any subset $\mathbb{I}\subseteq\Z$ and $\sigma\in\R$, we define $\mu_{\sigma,\mathbb{I}}$ as a product measure on $\Omega_{\mathbb{I}}$ whose one-dimensional marginals have expectation $\sigma$. If $\mathbb{I}$ is finite, we also define $\mu_{\sigma,\mathbb{I}}^{\mathrm{can}}$ as the uniform measure on the hyperplane in $\Omega_{\mathbb{I}}$ consisting of all configurations whose average of $\eta$-values in $\mathbb{I}$ is equal to $\sigma$. Next, provided any initial measure $\mu(0)$ on $\Omega_{\Z}$, we will define $\mu(\t)$ as the probability measure on $\Omega_{\Z}$ given by the law of the particle system after time $\t\geq0$ with initial measure $\mu(0)$. We additionally let $\mathsf{R}(\t)$ denote the Radon-Nikodym derivative of $\mu(\t)$ with respect to $\mu_{0,\Z}$, and for any $\t\geq0$, we define the time-averaged law 
\begin{align}
\bar{\mathsf{R}}(\t)=\mathbf{1}_{\t>0}\t^{-1}\int_{0}^{\t}\mathsf{R}(\s)\d\s + \mathbf{1}_{\t=0}\mathsf{R}(0).
\end{align}
Now, let us define important functionals based on the aforementioned invariant measures. These will be the relative entropy and Dirichlet form. When finite, define the following grand-canonical functionals for any probability density $\mathsf{R}$ with respect to $\mu_{0,\Z}$, in which we recall that $\mathscr{L}_{\x}$ is the speed-1 generator of symmetric simple exclusion on the bond $\{\x,\x+1\}$:
\begin{align}
\mathfrak{D}_{\mathrm{KL}}(\mathsf{R}) \ \overset{\bullet}= \ \E^{\mu_{0,\Z}}\mathsf{R}\log\mathsf{R} \and \mathfrak{D}_{\mathrm{Dir}}(\mathsf{R}) \ \overset{\bullet}= \ -{\sum}_{\x\in\Z}\E^{\mu_{0,\Z}}\mathsf{R}^{1/2}\mathscr{L}_{\x}(\mathsf{R}^{1/2}).
\end{align}
Now, we localize and define the following local canonical functionals. First, for any subset $\mathbb{I}\subseteq\Z$ we let $\Pi_{\mathbb{I}}$ denote the pushforward on probability measures/densities induced by the projection map $\Omega_{\Z}\to\Omega_{\mathbb{I}}$ that forgets all $\eta$-variables outside $\mathbb{I}$. For any finite set $\mathbb{I}\subseteq\Z$ and probability density $\mathsf{R}$ with respect to $\mu_{0,\Z}$ and $\sigma\in\R$, we will additionally let $\Pi_{\mathbb{I},\sigma}$ denote first projecting with $\Pi_{\mathbb{I}}$ and then conditioning on the hyperplane support of {$\mu_{\sigma,\mathbb{I}}^{\mathrm{can}}$}, ultimately giving a probability density with respect to {$\mu_{\sigma,\mathbb{I}}^{\mathrm{can}}$} if acting on a probability density with respect to $\mu_{0,\Z}$. Finally, we will define the following relative entropy and Dirichlet form functionals for canonical measures, given any $\sigma\in\R$ and finite $\mathbb{I}\subseteq\Z$ and probability density $\mathsf{R}$ with respect to $\mu_{0,\Z}$:
\begin{align}
\mathfrak{D}_{\mathrm{KL}}^{\sigma,\mathbb{I}}(\mathsf{R}) \ \overset{\bullet}= \ \E^{\sigma,\mathbb{I}}(\Pi_{\mathbb{I},\sigma}\mathsf{R})\log(\Pi_{\mathbb{I},\sigma}\mathsf{R}) \and \mathfrak{D}_{\mathrm{Dir}}^{\sigma,\mathbb{I}}(\mathsf{R}) \ \overset{\bullet}= \ -{\sum}_{\x,\y\in\mathbb{I}}\E^{\sigma,\mathbb{I}}(\Pi_{\mathbb{I},\sigma}\mathsf{R}^{1/2})\mathscr{L}_{\x,\y}(\Pi_{\mathbb{I},\sigma}\mathsf{R}^{1/2}),
\end{align}
where we recall $\mathscr{L}_{\x,\y}$ is the speed-1 symmetric simple exclusion process over the bond $\{\x,\y\}$ and $\E^{\sigma,\mathbb{I}}$ denotes expectation with respect to $\mu_{\sigma,\mathbb{I}}^{\mathrm{can}}$; this last piece of notation is something we will adopt for the rest of this section.
\end{definition}
%%%
%%%
\begin{lemma}\label{lemma:lp7}
\fsp Provided any strictly positive $\delta$ and probability density $\mathsf{Q}(0)$ with respect to $\mu_{0,\Z}$, there exists a probability density $\mathsf{R}(0)$ with respect to $\mu_{0,\Z}$ such that for all times $\t\in[0,1]$, we have the following for $\mathbb{I}_{N}$ the ball of radius $N$ in $\Z$ around the origin:
\begin{align}
\mathfrak{D}_{\mathrm{Dir}}(\bar{\mathsf{R}}(\t)) \ \leq \ \t^{-1}\int_{0}^{\t}\mathfrak{D}_{\mathrm{Dir}}(\mathsf{R}(\s))\d\s \ \lesssim \ N^{-1/2+\beta_{\star}+\delta}\t^{-1} \and \E^{\mu_{0,\Z}}|\Pi_{\mathbb{I}_{N}}\mathsf{R}(\t)-\Pi_{\mathbb{I}_{N}}\mathsf{Q}(\t)| \ \lesssim_{\delta} \ N^{-100}. \label{eq:lp7a}
\end{align}
Moreover, provided any uniformly bounded functional $\mathfrak{b}:\Omega\to\R$ with support contained in a finite subset $\mathbb{I}\subseteq\mathbb{I}_{N}\subseteq\Z$, we have
\begin{align}
\E^{\mu_{0,\Z}}\bar{\mathsf{Q}}(1)\mathfrak{b} \ \lesssim \ N^{-1/2+\beta_{\star}+\delta}|\mathbb{I}|^{2} + {\sup}_{\sigma\in\R}\E^{\sigma,\mathbb{I}}|\mathfrak{b}|. \label{eq:lp7b}
\end{align}
\end{lemma}
%%%
Roughly speaking, the estimates in \eqref{eq:lp7a} compares any probability density $\mathsf{Q}(0)$ to another probability measure that differs, in terms of statistics on a ball with radius $N$, in total variation distance by a negligibly small amount in the large-$N$ limit; for this auxiliary/additional probability measure, we then have a time-averaged Dirichlet form estimate. This comparison is basically done through the observation that because the particles in the system perform random walks with symmetric speed $N^{2}$ and asymmetric speed $N^{3/2}$, data beyond the length-scale of $N^{3/2}$ should be irrelevant for statistics in a ball of length of order $N$, so that we can compare a general density $\mathsf{Q}(0)$ with one that is stationary outside a ball of radius $N^{3/2+\delta}$, for example. Entropy production/Dirichlet form estimates for this auxiliary measure are then standard as in Appendix 1.9 of \cite{KL}; the Dirichlet form is controlled by its relative entropy of order $N^{3/2+\delta}$ times $N^{-2}$ to account for the fast $N^{2}$-speed to equilibrium, and then times $N^{\beta_{\star}}$ to account for the slow bond effect. The estimate \eqref{eq:lp7b} then employs the Dirichlet form estimate, a diffusive-in-length-scale log-Sobolev inequality from \cite{Yau}, and the classical entropy inequality in Appendix 1.8 of \cite{KL} to reduce time-averaged non-stationary expectations to invariant measure expectations up to a manageable cost. We clarify these details later in this section.

Let us now make Lemma \ref{lemma:lp7} useful and reduce our estimate of $\mathbf{H}^{\beta_{\star},\mathbf{T},2}$ to something of the form of the LHS of \eqref{eq:lp7b}. For this, we require a number of preliminary steps. First, we will replace $\bar{\mathfrak{q}}$ in $\mathbf{H}^{\beta_{\star},\mathbf{T},2}$ from Definition \ref{definition:lp2} by its time-average on a time-scale $\t^{\mathrm{av}}$ that is sufficiently above the microscopic time-scale in order to exploit the fluctuating behaviors of $\bar{\mathfrak{q}}$; this fluctuating behavior comes from the fact that $\eta$-values communicate and equilibrate locally in short times, and it can be made precise via the feature that $\E^{\sigma,\mathbb{I}}\bar{\mathfrak{q}}=0$ for \emph{any} $\sigma\in\R$ and $\mathbb{I}$ finite and containing the support of $\bar{\mathfrak{q}}$. The cost in replacement-by-time-average of $\bar{\mathfrak{q}}$ for time-scale $\t^{\mathrm{av}}$ is controlled by the time-regularity of the $\mathbf{H}^{\beta_{\star}}$ heat operator and the time-regularity of $\mathbf{Z}^{N}$, because we are integrating $\bar{\mathfrak{q}}$ times the heat kernel and $\mathbf{Z}^{N}$ in $\mathbf{H}^{\beta_{\star},\mathbf{T},2}$; this final point is basically a version of integration-by-parts to move time-gradients from one factor onto another factor. Time-regularity of the heat kernel is proven in Proposition \ref{prop:hke8}, and that of $\mathbf{Z}^{N}$ is proven in Lemma \ref{lemma:ste}. Roughly speaking, the order-$N$ factor in $\mathbf{H}^{\beta_{\star},\mathbf{T},2}$ is basically cancelled by the $N^{-1}$-factor in the heat kernel estimate of Proposition \ref{prop:hke3}. Time regularity on short time-scales then provides quantitative power-saving in $N$.
%%%
\begin{lemma}\label{lemma:lp8}
\fsp Define $\t^{\mathrm{av}}=N^{-2+100\beta_{\star}}$, where $100$ is just some big constant that we can take bigger if we are willing to take $\beta_{\star}$ smaller but still positive in this paper. For a positive constant $\gamma$ depending only on $\beta_{\star}$ as in \emph{Proposition \ref{prop:lpe3}}, we have the following upper bound outside an event of probability at most $N^{-100}$ times a uniformly bounded constant:
\begin{align}
\|\mathbf{H}^{\beta_{\star},\mathbf{T},2}\|_{\t_{\mathrm{st}};\Z;\max} \ \leq \ \|\mathbf{H}^{\beta_{\star}}_{\t,\x}(\mathbf{1}_{\s\leq\t^{\star}}2^{-1}c_{N}\mathbf{1}_{\y=0}\mathsf{I}^{\mathbf{T},\t^{\mathrm{av}}}(\bar{\mathfrak{q}}_{\s,\y})\mathbf{Z}_{\s,\y}^{N})\|_{\t_{\mathrm{st}};\Z;\max} + N^{-\gamma}\|\mathbf{Z}^{N}\|_{\t_{\mathrm{st}};\Z}. \label{eq:lp8a}
\end{align}
\end{lemma}
%%%
We will now control the first term on the RHS of \eqref{eq:lp8a}. To this end, we will forget about the {$\mathbf{Z}^{N}$ factor therein} at the cost of $\|\mathbf{Z}^{N}\|_{\t_{\mathrm{st}};\Z}$ similar to the proof of Lemma \ref{lemma:lp3}. We then control integrating against the heat kernel with forward spatial-variable $\y=0$ using the on-diagonal heat kernel estimate in Proposition \ref{prop:hke3}; we emphasize that the cutoff $\s\leq\t^{\star}$ allows us to control the heat kernel short-time singularity in Proposition \ref{prop:hke3} and therefore integrate just the $\bar{\mathfrak{q}}$-time-average. This integral of a non-negative absolute value can then be extended from $[0,\t]$ to $[0,1]$ given any $\t\in[0,1]$, giving Lemma \ref{lemma:lp9} below. Ultimately, we obtain the following estimate, which wraps up the calculations discussed in this paragraph by taking an expectation.
\begin{lemma}\label{lemma:lp9}
\fsp Recalling $\e>0$ arbitrarily small from \emph{Definition \ref{definition:lp2}} and $\t^{\mathrm{av}}=N^{-2+100\beta_{\star}}$ from \emph{Lemma \ref{lemma:lp8}}, we have the estimate
\begin{align}
\E\left(\|\mathbf{H}^{\beta_{\star}}_{\t,\x}(\mathbf{1}_{\s\leq\t^{\star}}2^{-1}c_{N}\mathbf{1}_{\y=0}\mathsf{I}^{\mathbf{T},\t^{\mathrm{av}}}(\bar{\mathfrak{q}}_{\s,\y})\mathbf{Z}_{\s,\y}^{N})\|_{\t_{\mathrm{st}};\Z;\max}\|\mathbf{Z}^{N}\|_{\t_{\mathrm{st}};\Z}^{-1}\right) \ \lesssim \ N^{5\beta_{\star}+\e/2}\E\int_{0}^{1}|\mathsf{I}^{\mathbf{T},\mathfrak{t}^{\mathrm{av}}}(\bar{\mathfrak{q}}_{\s,0})|\d\s. \label{eq:lp9a}
\end{align}
\end{lemma}
%%%
The third and final step towards turning our estimate of $\mathbf{H}^{\beta_{\star},\mathbf{T},2}$ into something of the form of the LHS of \eqref{eq:lp7b} is to unfold the RHS of \eqref{eq:lp9a} in the following fashion. Moving the expectation past the integral on the RHS of \eqref{eq:lp9a}, we may decompose the expectation of the time-average of $\bar{\mathfrak{q}}_{\s,0}$ in terms of an expectation with respect to the path-space measure on the Skorokhod path space $\mathscr{D}(\R_{\geq0},\Omega)$ induced by the particle system \emph{conditioning} on the initial configuration in said path-space expectation being the time-$\s$ configuration, afterwards taking an expectation over the time-$\s$ configuration according to the law of the particle system at time $\s$. In particular, the expectation/integral from the RHS of \eqref{eq:lp9a} can be viewed as a path-space expectation of a path-space functional that depends on randomness through its initial configuration that is sampled according to the distribution or law of the time-$\s$ configuration. After integrating $\s$ over $[0,1]$, we are then sampling the initial configuration for this path-space expectation according to the time-averaged law of the particle system, thus realizing the RHS of \eqref{eq:lp9a} as the LHS of \eqref{eq:lp7b} with $\mathfrak{b}$ equal to the path-space expectation of the time-average of $\bar{\mathfrak{q}}_{0,0}$ with initial condition sampled according to the time-averaged law of the particle system. The last piece of this final step/following result is then to address the following issue. Said choice of $\mathfrak{b}$, which is a functional of $\Omega_{\Z}$ through the initial configuration of its path-space expectation, has infinite support. Indeed, although $\bar{\mathfrak{q}}$ has finite support in $\Z$, for any positive time, the position of particles anywhere in $\Z$ have, in principle, an effect on the $\eta$-variables that $\bar{\mathfrak{q}}$ depends upon. This is an analog of the fact that a Brownian motion has an infinite propagation speed. However, particles that are sufficiently far from the support of $\bar{\mathfrak{q}}$, in a sense depending on the time $\t^{\mathrm{av}}$ that we consider the path-space for, should effect the time-average of $\bar{\mathfrak{q}}$ over a time-scale of $\t^{\mathrm{av}}$ with low probability; see the paragraph after Lemma \ref{lemma:lp7}. Making this precise ultimately gives:
%%%
\begin{lemma}\label{lemma:lp10}
\fsp Recalling {$\t^{\mathrm{av}}=N^{-2+100\beta_{\star}}$} from \emph{Lemma \ref{lemma:lp8}}, we have the following with notation to be defined afterwards:
\begin{align}
\E\int_{0}^{1}|\mathsf{I}^{\mathbf{T},{\t^{\mathrm{av}}}}(\bar{\mathfrak{q}}_{\s,0})|\d\s \ \lesssim \ \E^{\mu_{0,\Z}}\bar{\mathsf{Q}}(1)\E^{\mathrm{path},\t^{\mathrm{av}}}|\mathsf{I}^{\mathbf{T},{\t^{\mathrm{av}}}}(\bar{\mathfrak{q}}_{0,0})| + N^{-100}. \label{eq:lp10a}
\end{align}
Above, $\bar{\mathsf{Q}}(1)$ is the time-averaged law of the particle system over $[0,1]$ with the initial measure $\mathsf{Q}(0)$. We now define the $\E^{\mathrm{path},\t^{\mathrm{av}}}$ expectation. Let $\mathbb{I}(\t^{\mathrm{av}})$ be the ball in $\Z$ around the origin of radius $\iota(\t^{\mathrm{av}})=N(\t^{\mathrm{av}})^{1/2}\log^{100}N + N^{3/2}\t^{\mathrm{av}}\log^{100}N+N^{20\beta_{\star}}$.  Next, we let $\E^{\mathrm{path},\t^{\mathrm{av}}}$ be expectation with respect to the path-space measure on $\mathscr{D}(\R_{\geq0},\Omega_{\mathbb{I}(\t^{\mathrm{av}})})$ induced by a ``periodization" of the particle system on $\mathbb{I}(\t^{\mathrm{av}})$, which is the Markov process with state space $\Omega_{\mathbb{I}(\t^{\mathrm{av}})}$ defined by running the particle system but when a particle tries to jump from inside $\mathbb{I}(\t^{\mathrm{av}})$ to outside $\mathbb{I}(\t^{\mathrm{av}})$, we stipulate that it instead jumps to where it would if $\mathbb{I}(\t^{\mathrm{av}})$ were realized as a discrete torus with periodic boundary conditions $1+\sup\mathbb{I}(\t^{\mathrm{av}})\overset{\bullet}=\inf\mathbb{I}(\t^{\mathrm{av}})$. In particular, the time-average $\mathsf{I}^{\mathbf{T},\t^{\mathrm{av}}}(\bar{\mathfrak{q}}_{0,0})$ is the time-average of $\bar{\mathfrak{q}}_{0,0}$ evaluated at particle configurations given by taking an initial configuration on $\mathbb{I}(\t^{\mathrm{av}})$ sampled in the expectation on the RHS of \eqref{eq:lp10a} and evolving it under the aforementioned ``periodization" dynamic. We note that these particle configurations at which we evaluate $\bar{\mathfrak{q}}_{0,0}$ live in $\Omega_{\mathbb{I}(\t^{\mathrm{av}})}$, which contains the support of $\bar{\mathfrak{q}}_{0,0}$ because the support of $\bar{\mathfrak{q}}_{0,0}$, which is constructed in \emph{Definition \ref{definition:lp2}}, contains the origin and has length $N^{20\beta_{\star}}$ and is thus contained in $\mathbb{I}(\t^{\mathrm{av}})$.
\end{lemma}
%%%
The next step is we take to estimate the first term on the RHS of \eqref{eq:lp10a} via \eqref{eq:lp7b} in Lemma \ref{lemma:lp7} with $\mathfrak{b}$ equal to the $\E^{\mathrm{path},\t^{\mathrm{av}}}$ term in \eqref{eq:lp10a}; note this choice of $\mathfrak{b}$ only depends on particle configurations on $\mathbb{I}(\t^{\mathrm{av}})$ and thus $\mathbb{I}=\mathbb{I}(\t^{\mathrm{av}})$ in said application of \eqref{eq:lp7b} to \eqref{eq:lp10a}. We are then left to estimate canonical measure expectations of $\E^{\mathrm{path},\t^{\mathrm{av}}}$, which is the point of the next result. In short, the next result follows via the Kipnis-Varadhan inequality for stationary processes, which provides Brownian estimates for integrals of functionals that satisfy the same fluctuating property as $\bar{\mathfrak{q}}$. This Kipnis-Varadhan inequality is a key estimate in \cite{GJ15} for stationary models and can be found in \cite{GJ15,KL}; we explain this inequality and the proof of the following in detail at the end of this section.
%%%
\begin{lemma}\label{lemma:lp11}
\fsp Recalling $\t^{\mathrm{av}}=N^{-2+100\beta_{\star}}$ from \emph{Lemma \ref{lemma:lp8}} and retaining notation of \emph{Lemma \ref{lemma:lp10}}, we have
\begin{align}
{\sup}_{\sigma\in\R}\E^{\sigma,\mathbb{I}}\E^{\mathrm{path},\t^{\mathrm{av}}}|\mathsf{I}^{\mathbf{T},{\t^{\mathrm{av}}}}(\bar{\mathfrak{q}}_{0,0})| \ \lesssim \ N^{-1}(\t^{\mathrm{av}})^{-1/2}N^{20\beta_{\star}} \ \lesssim \ N^{-30\beta_{\star}}. \label{eq:lp11a}
\end{align}
\end{lemma}
%%%
%%%
\subsubsection{Wrapping up the proof of \emph{Proposition \ref{prop:lp4}}}
%%%
Let us now summarize the architecture for establishing Proposition \ref{prop:lp4} we have discussed in this subsection so far. First, we combine Lemma \ref{lemma:lp8}, Lemma \ref{lemma:lp9} and Lemma \ref{lemma:lp10} to deduce the following estimate, in which $\gamma$ is a positive constant depending only on $\beta_{\star}\leq1$ as in Proposition \ref{prop:lpe3}:
\begin{align}
\E\left(\|\mathbf{H}^{\beta_{\star},\mathbf{T},2}\|_{\t_{\mathrm{st}};\Z;\max}\|\mathbf{Z}^{N}\|_{\t_{\mathrm{st}};\Z}^{-1}\right) \ \lesssim \ N^{5\beta_{\star}+\e/2}\E^{\mu_{0,\Z}}\bar{\mathsf{Q}}(1)\E^{\mathrm{path},\t^{\mathrm{av}}}|\mathsf{I}^{\mathbf{T},\mathfrak{t}^{\mathrm{av}}}(\bar{\mathfrak{q}}_{0,0})| + N^{-\gamma}. \label{eq:lp41}
\end{align}
To estimate the first term on the RHS of \eqref{eq:lp41}, we employ \eqref{eq:lp7b} with $\mathfrak{b}$ equal to the $\E^{\mathrm{path},\t^{\mathrm{av}}}$ term on the RHS of \eqref{eq:lp41}, whose support is $\mathbb{I}=\mathbb{I}(\t^{\mathrm{av}})$ from Lemma \ref{lemma:lp10}. We note that this choice is legal because $\mathbb{I}(\t^{\mathrm{av}})\subseteq\Z$ is finite and the $\E^{\mathrm{path},\t^{\mathrm{av}}}$ term on the RHS of \eqref{eq:lp41} is uniformly bounded since it averages the uniformly bounded functional $\bar{\mathfrak{q}}$. Thus, we have, for $\delta=\e/2$ in \eqref{eq:lp7b},
\begin{align}
N^{5\beta_{\star}+\e/2}\E^{\mu_{0,\Z}}\bar{\mathsf{Q}}(1)\E^{\mathrm{path},\t^{\mathrm{av}}}|\mathsf{I}^{\mathbf{T},\mathfrak{t}^{\mathrm{av}}}(\bar{\mathfrak{q}}_{0,0})| \ \lesssim \ N^{-1/2+6\beta_{\star}+\e}|\mathbb{I}(\t^{\mathrm{av}})|^{2} + N^{-25\beta_{\star}+\e/2}, \label{eq:lp42}
\end{align}
where the supremum over $\sigma\in\R$ in \eqref{eq:lp7b} is estimated for our choice of $\mathfrak{b}$ via Lemma \ref{lemma:lp11}. Recalling $|\mathbb{I}(\t^{\mathrm{av}})|\lesssim\iota(\t^{\mathrm{av}})$ from the statement of Lemma \ref{lemma:lp10}, plugging this into the first term in \eqref{eq:lp42}, and then elementary power-counting shows the first term on the RHS of \eqref{eq:lp42} is controlled by $N^{-\gamma}$ for $\gamma$ positive and depending only on $\beta_{\star}$, and thus the same for the RHS of \eqref{eq:lp41}. Indeed, note $\iota(\t^{\mathrm{av}})$ is controlled by powers of $N^{\beta_{\star}}$ with uniformly bounded exponents, and therefore if $\beta_{\star}$ is sufficiently small though still positive and fixed, these powers cannot beat the $N^{-1/2}$ factor on the RHS of \eqref{eq:lp42}. To deduce Proposition \ref{prop:lp4} from estimating the expectation on the LHS of \eqref{eq:lp41} by $N^{-\gamma}$ times a uniformly bounded factor, it suffices to employ the Markov inequality. \qed
%%%
\subsection{Proof of Proposition \ref{prop:lp5}}
%%%
The proof of Proposition \ref{prop:lp5} will require a significantly fewer number of preliminary steps as the proof of Proposition \ref{prop:lp4}. First, we compute the difference $\mathfrak{d}=\mathfrak{q}^{\mathrm{tot}}-\bar{\mathfrak{q}}$ we defined in Definition \ref{definition:lp2} that is relevant in $\mathbf{H}^{\beta_{\star},\mathbf{T},3}$. Inspecting the definition of $\bar{\mathfrak{q}}$ in Definition \ref{definition:lp2} and the definition of $\mathfrak{q}^{\mathrm{tot}}$ in Proposition \ref{prop:mshe1}, we observe the following estimate:
\begin{align}
|\mathbf{H}^{\beta_{\star},\mathbf{T},3}_{\t,\x}| \ \lesssim \ N|\mathbf{H}^{\beta_{\star}}_{\t,\x}\left(\mathbf{1}_{\s\leq\t^{\star}}\mathbf{1}_{\y=0}|\mathsf{I}^{\mathbf{X},\mathfrak{l}_{N}}(\eta_{\s,\y+1})|\mathbf{Z}^{N}_{\s,\y}\right)| + N^{1/2}|\mathbf{H}^{\beta_{\star}}_{\t,\x}\left(\mathbf{1}_{\y=0}\mathbf{Z}_{\s,\y}^{N}\right)|. \label{eq:lp51}
\end{align}
We estimate both terms on the RHS of \eqref{eq:lp51}. We start with the second term on the RHS of \eqref{eq:lp51} because it is a straightforward deterministic estimate based on heat kernel estimates of Proposition \ref{prop:hke3}, similar to the proof of Lemma \ref{lemma:lp3}. In the following result, we also record a deterministic estimate, which is proved via very similar elementary ideas, for the first term on the RHS of \eqref{eq:lp51} that will set up our forthcoming analysis for said first term.
%%%
\begin{lemma}\label{lemma:lp12}
\fsp We have the following deterministic estimate:
\begin{align}
 N^{1/2}\|\mathbf{H}^{\beta_{\star}}_{\t,\x}\left(\mathbf{1}_{\y=0}\mathbf{Z}_{\s,\y}^{N}\right)\|_{\t_{\mathrm{st}};\Z;\max} \ \lesssim \ N^{-1/2+5\beta_{\star}/2}\|\mathbf{Z}^{N}\|_{\t_{\mathrm{st}};\Z}. \label{eq:lp12a}
\end{align}
Additionally, we recall $\t^{\star}=\t-N^{-5\beta_{\star}-\e}$ as a function of time $\t$ and our small $\e$ in \emph{Definition \ref{definition:lp2}} that we will eventually choose sufficiently small but still strictly positive and depending only on $\beta_{\star}$. We have the following deterministic estimate:
\begin{align}
 N\|\mathbf{H}^{\beta_{\star}}_{\t,\x}\left(\mathbf{1}_{\s\leq\t^{\star}}\mathbf{1}_{\y=0}|\mathsf{I}^{\mathbf{X},\mathfrak{l}_{N}}(\eta_{\s,\y+1})|\mathbf{Z}^{N}_{\s,\y}\right)\|_{\t_{\mathrm{st}};\Z;\max} \ \lesssim \ N^{5\beta_{\star}+\e}\int_{0}^{1}|\mathsf{I}^{\mathbf{X},\mathfrak{l}_{N}}(\eta_{\s,1})|\mathbf{Z}^{N}_{\s,0}\d\s. \label{eq:lp12b}
\end{align}
\end{lemma}
%%%
Let us now direct the rest of our attention to the term on the RHS of \eqref{eq:lp12b}. Our analysis of this term will require two steps. The first step translates the spatial average $\mathsf{I}^{\mathbf{X},\mathfrak{l}_{N}}$ in terms of regularity of the microscopic Cole-Hopf transform; this is inspired by \cite{DT}. Roughly speaking, the height function integrates $\eta$-variables, so the spatial average of length-scale $\mathfrak{l}_{N}$ of $\eta$-variables can be written as a length-scale $\mathfrak{l}_{N}$ gradient of the height function. However, partly inspiring the approach via Gartner transform $\mathbf{G}^{N}$ of this paper towards universality of the KPZ equation, we do not have sufficiently precise control over analytic properties of the height function. Thus, we will translate these analytic regularity estimates for $\mathbf{h}^{N}$ into those for its exponential $\mathbf{G}^{N}$ via elementary calculus, which then are translated into estimates for $\mathbf{Z}^{N}$ with high probability by Lemma \ref{lemma:mshe5}.
%%%
\begin{lemma}\label{lemma:lp13}
\fsp We have the following outside an event of probability of order $N^{-100}$; the implied constant is independent of $\s\geq0$:
\begin{align}
|\mathsf{I}^{\mathbf{X},\mathfrak{l}_{N}}(\eta_{\s,1})|\mathbf{Z}^{N}_{\s,0} \ \lesssim \ N^{1/2}|\mathfrak{l}_{N}|^{-1}|\grad_{\mathfrak{l}_{N}}^{\mathbf{X}}\mathbf{Z}_{\s,0}^{N}| + N^{-99}.
\end{align}
\end{lemma}
%%%
Thus, we are left to estimate the spatial gradient of the Gartner transform $\mathbf{Z}^{N}$ appearing in Lemma \ref{lemma:lp13}. We will do so in the following pointwise moment topology, which will be sufficient because taking expectations of both sides of \eqref{eq:lp12b} and appealing to Lemma \ref{lemma:lp13}, we are ultimately left with estimating pointwise moments of the gradient of $\mathbf{Z}^{N}$. In short, the proof of the following result amounts to a stochastic analysis like that in the proof of Lemma \ref{lemma:pc3}, for example. Therefore, it is a rather standard procedure modulo technical adjustments to deal with the $\mathbf{Z}^{N}$ space-time norm below. Although we only require a first moment-type bound for the spatial gradient of $\mathbf{Z}^{N}$ in Lemma \ref{lemma:lp13}, we state a second moment-type bound, which certainly implies a first moment-type bound, because it is more natural in terms of the proof of the following result.
%%%
\begin{lemma}\label{lemma:lp14}
\fsp There exists positive $\gamma$ depending only on $\beta_{\star}$ as in \emph{Proposition \ref{prop:lpe3}} such that
\begin{align}
\E\left(\mathbf{1}_{0\leq\t\leq{\t_{\mathrm{st}}}}\t|\grad_{\mathfrak{l}_{N}}^{\mathbf{X}}\mathbf{Z}_{\t,0}^{N}|^{2}(1+\|\mathbf{Z}^{N}\|_{{\t_{\mathrm{st}}};\Z}^{2})^{-2}\right) \ \lesssim \ N^{-1-10\beta_{\star}-\gamma}|\mathfrak{l}_{N}|^{2}.
\end{align}
\end{lemma}
%%%
We briefly remark that Lemma \ref{lemma:lp14} implies that the spatial gradient on the RHS of the estimate in Lemma \ref{lemma:lp13} is controlled by the $1+\|\mathbf{Z}^{N}\|_{\t_{\mathrm{st}};\Z}^{2}$ factor we allow for times an integrable singularity of $\s^{-1/2}$, at least in a high probability of expectation sense.
%%%
\subsubsection{Wrapping up the proof of \emph{Proposition \ref{prop:lp5}}}
%%%
By the estimate \eqref{eq:lp51} and Lemma \ref{lemma:lp12}, namely the estimate \eqref{eq:lp12a} therein, it suffices to estimate the first term on the RHS of \eqref{eq:lp51} with the appropriate high probability. To this end, we estimate its expectation; by Lemma \ref{lemma:lp12} and \eqref{eq:lp12b} therein combined with Lemma \ref{lemma:lp13} and Lemma \ref{lemma:lp14}, we deduce that the expectation of the $\|\|_{\t_{\mathrm{st}};\Z;\max}$-norm of the first term on the RHS of \eqref{eq:lp51} times $(1+\|\mathbf{Z}^{N}\|_{\t_{\mathrm{st}};\Z}^{2})^{-1}$ is at most $N^{-\gamma/2+\e}$ times a bounded constant. The Markov inequality then tells us that with high probability, the first term on the RHS of \eqref{eq:lp51} is controlled by a uniformly bounded constant times $(1+\|\mathbf{Z}^{N}\|_{\t_{\mathrm{st}};\Z}^{2})$ times a negative power of $N$ with an exponent $-\gamma/2+\e$ that is strictly positive and depends only on $\beta_{\star}$, if $\e>0$ is sufficiently small but still universal depending on $\beta_{\star}$. The triangle inequality, \eqref{eq:lp51}, and high-probability bounds for the RHS of \eqref{eq:lp51} finishes the proof as in the proof of Proposition \ref{prop:lp4}. \qed
%%%
\subsection{Proofs of Auxiliary Ingredients}
%%%
The remainder of this section is dedicated to estimates of the preliminary lemmas used to prove Proposition \ref{prop:lp4} and Proposition \ref{prop:lp5}, respectively. We provide proofs for these estimates in the order they were written.
%%%
\begin{proof}[Proof of \emph{Lemma \ref{lemma:lp7}}]
We define $\mathsf{R}(0)$ to be the tensor of $\Pi_{\mathbb{I}^{N}}\mathsf{Q}(0)$ with $\Pi_{\Z\setminus\mathbb{I}^{N}}\mu_{0,\Z}$ in which $\mathbb{I}^{N}\subseteq\Z$ is the ball in $\Z$ centered at the origin with radius $N^{3/2+\delta}$ for our arbitrarily small but fixed and positive $\delta$. In particular, $\mathsf{R}(0)$ cuts off $\mathsf{Q}(0)$ outside $\mathbb{I}^{N}$ and replaces data outside $\mathbb{I}^{N}$ with stationary grand-canonical measure data. Following standard entropy production as in Appendix 1.9 of \cite{KL} while accounting for the slow bond effect globally, which is sub-optimal but sufficient for our purposes, we have
\begin{align}
\partial_{\t}\mathfrak{D}_{\mathrm{KL}}(\mathsf{R}(\t)) \ \lesssim \ -N^{2-\beta_{\star}}\mathfrak{D}_{\mathrm{Dir}}(\mathsf{R}(\t)). \label{eq:lp71}
\end{align}
The advantage behind $\mathsf{R}(\t)$ is that $\mathfrak{D}_{\mathrm{KL}}(\mathsf{R}(\t))\leq\mathfrak{D}_{\mathrm{KL}}(\mathsf{R}(0))\lesssim N^{3/2+\delta}$, where the first inequality follows by standard monotonicity of relative entropy that is implied by \eqref{eq:lp71}, and the second inequality follows from observing $\mathsf{R}(0)$ is stationary outside the set $\mathbb{I}^{N}$ whose length is order $N^{3/2+\delta}$ and whose associated state space $\Omega_{\mathbb{I}^{N}}$ has size exponential in $N^{3/2+\delta}$; the relative entropy bound ultimately follows by taking logarithm of this $\Omega_{\mathbb{I}^{N}}$ size. Integrating \eqref{eq:lp71} in time and using this bound implies
\begin{align}
\int_{0}^{\t}\mathfrak{D}_{\mathrm{dir}}(\mathsf{R}(\s))\d\s \ \lesssim \ N^{-1/2+\beta_{\star}+\delta}. \label{eq:lp72}
\end{align}
Dividing by $\t$ on both sides of \eqref{eq:lp72} and convexity of the Dirichlet form gives the first estimate in \eqref{eq:lp7a}. To establish the second estimate in \eqref{eq:lp7a}, we first observe the LHS of said estimate as the total variation distance between the projection of the time-$\t$ law of $\mathsf{Q}(0)$ and $\mathsf{R}(0)$ initial data onto the set $\mathbb{I}_{N}$. In particular, to control said LHS/total variation distance, it suffices to find a coupling between two copies of the particle system, which we denote by $\mathsf{S}(\mathsf{Q})$ and $\mathsf{S}(\mathsf{R})$, which begin with initial measures given by $\mathsf{Q}(0)$ and $\mathsf{R}(0)$, respectively, and show the copies coincide on $\mathbb{I}_{N}$ with sufficiently high probability. We describe this coupling below.
%%%
\begin{itemize}
\item We first couple the initial data of $\mathsf{S}(\mathsf{Q})$ and $\mathsf{S}(\mathsf{R})$. First, sample a configuration according to $\mathsf{Q}(0)$. Second, take the projection of this configuration onto $\mathbb{I}^{N}$ constructed at the beginning of this proof, and glue it together with a configuration outside $\mathbb{I}^{N}$ sampled independently from the projection of $\mu_{0,\Z}$ outside $\mathbb{I}^{N}$. This gives a configuration whose law is given by $\mathsf{R}(0)$ by construction.
\item We now couple the dynamics/paths of these two species. Observe that the symmetric part of the particle random walks/exclusion process along a bond may be viewed as swapping $\eta$-values at the two points attached to the bond according to a constant speed clock of speed $N^{2}$. We will couple the symmetric parts of $\mathsf{S}(\mathsf{Q})$ and $\mathsf{S}(\mathsf{R})$ dynamics by coupling these ``spin-swap" clocks; the $\eta$-values along any bond in the $\mathsf{S}(\mathsf{Q})$ process will swap exactly when the $\eta$-values along the same bond in the $\mathsf{S}(\mathsf{R})$ swap values. As for the totally asymmetric part, we apply the basic coupling, so particles at the same location in the $\mathsf{S}(\mathsf{Q})$ and $\mathsf{S}(\mathsf{R})$ dynamics will move together according to a clock in the asymmetric part of the particle system dynamics whenever possible.
\end{itemize}
%%%
Observe now that the total variation distance between $\mathsf{R}(\t)$ and $\mathsf{Q}(\t)$ is controlled by the probability that at time $\t$ under the above coupling, the species $\mathsf{S}(\mathsf{Q})$ and $\mathsf{S}(\mathsf{R})$ disagree on $\eta$-values in $\mathbb{I}_{N}$ from the statement of Lemma \ref{lemma:lp7}.  To control this disagreement probability, let us make the following two observations concerning ``discrepancies" between $\mathsf{S}(\mathsf{Q})$ and $\mathsf{S}(\mathsf{R})$, which we define to be points on the lattice $\Z$ where $\mathsf{S}(\mathsf{Q})$ and {$\mathsf{S}(\mathsf{R})$} have disagreeing $\eta$-values.
%%%
\begin{itemize}
\item Initially, all discrepancies are supported outside of $\mathbb{I}^{N}$, the radius $N^{3/2+\delta}$ ball centered at the origin. Decomposing the complement in $\Z$ of $\mathbb{I}^{N}$ into adjacent and connected discrete intervals with length $N^{3/2+\delta}$, the disagreement probability in the previous paragraph is then controlled by the probability that any of the possible $N^{3/2+\delta}$-many discrepancies, which are at least $\mathfrak{l}N^{3/2+\delta}$ away from the radius $N$ neighborhood $\mathbb{I}_{N}$ from the statement of Lemma \ref{lemma:lp7}, travels a necessary distance of $2^{-1}\mathfrak{l}N^{3/2+\delta}$ in order to end up in $\mathbb{I}_{N}$, then summed over all $\mathfrak{l}\geq1$. Equivalently, we are taking a union bound over how far a discrepancy must travel to end up in $\mathbb{I}_{N}$ and over all possible discrepancies that must travel such a distance. Note we crucially use that discrepancies cannot be created under the coupling in the previous list of bullet points, but only either transported/moved or annihilated.
\item The motion of a discrepancy is given by a symmetric simple random walk with speed $N^{2}$ plus a random drift/killing with speed of order $N^{3/2}$; this can be readily checked and follows from our choice of coupling spin-swap clocks/dynamics rather than the basic coupling for the symmetric parts of $\mathsf{S}(\mathsf{Q})$ and $\mathsf{S}(\mathsf{R})$. Therefore, the disagreement probability we are trying to bound is at most the following, in which $\mathbf{P}(\mathfrak{l})$ is the probability that a simple random walk of symmetric speed $N^{2}$ and random asymmetric mechanism of speed $N^{3/2}$ has a maximal displacement of $N^{3/2+\delta}|\mathfrak{l}|$; in what follows, the exponent $p\geq1$ is chosen shortly:
\begin{align}
{\sum}_{\mathfrak{l}\geq1}N^{3/2+\delta}\mathbf{P}({\mathfrak{l}}) \ \lesssim_{p} \ {\sum}_{\mathfrak{l}\geq1}N^{3/2+\delta}\mathfrak{l}^{-2p}N^{-2p\delta} \ \lesssim \ N^{-2p\delta+3/2+\delta}. \label{eq:lp73}
\end{align}
The estimates in \eqref{eq:lp73} follow from controlling $\mathbf{P}({\mathfrak{l}})$ by controlling the asymmetric part of the aforementioned simple random walk by the number of times a Poisson clock rings in time $N^{3/2}$ and controlling the symmetric part by using the Doob maximal inequality for martingales. Provided any $\delta$ fixed and positive, we choose $p\geq1$ sufficiently large to control the RHS of \eqref{eq:lp73}, and therefore the total variation of interest, by $N^{-100}$ times a $\delta$-dependent factor. This completes the proof of \eqref{eq:lp7a}. 
\end{itemize}
%%%
We now move to \eqref{eq:lp7b}. To this end, provided the second estimate in \eqref{eq:lp7a}, we deduce the total variation distance between $\bar{\mathsf{Q}}(1)$ and $\bar{\mathsf{R}}(1)$ is also of order at most $N^{-100}$, because total variation distance is convex. Therefore, because $\mathfrak{b}$ on the LHS of \eqref{eq:lp7b} is uniformly bounded, we may exchange $\bar{\mathsf{Q}}(1)$ on the LHS of \eqref{eq:lp7b} with $\bar{\mathsf{R}}(1)$ if we allow for a cost of order at most $N^{-100}$, which is absorbed into the first term on the RHS of \eqref{eq:lp7b}. Having made this replacement, we note the following, in which $\mathbf{p}(\sigma)$ denotes the probability under $\bar{\mathsf{R}}(1)$ that the average of $\eta$-values on $\mathbb{I}$ is equal to $\sigma$; indeed, the estimate below follows by standard manipulations with conditioning on the $\eta$-average on $\mathbb{I}$, and the sum on the RHS is over a finite set of $\sigma\in\R$ for which $\mathbf{p}(\sigma)\neq0$:
\begin{align}
\E^{\mu_{0,\Z}}\bar{\mathsf{R}}(1)\mathfrak{b} \ = \ {\sum}_{\sigma\in\R}\mathbf{p}(\sigma)\left(\E^{\sigma,\mathbb{I}}\Pi_{\mathbb{I},\sigma}\bar{\mathsf{R}}(1) \cdot \mathfrak{b}\right).  \label{eq:lp74}
\end{align}
We now employ the entropy inequality from Appendix 1.8 of \cite{KL} per $\sigma$, and then sum over $\sigma\in\R$ against $\mathbf{p}(\sigma)$ weights:
\begin{align}
|\E^{\mu_{0,\Z}}\bar{\mathsf{R}}(1)\mathfrak{b}| \ \lesssim \ {\sum}_{\sigma\in\R}\mathbf{p}(\sigma)\mathfrak{D}_{\mathrm{KL}}^{\sigma,\mathbb{I}}(\bar{\mathsf{R}}(1)) \ + \ {\sum}_{\sigma\in\R}\mathbf{p}(\sigma)\log\E^{\sigma,\mathbb{I}}\exp\left(|\mathfrak{b}|\right). \label{eq:lp75}
\end{align}
Because $|\mathfrak{b}|$ is uniformly bounded, by the inequality $\log(1+|\x|)\leq|\x|$ and the uniform Lipschitz property of $\exp$ on compact domains, we deduce that the log-term on the RHS of \eqref{eq:lp75} is bounded by $|\mathfrak{b}|$ up to a uniformly bounded factor. Because $\mathbf{p}(\sigma)$ is a probability measure on $\sigma\in\R$ corresponding to probabilities of disjoint hyperplanes in $\Omega_{\mathbb{I}}$, we may therefore control the second term on the RHS of \eqref{eq:lp75} by a supremum over $\sigma\in\R$ to get an estimate of the following fashion:
\begin{align}
|\E^{\mu_{0,\Z}}\bar{\mathsf{R}}(1)\mathfrak{b}| \ \lesssim \ {\sum}_{\sigma\in\R}\mathbf{p}(\sigma)\mathfrak{D}_{\mathrm{KL}}^{\sigma,\mathbb{I}}(\bar{\mathsf{R}}(1)) \ + \ {\sup}_{\sigma\in\R}\E^{\sigma,\mathbb{I}}|\mathfrak{b}|. \label{eq:lp76}
\end{align}
To obtain \eqref{eq:lp7b}, it therefore remains to estimate the first term on the RHS of \eqref{eq:lp76}. First, because we take relative entropy with respect to a canonical measure, we may employ the log-Sobolev inequality of \cite{Yau}. This is quadratic in the length-scale, because exclusion processes are driven by random walks, and therefore
\begin{align}
{\sum}_{\sigma\in\R}\mathbf{p}(\sigma)\mathfrak{D}_{\mathrm{KL}}^{\sigma,\mathbb{I}}(\bar{\mathsf{R}}(1)) \ \lesssim \ |\mathbb{I}|^{2}{\sum}_{\sigma\in\R}\mathbf{p}(\sigma)\mathfrak{D}_{\mathrm{Dir}}^{\sigma,\mathbb{I}}(\bar{\mathsf{R}}(1)) \ \leq \ |\mathbb{I}|^{2}\mathfrak{D}_{\mathrm{Dir}}(\bar{\mathsf{R}}(1)), \label{eq:lp77}
\end{align}
where the far RHS of \eqref{eq:lp77} is the Dirichlet form of $\bar{\mathsf{R}}(1)$. Indeed, the last statement in \eqref{eq:lp77} follows from standard manipulations of conditional expectation that glues canonical measure Dirichlet forms into a grand-canonical measure Dirichlet form. Actually, because we only look at the Dirichlet form terms in $\mathbb{I}$ in the middle of \eqref{eq:lp77}, the last inequality in \eqref{eq:lp77} is an equality if we replace $\bar{\mathsf{R}}(1)$ by its projection onto $\mathbb{I}\subseteq\Z$, but convexity allows us to remove this projection if we only want an upper bound. At this point, it suffices to employ the first estimate in \eqref{eq:lp7a} to the far RHS of \eqref{eq:lp77} and plug in the resulting estimate for the far LHS of \eqref{eq:lp77} into the first term on the RHS of \eqref{eq:lp76}. This ultimately yields \eqref{eq:lp7b}, so we are done.
\end{proof}
%%%
%%%
\begin{proof}[Proof of \emph{Lemma \ref{lemma:lp8}}]
Let us first take the difference between $\mathbf{H}^{\beta_{\star},\mathbf{T},2}$ and the term inside the norm on the RHS of \eqref{eq:lp8a}. Where $\mathbf{H}^{\beta_{\star},\mathbf{T},2}$ has $\bar{\mathfrak{q}}$, the term inside the norm on the RHS of \eqref{eq:lp8a} has the time-average $\bar{\mathfrak{q}}$ on time-scale $\t^{\mathrm{av}}$ defined in the statement of Lemma \ref{lemma:lp8}. Thus, the difference mentioned in the first sentence of this proof is controlled by the heat operator acting before time $\t^{\star}$ and at $\y=0$ of the Gartner transform $\mathbf{Z}_{\s,0}^{N}$ times the difference between $\bar{\mathfrak{q}}$ and said time-average. Because the difference between $\bar{\mathfrak{q}}$ its scale-$\t^{\mathrm{av}}$ time-average is an average of time-gradients on non-negative time-scales $\r\leq\t^{\mathrm{av}}$, it suffices to control the following on some high probability event for a reason we explain in slightly more detail afterwards:
\begin{align}
{\sup}_{0\leq\r\leq\t^{\mathrm{av}}}\|\mathbf{H}^{\beta_{\star},\mathbf{T}}_{\t,\x}(\mathbf{1}_{\s\leq\t^{\star}}c_{N}\mathbf{1}_{\y=0}\grad_{\r}^{\mathbf{T}}\bar{\mathfrak{q}}_{\s,\y}\mathbf{Z}_{\s,\y}^{N})\|_{\t_{\mathrm{st}};\Z;\max}. \label{eq:lp81}
\end{align}
Indeed, by the triangle inequality, the LHS of \eqref{eq:lp8a} is controlled by the first term on the RHS of \eqref{eq:lp8a} plus the $\|\|_{\t_{\mathrm{st}};\Z;\max}$-norm of the heat operator difference from the first sentence of this proof, which is controlled further by \eqref{eq:lp81}. Let us now condition on the event in which $\|\grad_{\r}^{\mathbf{T}}\mathbf{Z}^{N}\|_{\t_{\mathrm{st}};\Z}\lesssim N^{-1/7}\|\mathbf{Z}^{N}\|_{\t_{\mathrm{st}};\Z} $ given $\r\in[0,\t^{\mathrm{av}}]$; note that $\t^{\mathrm{av}}\leq N^{-5/3}$ by construction in the statement of Lemma \ref{lemma:lp8} if $\beta_{\star}$ is small enough, and thus with Lemma \ref{lemma:ste}, we may condition on said time-gradient estimate for $\mathbf{Z}^{N}$ if we give up an event of probability at most a uniformly bounded constant times $N^{-100}$. To take advantage of this, we first observe the heat operator inside the norm in \eqref{eq:lp81} is controlled, in absolute value, by $\Phi(1)+\Phi(2)$, where
\begin{align}
\Phi(1) \ \overset{\bullet}= \ |\grad_{{-\mathrm{r}}}^{\mathbf{T}}\mathbf{H}_{\t,\x}^{\beta_{\star},\mathbf{T}}(\mathbf{1}_{\s\leq\t^{\star}}c_{N}\mathbf{1}_{\y=0}\mathbf{Z}^{N}_{\s,\y})| \and \Phi(2) \ \overset{\bullet}= \ \mathbf{H}_{\t,\x}^{\beta_{\star},\mathbf{T}}\left(c_{N}\mathbf{1}_{\y=0}|\grad_{-\r}^{\mathbf{T}}(\mathbf{1}_{\s\leq\t^{\star}}\mathbf{Z}_{\s,\y}^{N})|\right).
\end{align}
Indeed, the bound in terms of $\Phi(1)+\Phi(2)$ follows by a version of integrating by parts in the time-integral of the heat operator. We note that in $\Phi(1)$, the time gradient acts on $\mathbf{H}^{\beta_{\star},\mathbf{T}}$ in the time-variable $\t$. In $\Phi(2)$, the time-gradient acts on the product of the indicator function and Gartner transform in the $\s$-variable. We also clarify that the scale in the time-gradients in $\Phi(1)$ and $\Phi(2)$ has an additional sign; this parallels the fact that in integration-by-parts, a negative sign is picked up. Lastly, we note that the boundary terms from integration-by-parts are encoded in the time-gradient acting on the heat operator in $\Phi(1)$; recalling this heat operator depends on $\t$ both on the upper limit for time-integration and on the heat kernel forwards time-variable, letting the time-gradient act on the former dependence on $\t$ gives the boundary term. We will now estimate $\Phi(1)$ and $\Phi(2)$. Let us start with $\Phi(1)$. For this, we recall $c_{N}\lesssim N$ and then employ the regularity estimates in Proposition \ref{prop:hke3} and Proposition \ref{prop:hke8}. To make this precise, we again note that the time-gradient acting on the heat operator in $\Phi(1)$ is bounded above by a boundary term given by the integral of the heat kernel at forwards spatial-variable at $\y=0$ on a time-interval of length $\r$ plus the integral of the time-gradient of the heat kernel. In particular, we obtain the following for $\Phi(1)$ in which we are loose with dropping indicator functions for the sake of upper bounds; see the proof of time-regularity in Proposition 3.2 of \cite{DT} for details concerning time-gradients of space-time heat operators:
\begin{align}
\Phi(1) \ &\lesssim \ c_{N}\int_{\t-\r}^{\t}\mathbf{H}_{\s,\t,\x,0}^{\beta_{\star}}\mathbf{Z}^{N}_{\s,0}\d\s + c_{N}\int_{0}^{\t-\r}|\grad_{-\r}^{\mathbf{T}}\mathbf{H}_{\s,\t,\x,0}^{\beta_{\star}}|\mathbf{Z}^{N}_{\s,0}\d\s\ \\
&\lesssim \ N^{-1+\frac52\beta_{\star}}c_{N}|\r|^{1/2}\|\mathbf{Z}^{N}\|_{\t_{\mathrm{st}};\Z}  + N^{-2+99999\beta_{\star}}c_{N}\|\mathbf{Z}^{N}\|_{\t_{\mathrm{st}};\Z} + N^{-1+5\beta_{\star}}c_{N}|\r|^{1/2}\|\mathbf{Z}^{N}\|_{\t_{\mathrm{st}};\Z} , \label{eq:lp82}
\end{align}
where the second bound in \eqref{eq:lp82} follows by plugging in heat kernel estimates in Proposition \ref{prop:hke3} and Proposition \ref{prop:hke8} and straightforward integration. After power-counting with the observation $|\r|\leq\t^{\mathrm{av}}=N^{-2+100\beta_{\star}}$, the bound of \eqref{eq:lp82} is clearly controlled by $N^{-\gamma}\|\mathbf{Z}^{N}\|_{\t_{\mathrm{st}};\Z} $ with $\gamma$ positive and depending only on $\beta_{\star}$, given $\beta_{\star}$ is sufficiently small. We now bound $\Phi(2)$. To this end, we may control the time-gradient of the product between the indicator function and $\mathbf{Z}^{N}$ by the time-gradient of the indicator function times $\|\mathbf{Z}^{N}\|_{\t;\Z}$ given that we only look at $\mathbf{Z}^{N}$ at $\y=0$ plus the time-gradient of $\mathbf{Z}^{N}$; this is a discrete version of the Leibniz rule. In any case, the time-gradient of the indicator function is supported in a length $\r$-neighborhood in time of $\t^{\star}$, so when integrating the heat kernel on this time-scale, we get another term to control that is given by integrating the heat kernel on a length-$\r$ time-interval as in \eqref{eq:lp82}. Summarizing this brief discussion, we ultimately obtain the following estimate for $\Phi(2)$:
\begin{align}
\Phi(2) \ &\lesssim \ c_{N}\int_{\t^{\star}-\r}^{\t^{\star}+\r}\mathbf{1}_{\s\leq\t}\mathbf{H}_{\s,\t,\x,\y}^{\beta_{\star}}\d\s\|\mathbf{Z}^{N}\|_{\t_{\mathrm{st}};\Z} + c_{N}\int_{0}^{\t}\mathbf{H}_{\s,\t,\x,\y}^{\beta_{\star}}|\grad_{-\r}^{\mathbf{T}}\mathbf{Z}_{\s,\y}^{N}|\d\s \\
&\lesssim \ N^{-1+\frac52\beta_{\star}}c_{N}|\r|^{1/2}\|\mathbf{Z}^{N}\|_{\t_{\mathrm{st}};\Z} + N^{-1+\frac52\beta_{\star}}c_{N}N^{-1/7}\|\mathbf{Z}^{N}\|_{\t_{\mathrm{st}};\Z}. \label{eq:lp83}
\end{align}
The estimate \eqref{eq:lp83} follows by heat kernel estimates in Proposition \ref{prop:hke3} along with the time-gradient estimate we have conditioned on that is written explicitly after \eqref{eq:lp81}. Again, the bound \eqref{eq:lp83} is controlled by $N^{-\gamma}\|\mathbf{Z}^{N}\|_{\t_{\mathrm{st}};\Z}$ for $\gamma$ positive and depending only on $\beta_{\star}$ assuming that $\beta_{\star}$ is small enough. This shows the term inside the norm in \eqref{eq:lp81} is controlled by $N^{-\gamma}\|\mathbf{Z}^{N}\|_{\t_{\mathrm{st}};\Z}$ for such a $\gamma$ exponent. Because this is uniformly on $(\t,\x)\in[0,\t_{\mathrm{st}}]\times\Z$ and $\r\in[0,\t^{\mathrm{av}}]$, this shows \eqref{eq:lp81} is controlled by $N^{-\gamma}\|\mathbf{Z}^{N}\|_{\t_{\mathrm{st}};\Z}$ on the high probability event on which we have the time-gradient estimate for $\mathbf{Z}^{N}$ written after \eqref{eq:lp81}, so we are done.
\end{proof}
%%%
%%%
\begin{proof}[Proof of \emph{Lemma \ref{lemma:lp9}}]
Let us unfold the heat operator on the LHS of \eqref{eq:lp9a}. After doing so, we employ the heat kernel estimates in Proposition \ref{prop:hke3}. Replacing $\mathbf{Z}^{N}$ with its maximal value over the space-time set $[0,\t_{\mathrm{st}}]\times\{0\}$, which is controlled by $\|\mathbf{Z}^{N}\|_{\t_{\mathrm{st}};\Z}$, we obtain the following deterministic estimate in which we control the short-time singularity in the heat kernel estimate of Proposition \ref{prop:hke3} via $|\t-\s|\geq N^{-5\beta_{\star}-\e}$ for $\s\leq\t^{\star}$, which implies, for the same $\s$, that $|\t-\s|^{-1/2}\leq N^{5\beta_{\star}/2+\e/2}$:
\begin{align}
|\mathbf{H}^{\beta_{\star}}_{\t,\x}(\mathbf{1}_{\s\leq\t^{\star}}2^{-1}c_{N}\mathbf{1}_{\y=0}\mathsf{I}^{\mathbf{T},\t^{\mathrm{av}}}(\bar{\mathfrak{q}}_{\s,\y})\mathbf{Z}_{\s,\y}^{N})| \ &\lesssim \ c_{N}\|\mathbf{Z}^{N}\|_{\t_{\mathrm{st}};\Z}\int_{0}^{\t^{\star}}\mathbf{H}_{\s,\t,\x,0}^{\beta_{\star}}|\mathsf{I}^{\mathbf{T},\t^{\mathrm{av}}}(\bar{\mathfrak{q}}_{\s,0})|\d\s \\
&\lesssim \ N^{-1+5\beta_{\star}+\e/2}c_{N}\|\mathbf{Z}^{N}\|_{\t_{\mathrm{st}};\Z}\int_{0}^{1}|\mathsf{I}^{\mathbf{T},\t^{\mathrm{av}}}(\bar{\mathfrak{q}}_{\s,0})|\d\s. \label{eq:lp91}
\end{align}
We clarify that on the RHS of \eqref{eq:lp91}, we should actually have an upper limit of integration of $\t^{\star}$, but because $\t^{\star}\leq\t\leq\t_{\mathrm{st}}\leq1$ with probability 1 and because the integrand is non-negative, we may extend the integration to $[0,1]$. Recalling $c_{N}\lesssim N$, dividing by $\|\mathbf{Z}^{N}\|_{\t_{\mathrm{st}};\Z}$, and taking expectation, the desired estimate \eqref{eq:lp9a} follows.
\end{proof}
%%%
%%%
\begin{proof}[Proof of \emph{Lemma \ref{lemma:lp10}}]
By the Fubini theorem, we first deduce the identity
\begin{align}
\E\int_{0}^{1}|\mathsf{I}^{\mathbf{T},\t^{\mathrm{av}}}(\bar{\mathfrak{q}}_{\s,0})|\d\s \ = \ \int_{0}^{1}\E|\mathsf{I}^{\mathbf{T},\t^{\mathrm{av}}}(\bar{\mathfrak{q}}_{\s,0})|\d\s. \label{eq:lp101}
\end{align}
We now decompose the expectation on the RHS of \eqref{eq:lp101} as follows. We write said expectation as an expectation with respect to the time-independent path-space measure induced by the particle system, where time-independence comes from time-homogeneity of the particle system dynamic, with an initial configuration given by the time-$\s$ configuration in the ``actual particle system" that we are realizing. We then take an expectation of this path-space expectation over said initial configuration according to the time-$\s$ law of the particle system, which we denote by $\mathsf{Q}(\s)$. Thus, we have the following that uses notation explained afterwards:
\begin{align}
\int_{0}^{1}\E|\mathsf{I}^{\mathbf{T},\t^{\mathrm{av}}}(\bar{\mathfrak{q}}_{\s,0})|\d\s \ = \ \int_{0}^{1}\E^{\mu_{0,\Z}}\mathsf{Q}(\s)\E^{\mathrm{path}}|\mathsf{I}^{\mathbf{T},\t^{\mathrm{av}}}(\bar{\mathfrak{q}}_{0,0})|\d\s; \label{eq:lp102}
\end{align}
the path-space expectation $\E^{\mathrm{path}}$ is with respect to the path-space law of the particle system on $\Z$. We now employ another Fubini theorem to move the integral with respect to $\s$ inside the expectation. Because $\E^{\mathrm{path}}$ is time-independent, we can pull $\E^{\mathrm{path}}$ outside this $\d\s$-integral. By definition of the time-averaged law $\bar{\mathsf{Q}}(1)$ in Definition \ref{definition:lp6}, we deduce
\begin{align}
\int_{0}^{1}\E^{\mu_{0,\Z}}\mathsf{Q}(\s)\E^{\mathrm{path}}|\mathsf{I}^{\mathbf{T},\t^{\mathrm{av}}}(\bar{\mathfrak{q}}_{0,0})|\d\s \ = \ \E^{\mu_{0,\Z}}\bar{\mathsf{Q}}(1)\E^{\mathrm{path}}|\mathsf{I}^{\mathbf{T},\t^{\mathrm{av}}}(\bar{\mathfrak{q}}_{0,0})|. \label{eq:lp103}
\end{align}
If we could now replace $\E^{\mathrm{path}}$ by $\E^{\mathrm{path},\t^{\mathrm{av}}}$, then the proof would be complete; the latter of these path-space expectations is defined in the statement of Lemma \ref{lemma:lp10} as an expectation with respect to path-space measure induced by a ``periodization" of the particle system dynamic on $\mathbb{I}(\t^{\mathrm{av}})\subseteq\Z$, which is also defined in the statement of Lemma \ref{lemma:lp10}. We will make the replacement of $\E^{\mathrm{path}}$ by $\E^{\mathrm{path},\t^{\mathrm{av}}}$ up to a cost of order at most $N^{-100}$, thereby completing the proof of Lemma \ref{lemma:lp10}. For this, we observe that it suffices to find a coupling of path-space measures defining $\E^{\mathrm{path}}$ and $\E^{\mathrm{path},\t^{\mathrm{av}}}$, whose initial configurations on $\mathbb{I}(\t^{\mathrm{av}})$ agree, such that the $\eta$-values on the support of $\bar{\mathfrak{q}}_{0,0}$ of Definition \ref{definition:lp2} agree for times in $[0,\t^{\mathrm{av}}]$. Indeed, this would imply the existence of a coupling of $\E^{\mathrm{path}}$ and $\E^{\mathrm{path},\t^{\mathrm{av}}}$ path-space measures for which the value of $\bar{\mathfrak{q}}_{\r,0}$ for $\r\in[0,\t^{\mathrm{av}}]$, and thus the average of $\bar{\mathfrak{q}}$ on this time-interval, is the same for both processes with high probability. This is done similar to the proof of Lemma \ref{lemma:lp7} with minor adjustments below.
%%%
\begin{itemize}
\item Let $\mathsf{S}(\Z)$ denote the realization of the $\E^{\mathrm{path}}$ process on $\Omega_{\Z}$ to be used in the following construction, and let $\mathsf{S}(\t^{\mathrm{av}})$ denote the realization of the $\E^{\mathrm{path},\t^{\mathrm{av}}}$ process on $\Omega_{\mathbb{I}(\t^{\mathrm{av}})}$, which is defined in the statement of Lemma \ref{lemma:lp10}, to be used below.
\item Assume that the $\eta$-values of $\mathsf{S}(\Z)$ and $\mathsf{S}(\t^{\mathrm{av}})$ on $\mathbb{I}(\t^{\mathrm{av}})$ are initially all equal. To construct a coupling the dynamics of $\mathsf{S}(\Z)$ and $\mathsf{S}(\t^{\mathrm{av}})$, for any bond that is present in both $\Z$ and $\mathbb{I}(\t^{\mathrm{av}})$ realized as a torus with periodic boundary conditions, we first couple the associated ``spin-swap" dynamics in the symmetric parts of the $\mathsf{S}(\Z)$ and $\mathsf{S}(\t^{\mathrm{av}})$ dynamics. For any shared bond, we also employ the basic coupling for the totally asymmetric part of the $\mathsf{S}(\Z)$ and $\mathsf{S}(\t^{\mathrm{av}})$ processes. So far, this is basically the construction in the proof of Lemma \ref{lemma:lp7} for shared bonds in $\Z$ and $\mathbb{I}(\t^{\mathrm{av}})$. For bonds that are not shared, for example those induced by periodic boundary conditions in $\mathsf{S}(\t^{\mathrm{av}})$, we will let the $\mathsf{S}(\Z)$ and $\mathsf{S}(\t^{\mathrm{av}})$ processes act independently, though this is not too important.
\item Let us emphasize that discrepancies can only be created by bonds that are not shared between $\mathsf{S}(\Z)$ and $\mathsf{S}(\t^{\mathrm{av}})$; see the proof of Lemma \ref{lemma:lp7}. Thus, they must arise from Poisson clocks at the boundary of $\mathbb{I}(\t^{\mathrm{av}})$, in particular at bonds that comprise of a point inside $\mathbb{I}(\t^{\mathrm{av}})$ and either a point outside $\mathbb{I}(\t^{\mathrm{av}})$ or at the boundary of $\mathbb{I}(\t^{\mathrm{av}})$ because of the periodic boundary conditions in $\mathsf{S}(\t^{\mathrm{av}})$.
\item Because the particles can only perform jumps of length 1, we observe that any discrepancy between $\mathsf{S}(\Z)$ and $\mathsf{S}(\t^{\mathrm{av}})$ inside the subset $\mathbb{I}(\t^{\mathrm{av}})$ must appear in a radius-1 neighborhood of the boundary of $\mathbb{I}(\t^{\mathrm{av}})$. In particular, any discrepancy between $\mathsf{S}(\Z)$ and $\mathsf{S}(\t^{\mathrm{av}})$ in the support of $\bar{\mathfrak{q}}_{0,0}$ must originally be a discrepancy that showed up within distance one of the boundary of $\mathbb{I}(\t^{\mathrm{av}})$. Moreover, because these discrepancies that are created necessarily within distance one of the boundary of $\mathbb{I}(\t^{\mathrm{av}})$ are created by one of at most 10 Poisson clocks, where 10 is just a fixed finite number, and because these clocks all have speed $N^{2}$, until time $\t^{\mathrm{av}}\leq1$ there can be at most $N^{100}$ discrepancies outside an event of exponentially low probability in $N$, because there can be at most $N^{100}$ ringings of these 10 Poisson clocks until time 1 except on an exponentially-low-probability event. 
\item Let us additionally note that if a discrepancy propagates from near the boundary of $\mathbb{I}(\t^{\mathrm{av}})$ to the support of $\bar{\mathfrak{q}}$, which according to construction in Definition \ref{definition:lp2} is contained in the radius-$N^{20\beta_{\star}}$ ball centered at the origin, then this discrepancy must have traveled a distance $N(\t^{\mathrm{av}})^{1/2}\log^{50}N+N^{3/2}\t^{\mathrm{av}}\log^{50}N$ through bonds that are shared between $\Z$ and $\mathbb{I}(\t^{\mathrm{av}})$; we emphasize that according to construction of $\mathbb{I}(\t^{\mathrm{av}})$ in the statement of Lemma \ref{lemma:lp10}, the proposed traveled distance just mentioned is much smaller than the length of $\mathbb{I}(\t^{\mathrm{av}})$ even without the last $N^{20\beta_{\star}}$ in its length. Again, this last claim about how the discrepancy must travel in to the support of $\bar{\mathfrak{q}}$ follows from the length-1 jumps of the particles in the particle systems. However, similar to the proof of Lemma \ref{lemma:lp7}, the law of this discrepancy path along shared bonds between $\Z$ and $\mathbb{I}(\t^{\mathrm{av}})$ is a symmetric random walk of speed $N^{2}$ plus an asymmetric jump/killing mechanism of speed $N^{3/2}$. The extra log-powers in $N$ in the necessarily traveled distance $N(\t^{\mathrm{av}})^{1/2}\log^{50}N+N^{3/2}\t^{\mathrm{av}}\log^{50}N$ guarantee that this happens with probability at most $\exp(-\kappa\log^{50}N)$ for some strictly positive $\kappa$ by standard large deviation estimates for random walks. We then note $\exp(-\kappa\log^{50}N)\lesssim N^{-300}$. 
\end{itemize}
%%%
The previous bullet point provides a coupling of $\E^{\mathrm{path}}$ and $\E^{\mathrm{path},\t^{\mathrm{av}}}$ path-space measures under which outside an event of exponentially small probability in $N$, there are at most $N^{100}$ candidate discrepancies in $\eta$-values in the support of $\bar{\mathfrak{q}}_{0,0}$, where each of these candidates is actually realized as such a discrepancy before time $\t^{\mathrm{av}}$ with probability at most of order $N^{-300}$. Taking a union bound over these at most $N^{100}$-many discrepancies, we may find a coupling of $\E^{\mathrm{path}}$ and $\E^{\mathrm{path},\t^{\mathrm{av}}}$ path-space measures for which $\bar{\mathfrak{q}}_{\r,0}$-values agree for all $\r\in[0,\t^{\mathrm{av}}]$ outside an event with probability at most order $N^{-200}$. Therefore, the cost of replacing $\E^{\mathrm{path}}$ by $\E^{\mathrm{path},\t^{\mathrm{av}}}$ on the RHS of \eqref{eq:lp103} is at most order $N^{-200}$ times the worst-case upper bound for the time-average of $\bar{\mathfrak{q}}$ in \eqref{eq:lp103}, the last of which is uniformly bounded. Combining this cost of at most order $N^{-200}$ in the replacement of $\E^{\mathrm{path}}$ by $\E^{\mathrm{path},\t^{\mathrm{av}}}$ in \eqref{eq:lp103} with \eqref{eq:lp101} and \eqref{eq:lp102}, this completes the proof of the proposed inequality.
\end{proof}
%%%
%%%
\begin{proof}[Proof of \emph{Lemma \ref{lemma:lp11}}]
Per $\sigma\in\R$, we may apply the Kipnis-Varadhan inequality in Appendix 1.6 of \cite{KL} to obtain the inequality below, in which $\|\|_{-1}$ refers to a negative-Sobolev norm that we clarify afterwards:
\begin{align}
\E^{\sigma,\mathbb{I}}\E^{\mathrm{path},\t^{\mathrm{av}}}|\mathsf{I}^{\mathbf{T},\mathfrak{t}^{\mathrm{av}}}(\bar{\mathfrak{q}}_{0,0})|^{2} \ \lesssim \ (\t^{\mathrm{av}})^{-1}\|\bar{\mathfrak{q}}\|_{-1}^{2}. \label{eq:lp111}
\end{align}
We now note that the vanishing of canonical measure expectation $\E^{\sigma,\mathbb{I}}\bar{\mathfrak{q}}_{0,0}=0$ given any $\sigma\in\R$, which follows from the fact that $\eta$-variables are exchangeable under canonical measures since exchanging $\eta$-values preserves particle numbers and thus defines a bijection on the support of any canonical measure, which keeps uniform canonical measures invariant. In this case, we may appeal to Proposition 3.3 in \cite{GJ15}, which deals with grand-canonical measures but this is not a big deal because canonical measures span grand-canonical measures on $\Omega_{\mathbb{I}}$ for $\mathbb{I}\subseteq\Z$ finite, to deduce the following, in which we use uniform boundedness of $\bar{\mathfrak{q}}$:
\begin{align}
\|\bar{\mathfrak{q}}\|_{-1}^{2} \ \lesssim \ N^{-2}{N^{40\beta_{\star}}}. \label{eq:lp112}
\end{align}
We clarify that the factor {$N^{40\beta_{\star}}$} comes from the fact that the support of $\bar{\mathfrak{q}}$ has length $N^{20\beta_{\star}}$, and Proposition 3.3 in \cite{GJ15} provides an estimate for $\|\|_{-1}^{2}$ that is quadratic in the support length. We also clarify the extra factor $N^{-2}$ comes from the $N^{2}$-speed of the interacting particle system; the relevance of this speed can be interpreted as either viewing $\|\|_{-1}$-norm as a $2$-norm of the inverse-square-root of the negative of the generator acting on the functional at hand, or it can be equivalently seen via the reparameterization calculation given briefly after Proposition 3.1 in \cite{DT}. In any case, we plug \eqref{eq:lp112} into the RHS of \eqref{eq:lp111} to get a $\sigma$-uniform estimate for the LHS of \eqref{eq:lp111}. We then get the desired estimate by controlling a first-moment by square root of the second moment.
\end{proof}
%%%
%%%
\begin{proof}[Proof of \emph{Lemma \ref{lemma:lp12}}]
We first prove \eqref{eq:lp12a} with a straightforward heat kernel estimate. First, we may forget the $\mathbf{Z}^{N}$ factor on the LHS of \eqref{eq:lp12a} at the cost of the supremum of $\mathbf{Z}^{N}$ on the space-time set $[0,\t_{\mathrm{st}}]\times\{0\}$. Controlling this supremum by $\|\mathbf{Z}^{N}\|_{\t_{\mathrm{st}};\Z}$,
\begin{align}
N^{1/2}|\mathbf{H}_{\t,\x}^{\beta_{\star}}(\mathbf{1}_{\y=0}\mathbf{Z}_{\s,\y}^{N})| \ \leq \ N^{1/2}\|\mathbf{Z}^{N}\|_{\t_{\mathrm{st}};\Z}\int_{0}^{\t}\mathbf{H}_{\s,\t,\x,0}^{\beta_{\star}}\d\s \ \lesssim \ N^{-1/2+5\beta_{\star}/2}\|\mathbf{Z}^{N}\|_{\t_{\mathrm{st}};\Z}, \label{eq:lp121}
\end{align}
where the last estimate in \eqref{eq:lp121} follows by estimating the heat kernel via Proposition \ref{prop:hke3} and straightforward integration. The far RHS of \eqref{eq:lp121} controls the far LHS of \eqref{eq:lp121} uniformly in $[0,\t_{\mathrm{st}}]\times\Z$, so we deduce from this \eqref{eq:lp12a}. The proof of \eqref{eq:lp12b} follows from the proof of Lemma \ref{lemma:lp9} without taking expectations and replacing $\bar{\mathfrak{q}}$ by $\mathsf{I}^{\mathbf{X},\mathfrak{l}_{N}}(\eta_{\s,\y+1})$.
\end{proof}
%%%
%%%
\begin{proof}[Proof of \emph{Lemma \ref{lemma:lp13}}]
By definition of the height function as $N^{-1/2}$ times an integrated version of $\eta$-variables, we first observe
\begin{align}
|\mathsf{I}^{\mathbf{X},\mathfrak{l}_{N}}(\eta_{\s,1})| \ = \ |N^{1/2}\mathfrak{l}_{N}^{-1}\grad_{\mathfrak{l}_{N}}^{\mathbf{X}}\mathbf{h}_{\s,1}^{N}| \ = \ |N^{1/2}\mathfrak{l}_{N}^{-1}\grad_{\mathfrak{l}_{N}}^{\mathbf{X}}\log\mathbf{G}_{\s,1}^{N}| \ = \ N^{1/2}\mathfrak{l}_{N}^{-1}|\log\left(1+(\mathbf{G}_{\s,1}^{N})^{-1}\grad_{\mathfrak{l}_{N}}^{\mathbf{X}}\mathbf{G}_{\s,1}^{N}\right)|, \label{eq:lp131}
\end{align}
where the second identity in \eqref{eq:lp131} follows by the definition of the Gartner transform and the final identity in \eqref{eq:lp131} follows from basic properties of the logarithm; we note the spatially constant terms in $\mathbf{h}^{N}$ and $\mathbf{G}^{N}$ are irrelevant after taking spatial gradients. Now note $|(\mathbf{G}_{\s,1}^{N})^{-1}\grad_{\mathfrak{l}_{N}}^{\mathbf{X}}\mathbf{G}_{\s,1}^{N}|\lesssim N^{-1/2}|\mathfrak{l}_{N}|\ll1$ deterministically by Taylor expansion. Indeed, by definition of $\mathbf{G}^{N}$ and $\mathbf{h}^{N}$,
\begin{align}
|\grad_{\mathfrak{l}_{N}}^{\mathbf{X}}\mathbf{G}_{\s,1}^{N}| \ = \ |\exp(-\mathbf{h}_{\s,1+\mathfrak{l}_{N}}^{N}+\mathrm{R}_{N}\s)-\exp(-\mathbf{h}_{\s,1}^{N}+\mathrm{R}_{N}\s)| \ &= \ |\exp(-\grad_{\mathfrak{l}_{N}}^{\mathbf{X}}\mathbf{h}_{\s,1}^{N})-1|\exp(-\mathbf{h}_{\s,1}^{N}+\mathrm{R}_{N}\s) \\
&\leq \ \mathbf{G}_{\s,1}^{N}{\sum}_{\mathfrak{j}\geq1}(\mathfrak{j}!)^{-1}|\grad_{\mathfrak{l}_{N}}^{\mathbf{X}}\mathbf{h}_{\s,1}^{N}|^{\mathfrak{j}} \\
&\lesssim \ N^{-1/2}|\mathfrak{l}_{N}|\mathbf{G}_{\s,1}^{N},
\end{align}
where the final bound follows by noting that the $\grad^{\mathbf{X}}\mathbf{h}^{N}$-term above is equal to $N^{-1/2}$ times a sum of $\mathfrak{l}_{N}$-many uniformly bounded $\eta$-variables. This allows us to use the inequality $|\log(1+\x)|\lesssim|\x|$ for $|\x|\leq2^{-1}$, which follows by control of the gradient of $\log$ for such $\x$. In particular, this and \eqref{eq:lp131} give the following deterministic estimate:
\begin{align}
|\mathsf{I}^{\mathbf{X},\mathfrak{l}_{N}}(\eta_{\s,1})| \ \lesssim \ N^{1/2}\mathfrak{l}_{N}^{-1}|(\mathbf{G}_{\s,1}^{N})^{-1}\grad_{\mathfrak{l}_{N}}^{\mathbf{X}}\mathbf{G}_{\s,1}^{N}|. \label{eq:lp132}
\end{align}
Multiplying both sides of \eqref{eq:lp132} by $\mathbf{G}_{\s,1}^{N}$ completes the proof if we swap the gradient of $\mathbf{G}^{N}$ by that of $\mathbf{Z}^{N}$, which we can do with the required high probability by Lemma \ref{lemma:mshe5} up to a cost of order $N^{-99}$; this follows as $\mathbf{G}^{N}$ and $\mathbf{Z}^{N}$ are evaluated at points within $|\mathfrak{l}_{N}|\ll N$ of the origin, and Lemma \ref{lemma:mshe5} compares $\mathbf{G}^{N}$ and $\mathbf{Z}^{N}$ uniformly over macroscopic length-scales.
\end{proof}
%%%
%%%
\begin{proof}[Proof of \emph{Lemma \ref{lemma:lp14}}]
We will control each term on the RHS of the stochastic integral equation for $\mathbf{Z}^{N}$ from Definition \ref{definition:mshe4}. First, for any $\t\in[0,\t_{\mathrm{st}}]$, we consider the following identity consequence of linearity of $\grad^{\mathbf{X}}$ gradient operators:
\begin{align}
\grad_{\mathfrak{l}_{N}}^{\mathbf{X}}\mathbf{Z}_{\t,\x}^{N} \ = \ \grad_{\mathfrak{l}_{N}}^{\mathbf{X}}\mathbf{H}_{\t,\x}^{\beta_{\star},\mathbf{X}}(\mathbf{Z}^{N})+\grad_{\mathfrak{l}_{N}}^{\mathbf{X}}\mathbf{H}_{\t,\x}^{\beta_{\star},\mathbf{T}}(\mathbf{Z}^{N}\d\xi^{N})+\grad_{\mathfrak{l}_{N}}^{\mathbf{X}}\mathbf{H}_{\t,\x}^{\beta_{\star},\mathbf{T}}(2^{-1}c_{N}\mathbf{1}_{\y=0}\mathfrak{q}_{\s,\y}^{\mathrm{tot}}\mathbf{Z}_{\s,\y}^{N}). \label{eq:lp141}
\end{align}
Let us now estimate each term on the RHS of \eqref{eq:lp141}, beginning with the first and third terms therein, as these terms require only a deterministic analysis based on heat kernel estimates in Proposition \ref{prop:hke12} and otherwise elementary calculations.
%%%
\begin{itemize}
\item Unfolding the definition of the spatial heat operator, the first term on the RHS of \eqref{eq:lp141} is actually the spatial-sum action of the spatial gradient of the heat kernel on the $\mathbf{Z}^{N}$ factor. More precisely, we have the estimate
\begin{align}
|\grad_{\mathfrak{l}_{N}}^{\mathbf{X}}\mathbf{H}_{\t,\x}^{\beta_{\star},\mathbf{X}}(\mathbf{Z}^{N})| \ \leq \ {\sum}_{\y\in\Z}|\grad_{\mathfrak{l}_{N}}^{\mathbf{X}}\mathbf{H}_{0,\t,\x,\y}^{\beta_{\star}}|\mathbf{Z}^{N}_{0,\y} \ \leq \ \|\mathbf{Z}^{N}\|_{\t_{\mathrm{st}};\Z}{\sum}_{\y\in\Z}|\grad_{\mathfrak{l}_{N}}^{\mathbf{X}}\mathbf{H}_{0,\t,\x,\y}^{\beta_{\star}}|\exp(N^{-1}|\y|). \label{eq:lp142}
\end{align}
Now, multiplying by $\exp(-N^{-1}|\x|)$ and employing the tail estimates in Proposition \ref{prop:hke12} for the heat kernel as in the proof of Proposition 3.2 in \cite{DT}, we deduce the following deterministic estimate. Roughly speaking, when we sum here the tail bound in Proposition \ref{prop:hke12} for the spatial gradient of the heat kernel, the sub-exponential factor takes away a factor of $N^{-1}|\t-\s|^{-1/2}$ by a straightforward geometric series or integral comparison calculation; this was the case in the proof of Proposition 3.2 in \cite{DT} as well. We get the following in which $\beta_{\star}\lesssim\e_{i}\lesssim\beta_{\star}$, and the extension to $\|\|_{\t_{\mathrm{st}};\Z}$-bounds follows as the RHS of \eqref{eq:lp143} is constant:
\begin{align}
\|\t^{1/2}\grad_{\mathfrak{l}_{N}}^{\mathbf{X}}\mathbf{H}_{\t,\x}^{\beta_{\star},\mathbf{X}}(\mathbf{Z}^{N})\|_{\t_{\mathrm{st}};\Z}(1+\|\mathbf{Z}^{N}\|_{\t_{\mathrm{st}};\Z}^{2})^{-1} \ \lesssim \ N^{-1+2\beta_{\star}}|\mathfrak{l}_{N}|+N^{-1+99999\beta_{\star}+\e_{2}}|\mathfrak{l}_{N}|. \label{eq:lp143}
\end{align}
Let us briefly remark that employing Proposition \ref{prop:hke12} is technically illegal if $|\mathfrak{l}_{N}|> N\t^{1/2}\vee1$, which is required for Proposition \ref{prop:hke12}, is violated. Moreover, even if we could apply Proposition \ref{prop:hke12}, the second term on the RHS of \eqref{eq:lp143} would not have the $|\mathfrak{l}_{N}|$-factor. However, we can always control a length-$\mathfrak{l}_{N}$ gradient in terms of $|\mathfrak{l}_{N}|$-many length-1 gradients, for which we get \eqref{eq:lp143} for $\mathfrak{l}_{N}$ replaced by $1$, and then we multiply these bounds by $|\mathfrak{l}_{N}|$ to account for all of these length-1 gradients, which is why both terms on the RHS of \eqref{eq:lp143} admit $|\mathfrak{l}_{N}|$ factors. Ultimately, we deduce that the first term on the RHS of \eqref{eq:lp141} satisfies the proposed estimate for $\grad_{\mathfrak{l}_{N}}^{\mathbf{X}}\mathbf{Z}^{N}$ deterministically.
\item We now move to the third term on the RHS of \eqref{eq:lp141} and we employ a similar idea. Unfolding similar to the derivation of \eqref{eq:lp142} but for the space-time heat operator, we have the following deterministic inequality; recall $|c_{N}|\lesssim N$ from Proposition \ref{prop:mshe1}:
\begin{align}
|\grad_{\mathfrak{l}_{N}}^{\mathbf{X}}\mathbf{H}_{\t,\x}^{\beta_{\star},\mathbf{T}}(2^{-1}c_{N}\mathbf{1}_{\y=0}\mathfrak{q}_{\s,\y}^{\mathrm{tot}}\mathbf{Z}_{\s,\y}^{N})| \ \lesssim \ N\|\mathbf{Z}^{N}_{\t,\x}\mathbf{1}_{\x=0}\|_{\t_{\mathrm{st}};\Z;\max}\int_{0}^{\t}|\grad_{\mathfrak{l}_{N}}^{\mathbf{X}}\mathbf{H}_{\s,\t,\x,0}^{\beta_{\star}}|\d\s. \label{eq:lp144}
\end{align}
The norm on the RHS of \eqref{eq:lp144} is controlled by forgetting the $\max$-subscript because it localizes to $\x=0$. Now, plugging the on-diagonal heat kernel estimate of Proposition \ref{prop:hke8} into the RHS of \eqref{eq:lp144}, we deduce the RHS of \eqref{eq:lp144} is bounded above by the quantity below; we note both terms below have a factor of $|\mathfrak{l}_{N}|$ similar to the RHS of \eqref{eq:lp143}:
\begin{align}
N^{-1+2\beta_{\star}}|\mathfrak{l}_{N}|\|\mathbf{Z}^{N}\|_{\t_{\mathrm{st}};\Z}\int_{0}^{\t}|\t-\s|^{-1+\e_{1}}\d\s + N^{-1+99999\beta_{\star}+\e_{2}}|\mathfrak{l}_{N}|\|\mathbf{Z}^{N}\|_{\t_{\mathrm{st}};\Z}\int_{0}^{\t}|\t-\s|^{-1+\e_{2}}\d\s. \label{eq:lp145}
\end{align}
Integrating \eqref{eq:lp145} and observing the resulting bound is uniform in space-time, we obtain the estimate
\begin{align}
\|\grad_{\mathfrak{l}_{N}}^{\mathbf{X}}\mathbf{H}_{\t,\x}^{\beta_{\star},\mathbf{T}}(2^{-1}c_{N}\mathbf{1}_{\y=0}\mathfrak{q}_{\s,\y}^{\mathrm{tot}}\mathbf{Z}_{\s,\y}^{N})\|_{\t_{\mathrm{st}};\Z;\max}(1+\|\mathbf{Z}^{N}\|_{\t_{\mathrm{st}};\Z}^{2})^{-1} \ \lesssim \ N^{-1+2\beta_{\star}}|\mathfrak{l}_{N}| + N^{-1+99999\beta_{\star}+\e_{2}}|\mathfrak{l}_{N}|. \label{eq:lp146}
\end{align}
The above estimate \eqref{eq:lp146} implies that the third term on the RHS of \eqref{eq:lp141} also satisfies the required estimate for $\grad_{\mathfrak{l}_{N}}^{\mathbf{X}}\mathbf{Z}^{N}$ in the statement of Lemma \ref{lemma:lp14} deterministically, similar to the first term on the RHS of \eqref{eq:lp141} as noted in the previous bullet point.
\end{itemize}
%%%
Provided the two bullet points above, it suffices to prove that the second term on the RHS of \eqref{eq:lp141} satisfies the required expectation estimate for $\grad_{\mathfrak{l}_{N}}^{\mathbf{X}}\mathbf{Z}^{N}$. To this end, we must employ a somewhat more delicate argument than for the other terms in \eqref{eq:lp141} and even compared to the proof of Proposition 3.2 in \cite{DT} due to the $\|\|_{\t_{\mathrm{st}};\Z}$-norm we must control the second term on the RHS of \eqref{eq:lp141} by; this norm does not easily combine well with the martingale inequalities in Lemma \ref{lemma:mge}, for example. We will start with the following decomposition, where $\Phi$ is the second term on the RHS of \eqref{eq:lp141} at $\x=0$ and $\mathcal{E}(\mathfrak{k})$ is the event in which {$\mathfrak{k}\leq1+\|\mathbf{Z}^{N}\|_{\t_{\mathrm{st}};\Z}\leq\mathfrak{k}+1$}:
\begin{align}
\Phi^{2} \ = \ \left({\sum}_{\mathfrak{k}\geq1}\Phi\mathbf{1}(\mathcal{E}(\mathfrak{k}))\right)^{2} \ \lesssim \ {\sum}_{\mathfrak{k}\geq1}\Phi^{2}\mathbf{1}(\mathcal{E}(\mathfrak{k})). \label{eq:lp147}
\end{align}
The last inequality in \eqref{eq:lp147} follows from expanding the square in the middle of \eqref{eq:lp147} and the Schwarz inequality upon noting that $\mathbf{1}(\mathcal{E}(\mathfrak{k}))\mathbf{1}(\mathcal{E}(\mathfrak{l}))=0$ unless $|\mathfrak{l}-\mathfrak{k}|\leq1$. Multiplying both the far LHS and far RHS of \eqref{eq:lp147} by $(1+\|\mathbf{Z}^{N}\|_{\t_{\mathrm{st}};\Z}^{2})^{-2}$ and taking expectations, we obtain the following inequality in which we set $\mathbf{Z}=1+\|\mathbf{Z}^{N}\|_{\t_{\mathrm{st}};\Z}^{2}$ for notational convenience:
\begin{align}
\E\Phi^{2}\mathbf{Z}^{-2} \ \lesssim \ {\sum}_{\mathfrak{k}\geq1}\E\Phi^{2}\mathbf{Z}^{-2}\mathbf{1}(\mathcal{E}(\mathfrak{k})). \label{eq:lp148}
\end{align}
We will now estimate each summand on the RHS of \eqref{eq:lp148} in a summable fashion. First, we observe that on the event $\mathcal{E}(\mathfrak{k})$, we have $\mathbf{Z}^{-2}={(1+\|\mathbf{Z}^{N}\|_{\t_{\mathrm{st}};\Z}^{2})^{-2}\lesssim(1+\|\mathbf{Z}^{N}\|_{\t_{\mathrm{st}};\Z})^{-4}\lesssim \mathfrak{k}^{-4}}$. Moreover, on the event $\mathcal{E}(\mathfrak{k})$, we may rewrite $\Phi$ as follows:
\begin{align}
\Phi\mathbf{1}(\mathcal{E}(\mathfrak{k})) \ = \ \mathbf{1}(\mathcal{E}(\mathfrak{k}))\grad_{\mathfrak{l}_{N}}^{\mathbf{X}}\mathbf{H}_{\t,0}^{\beta_{\star},\mathbf{T}}(\mathbf{Z}^{N,\mathfrak{k}}\d\xi^{N}) \quad \mathrm{where} \quad \mathbf{Z}^{N,\mathfrak{k}}_{\s,\y} \ = \ \mathbf{Z}^{N}_{\s,\y}\mathbf{1}(|\mathbf{Z}^{N}_{\s,\y}|\lesssim \exp(N^{-1}|\y|){\mathfrak{k}}).
\end{align}
Indeed, on the event $\mathcal{E}(\mathfrak{k})$, the indicator function defining $\mathbf{Z}^{N,\mathfrak{k}}$ is redundant, but we emphasize that $\mathbf{Z}^{N,\mathfrak{k}}$ is still an adapted process. Combining \eqref{eq:lp148} with the observations made afterwards until now, we deduce the following estimate where we ultimately remove the $\mathbf{1}(\mathcal{E}(\mathfrak{k}))$ indicator functions for the sake of obtaining an upper bound and where we estimate the sum over $\mathfrak{k}\geq1$ directly:
\begin{align}
\E\Phi^{2}\mathbf{Z}^{-2} \ \lesssim \ \left({\sum}_{\mathfrak{k}\geq1}{\mathfrak{k}^{-2}}\right){\sup}_{\mathfrak{k}\geq1}\E|\grad_{\mathfrak{l}_{N}}^{\mathbf{X}}\mathbf{H}_{\t,0}^{\beta_{\star},\mathbf{T}}({\mathfrak{k}^{-1}}\mathbf{Z}^{N,\mathfrak{k}}\d\xi^{N})|^{2} \ \lesssim \ {\sup}_{\mathfrak{k}\geq1}\E|\grad_{\mathfrak{l}_{N}}^{\mathbf{X}}\mathbf{H}_{\t,0}^{\beta_{\star},\mathbf{T}}({\mathfrak{k}^{-1}}\mathbf{Z}^{N,\mathfrak{k}}\d\xi^{N})|^{2}. \label{eq:lp149}
\end{align}
To estimate these expectations on the far RHS of \eqref{eq:lp149}, we employ the approach of spatial regularity estimates for the stochastic integral in Proposition 3.2 of \cite{DT}. However, in this case we are not actually integrating the product of the original Gartner transform $\mathbf{G}^{N}$ against $\d\xi^{N}$ but instead ${\mathfrak{k}^{-1}}\mathbf{Z}^{N,\mathfrak{k}}$. Fortunately, the term ${\mathfrak{k}^{-1}}\mathbf{Z}^{N,\mathfrak{k}}$ is uniformly bounded by $\exp(N^{-1}|\y|)$ by construction. Therefore, we may follow the proof of Lemma 3.1 in \cite{DT}, which estimates heat operators acting against products of adapted process with martingale differentials, but with ${\mathfrak{k}^{-1}}\mathbf{Z}^{N,\mathfrak{k}}$ in place of the original Gartner transform and then replacing all ${\mathfrak{k}^{-1}}\mathbf{Z}^{N,\mathfrak{k}}$ factors with $\exp(N^{-1}|\y|)$ factors for the sake of an upper bound. Ultimately, we get
\begin{align}
\E|\grad_{\mathfrak{l}_{N}}^{\mathbf{X}}\mathbf{H}_{\t,0}^{\beta_{\star},\mathbf{T}}({\mathfrak{k}^{-1}}\mathbf{Z}^{N,\mathfrak{k}}\d\xi^{N})|^{2} \ \lesssim \ N\int_{0}^{\t}{\sum}_{\y\in\Z}|\grad_{\mathfrak{l}_{N}}^{\mathbf{X}}\mathbf{H}_{\s,\t,0,\y}^{\beta_{\star}}|^{2}\exp(2N^{-1}|\y|)\d\s. \label{eq:lp1410}
\end{align}
At this point, we employ the on-diagonal spatial gradient bound in Proposition \ref{prop:hke8} and the fact that the heat kernel is a probability measure on $\Z$ in the forwards spatial variable with sub-exponential tails of any weight/speed/power. At the end of the day, we have the following estimates in which $k_{N}=N^{-2+2\beta_{\star}}|\mathfrak{l}_{N}|+N^{-2+999999\beta_{\star}}|\mathfrak{l}_{N}|$:
\begin{align}
{\sum}_{\y\in\Z}|\grad_{\mathfrak{l}_{N}}^{\mathbf{X}}\mathbf{H}_{\s,\t,0,\y}^{\beta_{\star}}|^{2}\exp(2N^{-1}|\y|) \ &\lesssim \ k_{N}|\t-\s|^{-1+\e_{1}}{\sum}_{\y\in\Z}|\grad_{\mathfrak{l}_{N}}\mathbf{H}_{\s,\t,0,\y}^{\beta_{\star}}|\exp(2N^{-1}|\y|) \nonumber \\
&\lesssim \ k_{N}|\t-\s|^{-1+\e_{1}}{\sum}_{\mathfrak{j}=0,\mathfrak{l}_{N}}{\sum}_{\y\in\Z}\mathbf{H}_{\s,\t,\mathfrak{j},\y}^{\beta_{\star}}\exp(2N^{-1}|\y|) \nonumber \\
&\lesssim \ k_{N}|\t-\s|^{-1+\e_{1}}{\sum}_{\mathfrak{j}=0,\mathfrak{l}_{N}}{\sum}_{\y\in\Z}\mathbf{H}_{\s,\t,\mathfrak{j},\y}^{\beta_{\star}}\exp(2N^{-1}(|\y-\mathfrak{j}|)) \nonumber \\
&\lesssim \ k_{N}|\t-\s|^{-1+\e_{1}}. \nonumber
\end{align}
The second inequality follows by controlling the gradient of the heat kernel by the heat kernel at two points that differ in spatial coordinate by $\mathfrak{l}_{N}$. The third inequality follows by the triangle inequality in the exponential and recalling $|\mathfrak{j}|\leq|\mathfrak{l}_{N}|=N^{20\beta_{\star}}\ll N$, so the multiplicative error/cost in this triangle inequality is uniformly bounded. The last inequality follows from the off-diagonal estimate for the heat kernel $\mathbf{H}^{\beta_{\star}}$ we explained in the proof of Proposition \ref{prop:hke11}, namely that its tails are sub-exponential with any speed and therefore admits a Laplace transform of any exponent/parameter. Integrating the last bound in $\s\in[0,\t]$ and multiplying by $N$, we deduce via  \eqref{eq:lp149} and \eqref{eq:lp1410} the bound $\E\Phi^{2}\mathbf{Z}^{-2}\lesssim Nk_{N}$. This is certainly controlled by $N^{-1-10\beta_{\star}-\gamma}|\mathfrak{l}_{N}|^{2}$ after a brief calculation upon recalling $\mathfrak{l}_{N}=N^{20\beta_{\star}}$, assuming $\beta_{\star}$ is sufficiently small. This completes the proof of the lemma upon recalling that $\Phi$ is the second term on the RHS of \eqref{eq:lp141} evaluated at $\x=0$ and $\mathbf{Z}=1+\|\mathbf{Z}^{N}\|_{\t_{\mathrm{st}};\Z}^{2}$.
\end{proof}
%%%
\appendix
%%%
\section{Martingale Estimate}
%%%
The main result of this section is a technical martingale inequality, which generalizes the BDG inequality of Lemma 3.1 in \cite{DT}. Roughly speaking, the proof of Lemma 3.1 in \cite{DT} is exclusive to space-time stochastic-type integrals of $\mathbf{G}^{N}\d\xi^{N}$; the fact that we consider the Gartner transform factor $\mathbf{G}^{N}$ is crucial in the proof of Lemma 3.1 in \cite{DT}. This is because short-time behaviors of $\mathbf{G}^{N}$ on microscopic time-scales can be controlled explicitly in terms of the particle system. To generalize Lemma 3.1 in \cite{DT} to arbitrary products of $\d\xi^{N}$ and adapted processes, we need control on these processes on short time-scales. We ultimately get the following.
%%%
\begin{lemma}\label{lemma:mge}
\fsp Consider any $\phi:\R_{\geq0}\times\Z\to\R$ and its nearest-neighbor square defined below given fixed times $0\leq\mathfrak{t}_{1}\leq\mathfrak{t}_{2}$; in the following, we define $\lfloor\r\rfloor_{N}$ to be the largest element in $N^{-2}\Z_{\geq0}$ that is bounded above by $\r$:
\begin{align}
\wt{\phi}_{\r,\x,\w}^{\mathfrak{t}_{1},\mathfrak{t}_{2}} \ &\overset{\bullet}= \ {\sup}_{\mathfrak{r}'\in[\mathfrak{t}_{1},\mathfrak{t}_{2}): \ \lfloor\mathfrak{r}'\rfloor_{N}=\lfloor \r\rfloor_{N}}{\sup}_{|\mathfrak{j}|\leq1} |\phi_{\mathfrak{r}',\x,\w+\mathfrak{j}}\phi_{\mathfrak{r}',\x,\w}|.
\end{align}
We additionally consider any process $\mathbf{X}^{N}$ on $\R_{\geq0}\times\Z$ satisfying the following stochastic equation, where $\mathrm{V}=\{\mathrm{V}_{1},\mathrm{V}_{2},\mathrm{V}_{3}\}$ with $\mathrm{V}_{i}$ functionals on $\R_{\geq0}\times\Z\times\Omega_{\Z}$, and the operator $\mathscr{L}_{N,\mathrm{V}}$ is defined via the second equation below for fixed $\mathfrak{l}\in\Z$:
\begin{align}
\d\mathbf{X}_{\t,\x}^{N} \ &= \ \mathscr{L}_{N,\mathrm{V}}\mathbf{X}_{\t,\x}^{N}\d\t + \mathrm{V}_{3;\t,\x}\d\t + \mathbf{X}_{\t,\x}^{N}\d\xi_{\t,\x}^{N} \\
&= \ 2^{-1}N^{2}\Delta_{\x}\mathbf{X}_{\t,\x}^{N}\d\t+\mathrm{V}_{3;\t,\x}\d\t+\mathbf{X}_{\t,\x}^{N}\d\xi_{\t,\x}^{N} + \mathrm{V}_{1;\t,\x}\mathbf{X}_{\t,\x}^{N}\d\t + \grad_{\mathfrak{l}}^{\mathbf{X}}(\mathrm{V}_{2;\t,\x}\mathbf{X}_{\t,\x}^{N})\d\t.
\end{align}
Suppose we have $|\mathrm{V}_{1}|+|\mathrm{V}_{2}|+|\mathrm{V}_{3}|\leq N^{3/2}$ uniformly in space-time and the randomness of the particle system. Provided any $p\in\R_{\geq1}$ and any $0\leq\mathfrak{t}_{1}\leq\mathfrak{t}_{2}$ and $\phi:\R_{\geq0}\times\Z\to\R$, all of which are deterministic, we have the following where $\|\cdot\|_{\omega;2p}^{2p}=\E|\cdot|^{2p}$:
\begin{align}
\|\int_{\mathfrak{t}_{1}}^{\mathfrak{t}_{2}}\sum_{\w\in\Z}\phi_{\r,\x,\w}\mathbf{X}_{\r,\w}^{N}\d\xi_{\r,\w}^{N}\|_{\omega;2p}^{2} \ \lesssim_{p} \ \int_{\mathfrak{t}_{1}}^{\mathfrak{t}_{2}}N\left(\sup_{\w\in\Z}\|\mathbf{X}_{\lfloor\r\rfloor_{N},\w}^{N}\|_{\omega;2p}^{2}\right)\sum_{\w\in\Z}\wt{\phi}_{\r,\x,\w}^{\mathfrak{t}_{1},\mathfrak{t}_{2}}\d\r. \label{eq:mg}
\end{align}
In particular, \eqref{eq:mg} holds for $\mathbf{X}^{N}$ equal to the process $\mathbf{Z}^{N}$ from \emph{Definition \ref{definition:mshe4}}, the processes $\mathbf{Y}^{N}$ and $\mathbf{W}^{N}$ from \emph{Definition \ref{definition:pc1}}, and the processes $\mathbf{D}^{N,1}$ and $\mathbf{D}^{N,1,1}$ from \emph{Definition \ref{definition:pc7}}.
\end{lemma}
%%%
%%%
\begin{remark}\label{remark:mge}
\fsp In the estimate \eqref{eq:mg}, we can include a factor of $\exp(-2\kappa N^{-1}|\w|)$ inside the supremum on the RHS for any $\kappa\geq0$ if we include $\exp(2\kappa N^{-1}|\w|)$ inside the sum over $\w\in\Z$ on the RHS; see Lemma 3.1 and the proof of Proposition 3.2 in \cite{DT}.
\end{remark}
%%%
The proof of Lemma \ref{lemma:mge} is not too much different than the proof of its analog of Lemma 3.1 in \cite{DT} beyond a few observations. Before we proceed with the proof, let us first introduce the following useful notation for a \emph{deterministic} kernel.
%%%
\begin{definition}\label{definition:mge2}
\fsp Let $\mathbf{J}$ with variables in $\R_{\geq0}^{2}\times\Z^{2}$ solve the following equation with initial data at $\t=\s$ given by $\mathbf{J}_{\s,\s,\x,\y}=\mathbf{1}_{x=y}$:
\begin{align}
{\partial_{\t}}\mathbf{J}_{\s,\t,\x,\y} \ = \ 2^{-1}N^{2}\Delta_{\x}\mathbf{J}_{\s,\t,\x,\y} + 99N^{3/2}{\sum}_{|\mathfrak{j}|\leq|\mathfrak{l}|}\mathbf{J}_{\s,\t,\x+\mathfrak{j},\y} + 99N^{3/2}.
\end{align}
It is easy to see that $\mathbf{J}$ is space-time-homogeneous because its coefficients are independent of space-time. Moreover, we require only standard PDE procedure to see that $\mathbf{J}$ satisfies the following estimate for any $\kappa\geq0$, in which $\mathbf{J}_{\w,\z}=\sup_{0\leq\t\leq N^{-2}}\mathbf{J}_{0,\t,\w,\z}$:
\begin{align}
{\sup}_{\w\in\Z}{\sum}_{\z\in\Z}\exp(\kappa|\w-\z|)|\mathbf{J}_{\w,\z}| \ \lesssim_{\kappa} \ 1.
\end{align}
Indeed, the operator defining the RHS of the $\mathbf{J}$ equation, dropping the constant term therein, relates the growth of $\mathbf{J}$ to its nearest-neighbors, which gives the sub-exponential tail behavior following the same Gronwall inequality procedure in the proof of Proposition 3.1 of \cite{DT} without the stochastic term, and this operator has norm uniformly bounded by $N^{2}$ when viewed as an operator from $\mathscr{L}^{\infty}(\Z)$ to itself. Thus, the growth of $\mathbf{J}$ is uniformly controlled on time-scales $\t\leq N^{-2}$ relevant for $\mathbf{J}_{\w,\z}$.
\end{definition}
%%%
%%%
\begin{proof}
The fact that we can take the proposed choices of $\mathbf{X}^{N}$ in the last sentence of the statement of Lemma \ref{lemma:mge} comes from noting each of these proposed choices satisfies a stochastic equation of the type assumed for $\mathbf{X}^{N}$ in Lemma \ref{lemma:mge}. Thus, we move to proving the estimate \eqref{eq:mg}. To this end, we inspect the proof of Lemma 3.1 in \cite{DT}, beginning with (3.7) and (3.8) therein. In particular, we first claim that it suffices to prove the following, in which $\mathfrak{n}_{\r,\mathfrak{t}}$ denotes the total number of times that a Poisson clock at $\w\in\Z$ rings in the time-interval $\lfloor\r\rfloor_{N}+[0,\mathfrak{t}]$, in which the time $\t\geq0$ satisfies $\mathfrak{t}\leq N^{-2}$, in which the constant $\kappa\geq0$ is fixed, and in which $\mathscr{F}_{\cdot}$ is the canonical filtration associated to the interacting particle system that we condition on below:
\begin{align}
\E_{\mathscr{F}_{\lfloor\r\rfloor_{N}}}\left(\mathfrak{n}_{\r,\mathfrak{t}}^{2p}|\mathbf{X}_{\lfloor\r\rfloor_{N}+\mathfrak{t},\w}^{N}|^{2p}\right) \ \lesssim_{p} \ {\sum}_{\z\in\Z}\exp(\kappa|\w-\z|)\mathbf{J}_{\w,\z}|\mathbf{X}_{\lfloor\r\rfloor_{N},\z}^{N}|^{2p}. \label{eq:mg1}
\end{align}
Provided \eqref{eq:mg1}, we may now take expectations of both sides over the randomness in the $\lfloor\r\rfloor_{N}$-algebra. Because $\mathbf{J}$ is deterministic, we can push said expectation past everything on the RHS of \eqref{eq:mg1} and have it only act on the $2p$-moment of $\mathbf{X}^{N}$. We then get
\begin{align}
\E\E_{\mathscr{F}_{\lfloor\r\rfloor_{N}}}\left(\mathfrak{n}_{\r,\mathfrak{t}}^{2p}|\mathbf{X}_{\lfloor\r\rfloor_{N}+\mathfrak{t},\w}^{N}|^{2p}\right) \ \lesssim_{p} \ {\sum}_{\z\in\Z}\exp(\kappa|\w-\z|)\mathbf{J}_{\w,\z}\E|\mathbf{X}_{\lfloor\r\rfloor_{N},\z}^{N}|^{2p} \ \lesssim \ {\sup}_{\w\in\Z}\E|\mathbf{X}_{\lfloor\r\rfloor_{N},\z}^{N}|^{2p}. \label{eq:mg2}
\end{align}
The last estimate in \eqref{eq:mg2} follows from the estimate for $\mathbf{J}$ noted and explained in Definition \ref{definition:mge2}. At this point, we may then follow the proof of Lemma 3.1 in \cite{DT} after (3.8) therein but without the adjustment of the time-variable $\lfloor\r\rfloor_{N}$ therein; this completes the proof modulo the proof of the estimate \eqref{eq:mg1}. To obtain \eqref{eq:mg1}, let us note that if $\mathbf{X}^{N}$ did not have any stochastic $\mathbf{X}^{N}\d\xi^{N}$ term in its defining equation that we assume for it, then \eqref{eq:mg1} would follow by noting that $\mathbf{X}^{N}$ evolves according to an equation whose fundamental solution is controlled by $\mathbf{J}_{\w,\z}$ along with a moment bound for the uniformly-bounded-rate Poisson clock counter $\mathfrak{n}$ that was used for (3.8) in \cite{DT}. We clarify that this Poisson random variable has uniformly bounded rate, because it is determined by a speed $N^{2}$ Poisson process on a time-interval of length $N^{-2}$. To deal with the $\mathbf{X}^{N}\d\xi^{N}$-term, we may treat $\d\xi^{N}$ as a multiplicative potential. Precisely, the fundamental solution for the equation defining $\mathbf{X}^{N}$ is controlled by $\mathbf{J}_{\w,\z}$ times the exponential of something that counts, up to some factors of $N^{-1/2}$ that come in front of $\d\xi^{N}$ to account for the size of the jumps in the Gartner transform, the number of jumps in the particle system at the moving location of a random walk trajectory whose transition probability is bounded by the $\mathbf{J}$ kernel. This is the Feynman-Kac formula, and it can be checked readily. In any case, we are then left to control from above the same moment as in the LHS of (3.8) in \cite{DT} but without the $Z$-factor therein and the exponential therein is replaced by itself, viewed as a function of the spatial variable $\x'$ convolved with the $\mathbf{J}_{\w,\z}$ kernel. This last replacement for the exponential comes from the fact that the aforementioned Feynman-Kac interpretation of $\mathbf{X}^{N}$ looks at all jumps in the particle system weighted by probabilities controlled by $\mathbf{J}_{\w,\z}$, not just jumps at a single point. Indeed, the advantage of the Gartner transform is explicit control on $\mathbf{G}^{N}$ in terms of jumps at a single point coming via the explicit formula for the Gartner transform in terms of the particle system, and it is why the RHS in \eqref{eq:mg} appears different than the RHS of the estimate in Lemma 3.1 in \cite{DT}. In any case, sub-exponential bounds for $\mathbf{J}_{\w,\z}$ noted in Definition \ref{definition:mge2} and elementary Poisson estimates like (3.8) in \cite{DT} complete the proof of \eqref{eq:mg1}.
\end{proof}
%%%
%
%
%
%%%
\section{Short-Time and Localization Estimates}
%%%
This section provides important short-time behavior estimates for solutions to general stochastic equations. By short-time, we mean sub-microscopic time-scales on which we expect nothing to happen. In particular, the following result is not saying too much, but it is a technically useful reduction to analysis on semi-discrete sets to analysis on fully discrete sets. The proof of the following set of estimates will not be sophisticated, so we avoid giving too many details.
%%%
\begin{lemma}\label{lemma:ste}
\fsp Consider any process $\mathbf{X}^{N}$ on $\R_{\geq0}\times\Z$ satisfying the following stochastic equation, in which $\mathrm{V}=\{\mathrm{V}_{1},\mathrm{V}_{2}\}$ where $\mathrm{V}_{1}$ and $\mathrm{V}_{2}$ are functionals on $\R_{\geq0}\times\Z\times\Omega_{\Z}$, and the operator $\mathscr{L}_{N,\mathrm{V}}$ is defined via the second equation below for $\mathfrak{l}\in\Z$:
\begin{align}
\d\mathbf{X}_{\t,\x}^{N} \ = \ \mathscr{L}_{N,\mathrm{V}}\mathbf{X}_{\t,\x}^{N}\d\t + \mathbf{X}_{\t,\x}^{N}\d\xi_{\t,\x}^{N} \ = \ 2^{-1}N^{2}\Delta_{\x}\mathbf{X}_{\t,\x}^{N}+\mathbf{X}_{\t,\x}^{N}\d\xi_{\t,\x}^{N} + \mathrm{V}_{1;\t,\x}\mathbf{X}_{\t,\x}^{N}\d\t + \grad_{\mathfrak{l}}^{\mathbf{X}}(\mathrm{V}_{2;\t,\x}\mathbf{X}_{\t,\x}^{N})\d\t.
\end{align}
Suppose that $|\mathrm{V}_{1}|+|\mathrm{V}_{2}|\lesssim N^{2}$ uniformly in all variables, and suppose that $|\mathfrak{l}|\lesssim1$. If $\mathbf{X}^{N}\not\equiv0$, we have the following estimate outside an event of probability $N^{-100}$, in which $\e$ is an arbitrarily small but fixed positive constant:
\begin{align}
\sup_{\mathfrak{s}\in[0,N^{-99}]}\|\mathbf{X}^{N}\|_{\mathfrak{t};\Z}^{-1}\|\grad_{\mathfrak{s}}^{\mathbf{T}}\mathbf{X}^{N}\|_{\t;\Z} \ \leq \ N^{-1/2+\e}. \label{eq:st2}
\end{align}
We also have the following estimate on the same event uniformly in $\t\in[0,1]$, in which $\mathbb{I}^{\mathbf{T}}=\{\mathfrak{j}N^{-100}\}_{\mathfrak{j}=0,\ldots,N^{100}}$:
\begin{align}
\|\mathbf{X}^{N}\|_{\t;\Z} \ \lesssim \ \sup_{\t\in\mathbb{I}^{\mathbf{T}}}\sup_{\x\in\Z}\exp(-|\x|/N)|\mathbf{X}^{N}_{\t,\x}|. \label{eq:st3}
\end{align}
Finally, both estimates \eqref{eq:st2} and \eqref{eq:st3} hold for $\mathbf{X}^{N}$ equal to the process $\mathbf{Z}^{N}$ from \emph{Definition \ref{definition:mshe4}}, the processes $\mathbf{Y}^{N}$ and $\mathbf{W}^{N}$ from \emph{Definition \ref{definition:pc1}}, and the processes $\mathbf{D}^{N,1}$ and $\mathbf{D}^{N,1,1}$ from \emph{Definition \ref{definition:pc7}}. Actually, for the choice of $\mathbf{X}^{N}=\mathbf{Z}^{N}$ and $\mathbf{X}^{N}=\mathbf{G}^{N}$, we also have the following longer-time-scale estimate with the same probability for any possibly random time $\t_{\mathrm{st}}\geq0$:
\begin{align}
\sup_{\s\in[0,N^{-5/3}]}\|\grad_{\s}^{\mathbf{T}}\mathbf{X}^{N}_{\t,\x}\mathbf{1}_{\x=0}\|_{\t_{\mathrm{st}};\Z} \ \lesssim \ N^{-1/6+\e}\|\mathbf{X}^{N}_{\t,\x}\mathbf{1}_{\x=0}\|_{\t_{\mathrm{st}};\Z} + N^{-100}.
\end{align}
\end{lemma}
%%%
%%%
\begin{proof}
The fact that \eqref{eq:st2} and \eqref{eq:st3} hold for the proposed choices of $\mathbf{X}^{N}$ in the last sentence in the statement of Lemma \ref{lemma:ste} follow as in the beginning of the proof of Lemma \ref{lemma:mge}. Let us first focus on proving \eqref{eq:st2}. To this end, we observe the following operator-norm estimate for $\mathscr{L}_{N,\mathrm{V}}$, whose validity we explain afterwards; for this, $\phi:\Z\to\R$ is a generic function, and $\|\phi\|_{\Z}$ denotes the supremum of $\exp(-|\x|/N)|\phi_{\x}|$ over $\x\in\Z$, namely the norm $\|\|_{\t;\Z}$ but for time-independent functions:
\begin{align}
\|\mathscr{L}_{N,\mathrm{V}}\phi\|_{\Z} \ \lesssim \ N^{2}\|\phi\|_{\Z}. \label{eq:st21}
\end{align}
Indeed, the estimate \eqref{eq:st21} would hold if the $\|\|_{\Z}$-norms did not include the sub-exponential factor because $\mathscr{L}_{N,\mathrm{V}}\phi$ controls $\phi$ in terms of $\phi$ and its values at neighboring points and then times a factor of order $N^{2}$. To account for the sub-exponential factor, it suffices to note that $\mathscr{L}_{N,\mathrm{V}}\phi$ controls $\phi$ only in terms of its values at neighboring points, and the sub-exponential factor in $\|\|_{\Z}$-norms is controlled by its values at neighboring points times uniformly bounded factors. Now, if the stochastic equation for $\mathbf{X}^{N}$ did not have the $\mathbf{X}^{N}\d\xi^{N}$ term, then \eqref{eq:st2} would follow from the representation of $\mathbf{X}^{N}$ at time $\s+\t$ as $\exp(\s\mathscr{L}_{N,\mathrm{V}})$ acting on $\mathbf{X}^{N}$ at time $\t$ and then Taylor expansion of this exponential operator for $\s\leq N^{-99}$, which is certainly enough to beat the $N^{2}$ factor in \eqref{eq:st21}. To account for the $\mathbf{X}^{N}\d\xi^{N}$ term, we simply note that this only contributes by multiplying $\mathbf{X}^{N}$ by $\exp(\kappa N^{-1/2})$ for some positive $\kappa$ for every jump in the $\d\xi^{N}$ term. It now suffices to note that outside an event of probability $N^{-100}$ times an $\e$-dependent factor, we only see $N^{\e}$-many jumps over all $\d\xi^{N}$ in a time-interval of length $N^{-99}$, because this total number of $\d\xi^{N}$-jumps is distributed according to a Poisson process of speed of order $N^{30}$, given by the number of non-zero $\d\xi^{N}$ terms, times the time-interval-length of $N^{-99}$. Indeed, we run the deterministic equation argument based on $\exp(\s\mathscr{L}_{N,\mathrm{V}})$ between jumps, and these jumps could only provide an additional factor of $\exp(\kappa N^{-1/2})$, and therefore a jump of order $N^{-1/2}$. This establishes \eqref{eq:st2}. To establish \eqref{eq:st3}, we condition on \eqref{eq:st2} and consider any $\t\in[0,1]$. Let $\t^{\mathbb{I}}$ denote the largest element in $\mathbb{I}^{\mathbf{T}}$ less than or equal to $\t$. The value of $\mathbf{X}^{N}$ at time $\t$ is controlled by the value of $\mathbf{X}^{N}$ at time $\t^{\mathbb{I}}$ plus a time-gradient of $\mathbf{X}^{N}$ with time-scale at most $N^{-100}$. This latter time-gradient is controlled, because of \eqref{eq:st2}, by the maximal value of $\mathbf{X}^{N}$ at time $\t^{\mathbb{I}}$ times a small factor $N^{-1/2+\e}$. This ultimately controls $\mathbf{X}^{N}$ at time $\t$ by the maximal value of $\mathbf{X}^{N}$ times the relevant sub-exponential factor over points in the space-time discrete set $\mathbb{I}^{\mathbb{T}}\times\Z$, which is exactly what is needed to prove \eqref{eq:st3}, so we are done for estimates \eqref{eq:st2} and \eqref{eq:st3} for the general choice of $\mathbf{X}^{N}$. 

Let us now specialize to $\mathbf{X}^{N}=\mathbf{Z}^{N}$. By Lemma \ref{lemma:mshe5}, it suffices to prove the estimate for $\mathbf{G}^{N}$ as $\mathbf{Z}^{N}$ and $\mathbf{G}^{N}$ differ by at most order $N^{-100}$ at $\x=0$ with the required probability in the statement of this lemma. In order to prove the estimate for $\mathbf{G}^{N}$, let us recall the construction of the Gartner transform $\mathbf{G}^{N}$ as a functional of the particle system. Based upon this representation for $\mathbf{G}^{N}$, it is easy to deduce that $\mathbf{G}^{N}$ evolves over a time-scale $\s$ by a factor of at most $\exp(\kappa N\s)$ times $\exp(\kappa N^{-1/2})$ for every jump at a single point in time-scale $\s$; in both factors, the constant $\kappa\geq0$ is fixed and uniformly bounded. For $|\s|\leq N^{-5/3}$, the first factor is easily controlled. Moreover, by the argument given for general $\mathbf{X}^{N}$ in the previous paragraph, with high probability we have only $N^{2}N^{-5/3+\e}=N^{1/3+\e}$ many jumps at a single point in $\Z$ outside an event of probability at most $N^{-100}$ times a uniformly bounded constant. Thus, if we multiply by at most $\exp(\kappa N^{-1/2})$ a total of $N^{1/3+\e}$-many times, we accumulate a total change of $\exp(\kappa N^{-1/2+1/3+\e})-1$, which is uniformly bounded by order $N^{-1/6+\e}$. This completes the proof.
\end{proof}
%%%
%
%
%
%%%
\section{Index for Notation}
%%%
%%%
\subsection{Expectation Operators}
%%%
We will let $\E^{\mu}$ denote expectation with respect to the probability measure $\mu$.
%%%
\subsection{Gradient Operators}
%%%
Consider any space-time function $\phi:\R_{\geq0}\times\Z\to\R$. Provided any $\mathfrak{k}\in\Z$, we define the length-scale $\mathfrak{k}$ spatial-gradient $\grad_{\mathfrak{k}}^{\mathbf{X}}\phi_{\t,\x}=\phi_{\t,\x+\mathfrak{k}}-\phi_{\t,\x}$. Provided any $\s\in\R$, define the time-scale $\s$ time-gradient $\grad_{\s}^{\mathbf{T}}\phi_{\t,\x}=\phi_{\t+\s,\x}-\phi_{\t,\x}$ if $\t+\s\geq0$, and otherwise $\grad_{\s}^{\mathbf{T}}\phi_{\t,\x}=\phi_{0,\x}-\phi_{\t,\x}$ if $\t+\s\leq0$. Indeed, we only consider $\phi$ at non-negative times.
%%%
\subsection{Averaging}
%%%
For any finite set $\mathscr{I}$ and any function $\psi:\mathscr{I}\to\R$, define the averaged sum $\wt{\sum}_{\x\in\mathscr{I}}\psi(\x)=|\mathscr{I}|^{-1}{\sum}_{\x\in\mathscr{I}}\psi(\x)$.
%%%
\subsection{Landau Notation}
%%%
For any generic set $\mathscr{I}$, the notation $\x\lesssim_{\mathscr{I}}\y$ is synonymous with $\x=\mathscr{O}(\y)$ where the implied constant is only allowed to depend on elements in $\mathscr{I}$. If we write $\x\lesssim\y$, the implied constant is assumed to be universal.
%%%
\subsection{Norms}
%%%
For any $p\geq1$ and any random variable $\mathrm{X}$, we define the $p$-th moment norm $\|\mathrm{X}\|_{\omega;p}^{p}=\E|\mathrm{X}|^{p}$. For any space-time function $\phi:\R_{\geq0}\times\Z\to\R$ along with any $\t\geq0$ and $\mathbb{K}\subseteq\R$, we also define $\|\phi\|_{\t;\mathbb{K};\max}={\sup}_{(\s,\y)\in[0,\t]\times\mathbb{K}}|\phi_{\s,\y}|$ along with the exponentially-weighted version $\|\phi\|_{\t;\mathbb{K}}={\sup}_{(\s,\y)\in[0,\t]\times\mathbb{K}}\exp(-N^{-1}|\y|)|\phi_{\s,\y}|$.
%%%
\subsection{Miscellaneous}
%%%
Provided any $\alpha_{1},\alpha_{2}\in\R$, we define the interval discretization $\llbracket\alpha_{1},\alpha_{2}\rrbracket=[\alpha_{1},\alpha_{2}]\cap\Z$.

%%%

%%%
\end{document}